\documentclass[12pt,letterpaper]{article}

\usepackage{graphics,epsfig,graphicx,color}
\graphicspath{{/EPSF/}{../Figures/}{Figures/}}
\usepackage{subfig}
\newsubfloat{figure}

\usepackage{amssymb,latexsym}

\usepackage{setspace}

\usepackage{amsmath}

\usepackage{mathrsfs,amsfonts}

\usepackage{amsthm,amsxtra}

\newif\ifPDF
\ifx\pdfoutput\undefined
	\PDFfalse
\else
	\ifnum\pdfoutput > 0
		\PDFtrue
	\else
		\PDFfalse
	\fi
\fi

\ifPDF
	\usepackage{pdftricks}
	\begin{psinputs}
		\usepackage{pstricks}
		\usepackage{pstcol}
		\usepackage{pst-plot}
		\usepackage{pst-tree}
		\usepackage{pst-eps}
		\usepackage{multido}
		\usepackage{pst-node}
		\usepackage{pst-eps}
	\end{psinputs}
\else
	\usepackage{pstricks}
\fi

\ifPDF
	\usepackage[debug,pdftex,colorlinks=true, 
	linkcolor=blue, bookmarksopen=false,
	plainpages=false,pdfpagelabels]{hyperref}
\else
	\usepackage[dvips]{hyperref}
\fi

\usepackage{fancyhdr}

\newcommand{\dsum}{\displaystyle\sum}

\newcommand{\bzero}{{\mathbf 0}}

\newcommand{\fT}{\mathfrak{T}}
\newcommand{\partite}{podal\ }

\newcommand{\bbR}{\mathbb R}

\newcommand{\cA}{\mathcal A} 
\newcommand{\cC}{\mathcal C}

\newcommand{\cS}{\mathcal S}

\newcommand{\E}{\epsilon} \newcommand{\T}{\tau}
\newcommand{\nd}{\noindent} \newcommand{\eps}{\eta}

\title{The Asymptotics of Large Constrained Graphs}

\author{
Charles Radin\thanks{Department of Mathematics, University of Texas, Austin, TX 78712; radin@math.utexas.edu}  
\and Kui Ren \thanks{Department of Mathematics, University of Texas, Austin, TX 78712; ren@math.utexas.edu} 
 \and Lorenzo Sadun\thanks{Department of Mathematics, University of Texas, Austin, TX 78712; sadun@math.utexas.edu} 
}

\begin{document}

\maketitle



\begin{abstract}
We show, through local estimates and simulation, that if one constrains simple
graphs by their densities $\E$ of edges and $\T$ of triangles, then
asymptotically (in the number of vertices) for over $95\%$ of the
possible range of those densities there is a well-defined typical
graph, and it has a very simple structure: the vertices are decomposed
into two subsets $V_1$ and $V_2$ of fixed relative size $c$  and $1-c$, and there
are well-defined probabilities of edges, $g_{jk}$, between  $v_j\in V_j$, and
$v_k\in V_k$. Furthermore the four parameters $c, g_{11},
g_{22}$ and $g_{12}$ are smooth functions of $(\E,\T)$ except at two smooth `phase
transition' curves.

\end{abstract}






\section{Introduction}
\label{SEC:Intro}

We consider large, simple graphs subject to two constraints: fixed
values of the density of edges, $\E$, and the density of triangles,
$\T$, where the densities are normalized so as to have value 1 for a
complete graph. Our goal is to obtain a qualitative understanding of
the asymptotic structure of such graphs, as the vertex number
diverges. (A graph is `simple' if the vertices are labelled, the
edges are undirected, there is at most one edge between any two
vertices, and there are no edges from a vertex to itself.)

We show that asymptotically there is a `unique typical graph'
characterized by a small number of parameters each of which is a
function of $\E$ and $\T$, and that the parameters vary smoothly in
$\E$ and $\T$ except across certain smooth (`phase transition')
curves. In particular we show that
more than $95\%$ of the `phase space' of possible pairs $(\E,\T)$ consists
of three phases separated by two smooth transition curves, and that within
these three regions the typical graph requires at most four parameters for
its description. Our strategy, and evidence, is a combination of local 
estimates and numerical simulation. The two parts are presented in separate sections but
they were intertwined in obtaining the results and we do not see how
to obtain them otherwise. In particular, in Section~\ref{SEC:numer} and the beginning of 
Section~\ref{SEC:Analysis} we present evidence (not proof) that typical graphs have a
very simple structure which we call multipodal. In the remainder of
Section~\ref{SEC:Analysis} we {\bf assume} this multipodal structure and derive two
things: the boundary between two of the resultant phases, and the behavior of
typical graphs near another phase boundary.

The precise definition of `unique typical graph' requires some technical 
clarification, given below, for two reasons:
to allow for the degeneracy associated with relabelling of vertices,
and to use a probabilistic description of graphs appropriate for
asymptotics. (This formalism of constrained graphs is natural if one is
motivated to obtain the sort of emergent behavior one sees in
thermodynamics; the constrained graphs are analogous to the
configurations in a microcanonical ensemble of a mean-field
statistical mechanical system, with $\E$ and $\T$ corresponding to mass
and energy densities \cite{RS1, RS2}.) 
We begin with a review of what is known about our constraint problem
from extremal graph theory. Important recent results are in \cite{Ra,
  PR}, which are also useful for background on earlier results.

\begin{figure}[!ht]
\centering
\includegraphics[angle=0,width=0.56\textwidth]{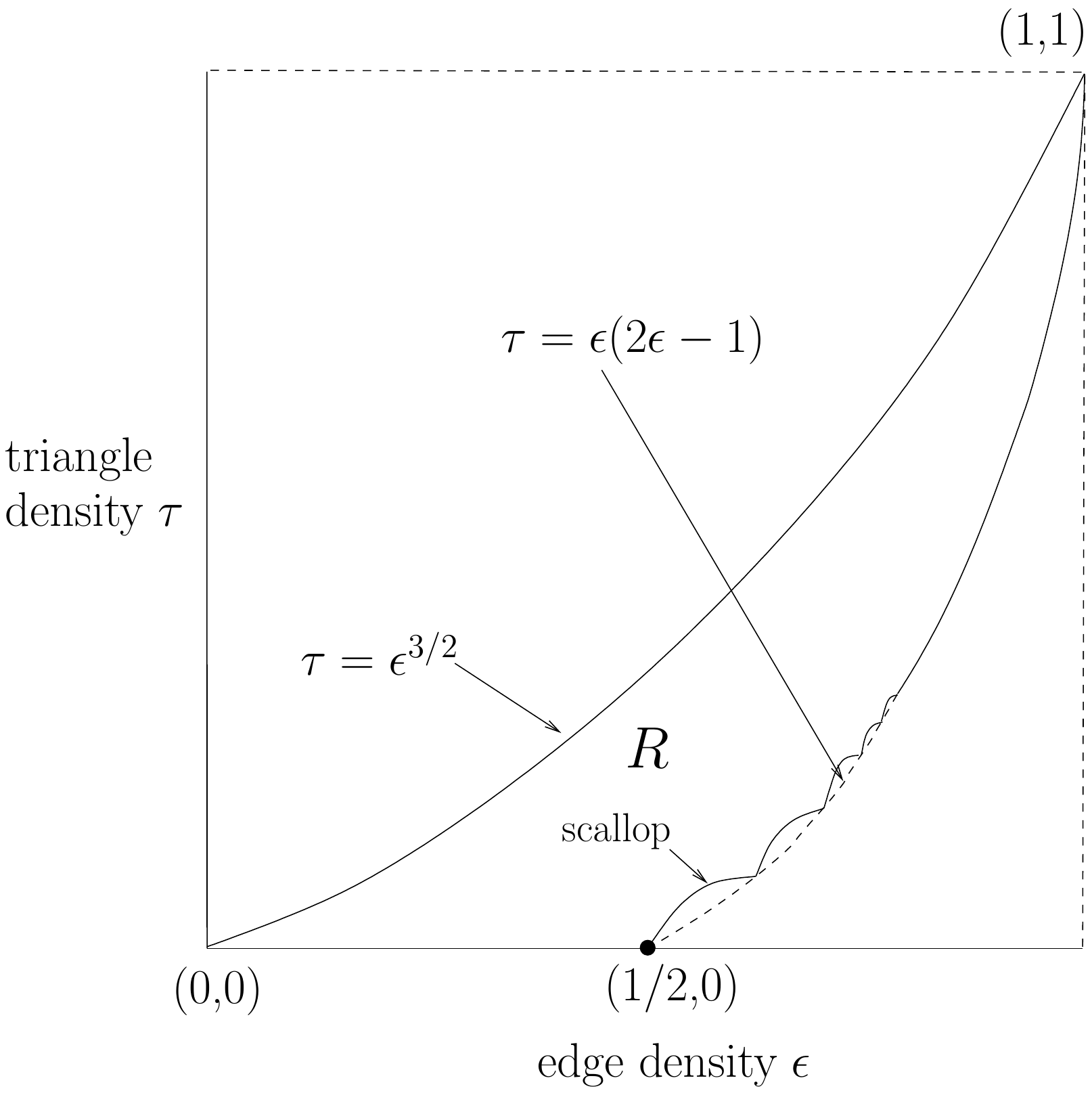}
\caption{The phase space $R$, outlined in solid lines. Features are exaggerated for clarity.}
\label{FIG:phasefig}
\end{figure}
Within the parameter space
$\{(\E, \T)\}\subset [0,1]^2$, for $1/2\le \E<1$
the lowest possible triangle density
$\T$ lies above the parabola $\T=\E(2\E-1)$, except for $\E_k=
(k-1)/k,\ k \ge 2$ when $\T_k$ lies on the parabola
and is uniquely
achieved by the complete balanced $k$-partite graph.
It is convenient to add two more points to the above sequence
$\{(\E_k, \T_k)\}$, namely $(\E_1, \T_1)=(0,0)$ corresponding to the
empty graph, and $(\E_\infty, \T_\infty)=(1,1)$ corresponding to the
complete graph. Furthermore, for $\E$ satisfying $\E_{k-1}<\E<
\E_{k},\ k\ge 1$, the lowest possible triangle density $\T$ lies on a
known curve (a `scallop'),
and the optimal graphs have a known structure 
\cite{Ra, PR, RS1, RS2}.
For all $0 < \E < 1$ the maximization of triangle density $\T$ for
given $\E$ is simpler than the above minimization: the maximum lies on
$\T=\E^{3/2}$ and is  uniquely achieved by a clique on enough vertices
to give the appropriate value of $\E$.
So extremal graph theory has determined the shape of the achievable
part of the parameter space as the region $R$ in
Fig.~\ref{FIG:phasefig}, and determined the structure of the graphs with
densities on the boundary.

\begin{figure}[!ht]
\centering
\rotatebox{90}{\hspace*{0.23\textwidth}{\Large $\tau$}}
\includegraphics[angle=0,width=0.46\textwidth]{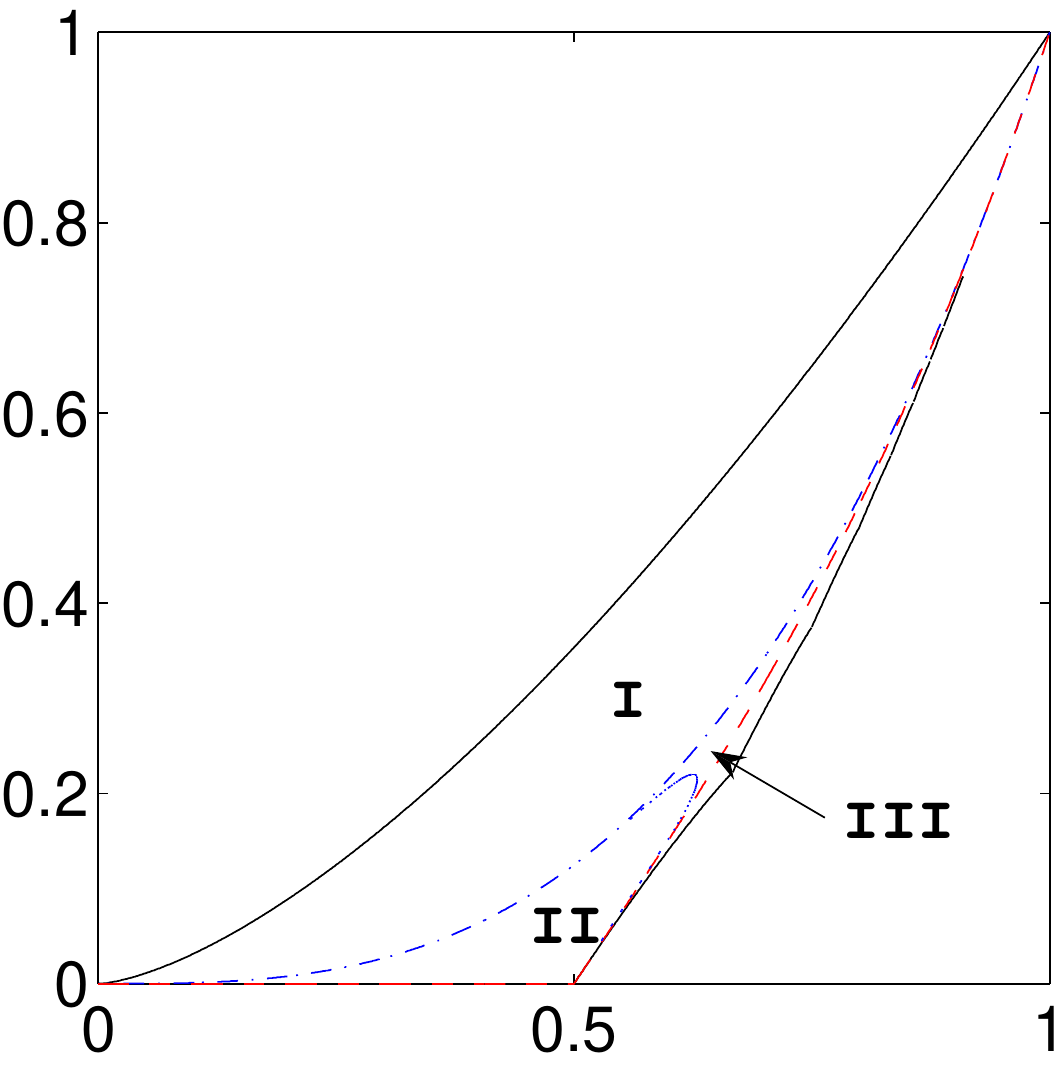}\\
\vspace*{-0.1cm}\hspace*{0.05\textwidth} {\Large $\epsilon$}
\caption{Boundary of phases I, II and III in the phase space. The blue dash-dotted
  line in the middle is the ER curve, the lower boundary of phase I.}
\label{FIG:Phase-Boundary}
\end{figure}

Our goal is to extend the study to the interior of
$R$ in a limited sense: for $(\E, \T)$ in the interior we wish to
determine, asymptotically in the size of the graph, what most graphs
are (or what a typical graph is) with these densities. (We clarify
`most' and `typical' below.) The typical graph
associated with $(\E, \T)$ will be described by probabilities of edges
between, and within, a finite number of vertex sets, parameters
which vary smoothly in $(\E, \T)$ except across phase transition
curves. Phases are the maximal connected regions in which the
parameters vary smoothly. One transition curve was determined in
\cite{RS2}, namely $(\E, \T)=(\E, \E^3),\ 0\le \E\le 1$. On this curve the typical
graph corresponds to edges chosen independently with probability $\E$. This
(`Erd\" os-R\' enyi' or `ER') curve separates $R$ into a high $\T$ region, a
single phase, and a low $\T$ region consisting of infinitely many
phases, at least one related to each scallop. We will present
simulation evidence for our determination of the typical graphs of the
high $\T$ phase (phase I) and for the two phases, II and III, which
are associated with the first scallop; see
Fig.~\ref{FIG:Phase-Boundary}.

There are two key features of these results: that by use of
probability one can capture the structure of most 
$(\E,\T)$-constrained graphs uniquely; and that their
structure is remarkably simple, for instance requiring at most four
parameters for over $95\%$ of the phase space.


We conclude this section with some notation which we use to clarify
our use above of the terms `most' and `typical' graphs.

Consider simple graphs $G$ with vertex set $V(G)$ of (labeled)
vertices, edge set $E(G)$ and triangle set $T(G)$, with cardinality
$|V(G)|=n$. Our main tool is $\displaystyle Z^{n,\alpha}_{\E,\T}$, the
number of graphs with densities:
\begin{equation}
e(G)\equiv \frac{|E(G)|}{{n \choose 2}}
\in (\E-\alpha,\E+\alpha); \quad
t(G)\equiv {\frac{|T(G)|}{{n\choose 3}}} \in (\T-\alpha,\T+\alpha).
\end{equation}

\nd From this we define the entropy density, the exponential rate of growth of 
$Z^{n,\alpha}_{e,t}$ as a function of $n$. First consider
\begin{equation}
s^{n,\alpha}_{\E,\T}={\ln(Z^{n,\alpha}_{\E,\T})\over n^2}, \hbox{ and
then } 
s(\E,\T)=\lim_{\alpha\downarrow 0}\lim_{n\to \infty}s^{n,\alpha}_{\E,\T}.
\end{equation}

\nd The double limit defining the entropy density $s(\E,\T)$ was proven to 
exist in \cite{RS1}. The objects of interest for us are the qualitative
features of $s(\E,\T)$ in the interior of $R$. To analyze them we make
use of a variational characterization of $s(\E,\T)$.

Let us first introduce some notation. Assume our vertex set $V$ is
decomposed into $M$ subsets: $V_1,\ldots,V_M$. We consider
(probabilistic) `multi\partite graphs' $G_A$ described by matrices 
$A=\{A_{ij}\,|\,
i,j=1,\ldots,M\}$ such that there is probability $A_{ij}$ of an edge
between any $v_i\in V_i$ and any $v_j\in V_j$. Special cases in which
each $A_{ij}\in \{0,1\}$ allow one to include nonprobabilistic graphs
in this setting.

This setting has been further generalized to analyze limits of
graphs as $n\to\infty$. (This work was recently
developed in \cite{LS1,LS2,BCLSV,BCL,LS3}; see also the recent book \cite{Lov}.)
The (symmetric) matrices $A_{ij}$ are
replaced by symmetric, measurable functions $g:(x,y)\in [0,1]^2\to
g(x,y)\in [0,1]$; the former are recovered by using a partition of
$[0,1]$ into consecutive subintervals to represent the partition of
$V$ into $V_1,\ldots,V_M$. The functions $g$ are called graphons, and
using them the following was proven in \cite{RS1} (adapting a
proof in \cite{CV}):
\vskip.1truein

\nd {\bf Variational Principle.}
  For any possible pair $(\E,\T)$, $s(\E,\T) = \max [-I(g)]$, where the
  maximum is over all graphons $g$ with $e(g)=\E$ and $t(g)=\T$,
  where 
\begin{equation}\label{EQ:E T}
e(g)=\int_{[0,1]^2}  g(x,y) \, dxdy, \qquad t(g) =
  \int_{[0,1]^3} g(x,y) g(y,z) g(z,x) \, dxdydz
\end{equation}
\nd  and the rate function is 
\begin{equation}\label{EQ:I}
I(g) = \int_{[0,1]^2} I_0[g(x,y)] \,dxdy,
\end{equation}
with the function $I_0(u)= {1\over 2} \left[ u \ln(u) + (1-u)\ln(1-u)\right]$. The existence of a maximizing 
graphon $g=g_{\E,\T}$ for any $(\E,\T)$ was proven in \cite{RS1},
again adapting a proof in \cite{CV}.

We want to consider two graphs \emph{equivalent} if they are obtained from
one another by relabelling the vertices. There is a generalized
version of this for graphons, with the relabelling replaced by
measure-preserving maps of $[0,1]$ into itself \cite{Lov}. The equivalence
classes of graphons are called reduced graphons, and on this space
there is a natural `cut metric' \cite{Lov}.

We can now clarify the notions of `most' and `typical' in the
introduction: if $g_{\E,\T}$ is the only reduced graphon maximizing
$s(\E,\T)$, then as the number $n$ of vertices diverges and $\alpha_n\to
0$, exponentially most graphs with densities
$e(G)
\in (\E-\alpha_n,\E+\alpha_n)$ and 
$t(G)\in (\T-\alpha_n,\T+\alpha_n)$
will have reduced graphon close to $g_{\E,\T}$ \cite{RS1}.

\vfill 
\eject

\section{Numerical Computation}
\label{SEC:numer}

To perform the study we outlined in the previous section, we combine numerical computation with local analysis. Let us first introduce our computational framework.

We need to solve the following constrained minimization problem:
\begin{equation}
	\min_{g} I(g)
\end{equation}
subject to the constraints
\begin{equation}\label{constraint-e-t}
	e(g)=\E, \quad\mbox{and}\quad t(g)=\T,
\end{equation}
where the functionals $I(g)$, $e(g)$ and $t(g)$ are defined in ~\eqref{EQ:I} and~\eqref{EQ:E T}.

To solve this minimization problem numerically, we need to represent the continuous functions with discrete values. We restrict ourselves to the class of piecewise constant functions.

For each integer $N\ge 1$ let $P_N: 0=p_0<p_1<p_2<\cdots<p_{N}=1$ be a partition of the interval
$[0,1]$. We denote by $c_i=p_i-p_{i-1}$ ($1\le i\le N$) the sizes of
the subintervals in the partition. Then we can form a partition of the
square $[0,1]^2$ using the (Cartesian) product 
$P_N\times P_N$. We define the class, $G(P_N)$, of those
symmetric functions $g:[0,1]^2\to[0,1]$ which are piecewise constant on the
subsets $P_N\times P_N$, and introduce the notation:
\begin{equation}
	g(x,y)=g_{ij},\ \ (x,y)\in(p_{i-1},p_i)\times (p_{j-1},p_j), \quad 1\le i,j\le N,
\end{equation}
with $g_{ij}=g_{ji}$. They are probabilistic generalizations of
multipartite graphs, which we call multi\partite\!\!: bi\partite for $N=2$,
tri\partite for $N=3$, etc. We will call such a multi\partite graphon `symmetric' if all
$c_j$ are equal and all $g_{jj}$ are equal, and otherwise `asymmetric'.

It is easy to check that with this type of $g$, the functionals $I$, $e$ and $t$ become respectively
\begin{equation}
	I(g)=\dsum_{1\le i,j\le N}I_0(g_{ij})c_ic_j=\dfrac{1}{2} \dsum_{1\le i,j\le N} [g_{ij}\ln g_{ij}+(1-g_{ij})\ln(1-g_{ij})]c_i c_j,
\end{equation}
\begin{equation}
	e(g)=\dsum_{1\le i, j\le N} g_{ij}c_i c_j,\qquad t(g)=\dsum_{1\le i,j,k\le N} g_{ij}g_{jk} g_{ki}c_i c_j c_k.
\end{equation}

Using the fact that $\cup_N G(P_N)$ is dense in the space of graphons, 
our objective is to solve the minimization problem in $G(P_N)$, for
all $N$:
\begin{equation}
\min_{\{c_j\}_{1\le j\le N}, \{g_{ij}\}_{1\le i,j\le N}} I(g),
\end{equation}
with the constraints (\ref{constraint-e-t}) and
\begin{equation}\label{EQ:Constraints}
	0\le g_{ij}, c_j\le 1, \qquad \dsum_{1\le j\le N}c_j=1, \qquad \mbox{and} \qquad g_{ij}=g_{ji}
\end{equation}
for any given pair of $\E$ and $\T$ values.

\paragraph{Minimum without the $\E$ and $\T$ constraints.} 
Let us first look at the minimization problem without the $e(g)=\E$ and $t(g)=\T$ constraints. In this case, we can easily calculate the first-order variations of the objective function:
\begin{equation}\label{EQ:I Grad g}
	I_{g_{ij}}=I_0'(g_{ij}) c_i c_j=\dfrac{1}{2}\ln[g_{ij}/(1-g_{ij})] c_ic_j,
\end{equation}
\begin{equation}\label{EQ:I Grad c}
	I_{c_p}=\sum_{1\le i, j\le N} I_0(g_{ij}) (\dfrac{\partial c_i}{\partial c_p} c_j+ c_i \dfrac{\partial c_j}{\partial c_p}),\ \ \mbox{with}\ \ 
	\sum_{1\le j\le N} \dfrac{\partial c_j}{\partial c_p} = 0. 
\end{equation}
where we used the constraint $\sum_{1\le j\le N} c_j=1$ to get the last equation. 

We observe from~\eqref{EQ:I Grad g} that when no constraint is imposed, the minimum of $I(g)$ is located at the constant graphon with $g_{ij}=1/2$,\ $1\le i, j \le N$. The value of the minimum is $-(\ln 2)/2=-0.34657359027997$. At the minimum, $\E=1/2$ and $\T=1/8$ (on the ER curve $\T=\E^3$). We can check also from~\eqref{EQ:I Grad c} that indeed $I_{c_p}=0$ at the minimum for arbitrary partition $\{c_j\}$ that satisfies $\sum_{1\le j\le N} c_j=1$. This is obvious since the graphon is constant and thus the partition does not play a role here.

The second-order variation can be calculated as well:
\begin{equation}\label{EQ:I Hess gg}					
I_{g_{ij}g_{lm}}= I_0''(g_{ij}) c_ic_j \delta_{il}\delta_{jm}=\dfrac{1}{2}\frac{c_ic_j}{g_{ij}(1-g_{ij})}\delta_{il}\delta_{jm}.
\end{equation}
\begin{equation}\label{EQ:I Hess gc}					
I_{g_{ij}c_p}=I_{c_p g_{ij}}
=I_0'(g_{ij}) (\dfrac{\partial c_i}{\partial c_p} c_j+ c_i \dfrac{\partial c_j}{\partial c_p}),\ \ \mbox{with}\ \ 
	\sum_{1\le i\le N} \dfrac{\partial c_i}{\partial c_p} = 0. 
\end{equation}
\begin{multline}\label{EQ:I Hess cc}
	I_{c_p c_q}=\sum_{1\le i, j\le N} I_0(g_{ij}) \Big(\dfrac{\partial^2 c_i}{\partial c_pc_q} c_j+\dfrac{\partial c_i}{\partial c_p} \dfrac{\partial c_j}{\partial c_q}+ \frac{\partial c_i}{\partial c_q} \dfrac{\partial c_j}{\partial c_p}+ c_i \dfrac{\partial^2 c_j}{\partial c_pc_q}\Big),\\ \mbox{with}\ \ 
	\sum_{1\le i\le N} \dfrac{\partial c_i}{\partial c_p} = 0,\ \ 
	\sum_{1\le i\le N} \dfrac{\partial^2 c_i}{\partial c_p c_q} = 0. 
\end{multline}
Thus the unconstrained minimizer is stable with respect to perturbation in $g$ values since the second variation $I_{g_{ij}g_{lm}}$ is positive definite. At the minimizer, however, the second variation with respect to $c$ is zero! This is consistent with what we know already.

We now propose two different algorithms to solve the constrained minimization problem. The first algorithm is based on Monte Carlo sampling while the second algorithm is a variant of Newton's method for nonlinear minimization. All the numerical results we present later are obtained with the sampling algorithm and confirmed with the Newton's algorithm; see more discussion in Section~\ref{SEC:Strategy}.

\subsection{Sampling Algorithm}
\label{SEC:Sampling}

In the sampling algorithm, we construct random samples of values of $I$ in the parameter space. We then take the minimum of the sampled values. The algorithm works as follows.
\begin{itemize}
	\item[{[0]}] Set the number $N$ and the number of samples to be constructed ($L$); Set counter $\ell =1$;
	\item[{[1]}] Generate the sizes of the blocks: $\{c_j^\ell\}\in[0, 1]^N$ and normalize so that $\sum_j c_j^\ell =1$;
	\item[{[2]}] Generate a sample $\{g_{ij}^\ell\}\in[0, 1]^{N\times N}$ with $g_{ij}^\ell=g_{ji}^\ell$, $1\le i,j\le N$ and rescale such that:
		\begin{itemize}
			\item[{[A]}] $\dsum_{1\le i,j\le N}g_{ij}^\ell c_i^\ell c_j^\ell=\E$;
			\item[{[B]}] $\dsum_{1\le i,j,k\le N}^\ell g_{ij}g_{jk}^\ell g_{ki}^\ell c_i^\ell c_j^\ell c_k^\ell=\T$;
		\end{itemize}
	\item[{[3]}] Evaluate $I_\ell=\dsum_{1\le i,j\le N} I_0(g_{ij}^\ell)c_i^\ell c_j^\ell$;
	\item[{[4]}] Set $\ell=\ell+1$; Go back to [1] if $\ell\le L$;
	\item[{[5]}] Evaluate $I_N=\min\{I_\ell\}_{\ell=1}^L$.
\end{itemize}
We consider two versions of Step [2] which work equally well (beside a slight difference in computational cost) in practice. In the first version, we generate a sample $\{g_{ij}^\ell\}\in[0, 1]^{N\times N}$ with $g_{ij}^\ell=g_{ji}^\ell$, $1\le i,j\le N$. We then rescale the sample using the relation $g_{ij}^\ell=\gamma g_{ij}^\ell+\tilde\gamma$. The relations [A] and [B] then give two equations for the parameter $\gamma$ and $\tilde\gamma$. We solve the equations for $\gamma$ and $\tilde\gamma$. If at least one of $\{g_{ij}^\ell\}$ violate the condition $g_{ij}^\ell\in[0,1]$ after rescaling, we re-generate a sample and repeat the process until we find a sample that satisfies $g_{ij}^\ell\in[0,1]$ ($1\le i,j\le N$) after rescaling. In the second version, we simply find a sample $\{g_{ij}^\ell\}$ by solving [A] and [B] as a third-order algebraic system using a multivariate root-finding algorithm, with the constraint that $\{g_{ij}^\ell\}\in[0, 1]^{N\times N}$. When multiple roots are found, we take them as different qualified samples.

This sampling algorithm is a global method in the sense that the algorithm will find a good approximation to the global minimum of the functional $I$ when sufficient samples are constructed. The algorithm will not be trapped in a local minimum. The algorithm is computationally expensive. However, it can be parallelized in a straightforward way. We will discuss the issues of accuracy and computational cost in Section~\ref{SEC:Strategy}.

\subsection{SQP Algorithm}
\label{SEC:Min}

The second algorithm we used to solve the minimization problem is a sequential quadratic programming (SQP) method for constrained optimization. To briefly describe the algorithm, let us denote the unknown by $x=(c_1,\cdots,c_N,g_{11},\cdots,g_{1N},g_{21},\cdots,g_{2N},\cdots,g_{N1},\cdots, g_{NN})$. Following~\cite{GMSW}, we rewrite the optimization problem as
\begin{equation}
	\min_{x\in[0,1]^{N+N^2}}I(x),\quad \mbox{subject to}\quad l\le r(x)\le u,\qquad r(x)=\left(
\begin{array}{c}x\\ \cA x\\ \cC(x)\end{array}\right)
\end{equation}
where $l$ and $u$ are the lower and upper bounds of $r(x)$ respectively, $\cC_1(x)=e(x)$ and $\cC_2(x)=t(x)$. The matrix $\cA$ is used to represent the last two linear constraints in~\eqref{EQ:Constraints}:
\begin{equation}
	\cA=\left(\begin{array}{cc} \Sigma & \bzero_{1\times N^2}\\ \bzero_{N^2 \times N} & \cS \end{array}\right)
\end{equation}
where $\Sigma=(1,\cdots,1)\in\bbR^{1\times N}$ is used to represent
the constraint $\sum c_j=1$, $\bzero_{1\times N^2} \in\bbR^{1\times
  N^2}$ is a matrix with all zero elements, and $\cS\in\bbR^{N^2\times
  N^2}$ is a matrix used to represent the symmetry constraint
$g_{ij}=g_{ji}$. The elements of $\cS$ are given as follows. For any
index $k$, we define the conjugate index $k'$ as $k'=(r-1)N+q$ with
$q$ and $r$ the unique integers such that $k=(q-1)N+r$ ($0\le r<N$). Then for all $1\le k\le N^2$, $\cS_{kk}=-\cS_{kk'}=1$ if $k\ne (j-1)N+1$ for some $j$, and $\cS_{kk}=0$ if $k\ne (j-1)N+1$ for some $j$. All other elements of $\cS$ are zero.

The SQP algorithm is characterized by the following iteration
\begin{equation}
	x_{\ell+1}=x_{\ell}+\alpha_\ell p_\ell,\qquad \ell\ge 0
\end{equation}
where $p_\ell$ is the search direction of the algorithm at step $\ell$ and $\alpha_\ell$ is the step length. The search direction $p_\ell$ in SQP is obtained by solving the following constrained quadratic problem
\begin{equation}
	\min I(x_\ell)+\frak g(x_\ell)^\fT p_\ell+\frac{1}{2} p_\ell^\fT H_\ell p_\ell, \qquad \mbox{subject to}\qquad l\le r(x_\ell)+ J(x_\ell)p_\ell\le u.
\end{equation}
Here $\frak g(x_\ell)=\nabla_x I$, whose components are given analytically in~\eqref{EQ:I Grad g} and \eqref{EQ:I Grad c}, $J(x_\ell)=\nabla_x r(x_\ell)$, whose components are given by
\begin{eqnarray}
\label{EQ:e Grad g}	e_{g_{ij}}&=& c_i^\ell  c_j^\ell,\\
\label{EQ:t Grad g}	t_{g_{ij}}&=& 3\dsum_{1\le k\le N}g_{jk}^\ell g_{ki}^\ell c_i^\ell c_j^\ell  c_k^\ell ,\\
\label{EQ:e Grad c}	e_{c_p}&=& \sum_{1\le i, j\le N}g_{ij}^\ell \Big(\dfrac{\partial c_i^\ell }{\partial c_p^\ell } c_j^\ell + c_i^\ell  \dfrac{\partial c_j^\ell }{\partial c_p^\ell }\Big),\ \ \mbox{with}\ \ 
	\sum_{1\le i\le N} \dfrac{\partial c_i^\ell }{\partial c_p^\ell } = 0\\
\label{EQ:t Grad c}	t_{c_p}&=& \dsum_{1\le i,j, k\le N}g_{ij}^\ell g_{jk}^\ell g_{ki}^\ell \Big(\dfrac{\partial c_i^\ell }{\partial c_p^\ell } c_j^\ell c_k^\ell + \dfrac{\partial c_j^\ell }{\partial c_p^\ell } c_i^\ell c_k^\ell +\dfrac{\partial c_k^\ell }{\partial c_p^\ell } c_i^\ell c_j^\ell \Big).
\end{eqnarray}
$H(x_\ell)$ is a positive-definite quasi-Newton approximation to the Hessian of the objective function. We take the BFGS updating rule to form $H(x_\ell)$ starting from the identity matrix at the initial step~\cite{NW}.

We implemented the SQP minimization algorithm using the software package given in~\cite{GMSW} which we benchmarked with the \verb|fmincon| package in MATLAB R2012b. 

\subsection{Computational Strategy}
\label{SEC:Strategy}

The objective of our calculation is to minimize the rate function $I$
for a fixed $(\E,\T)$ pair over the space of all graphons. Our
computational strategy is to first minimize for $g\in G(P_N)$, for a
fixed number of blocks $N$, and then minimize over the number of
blocks. Let $I_N$ be the minimum achieved by the graphon $g_N\in G(P_N)$,
then the minimum of the original problem is $I=\min_N\{I_N\}$.  Due to
limitations on computational power we can only solve up to $N=16$. The algorithms, however, are not limited by this. 

By construction we know that $I_2\ge I_3\ge \cdots\ge I_{N}$. What is surprising is that our computation suggests that the minimum is always achieved with bi\partite graphons in the phases that we are considering. In other words, $I_2=\min\{I_N\}_{N=2}^{N=16}$. 

To find the minimizers of $I_N$ for a fixed $N$, we run both the
sampling algorithm and the SQP algorithm. We observe that both
algorithms give the same results (to precision $10^{-6}$) in all cases
we have simulated. In the sampling algorithm, we observe
from~\eqref{EQ:I Grad g} and ~\eqref{EQ:I Grad c} that the changes in
$I$ caused by perturbations in $g_{ij}$ and $c_p$ are given
respectively by $\delta I/\delta g_{ij}\sim \ln[\delta
  g_{ij}/(1-\delta g_{ij})]c_i c_j$ and $\delta I/\delta
  c_p\sim I_0(g_{ij})\delta c_p$. These mean that to get an accuracy
of order $\eta$, our samples have to cover a grid of parameters
with mesh size $\delta g_{ij}\sim
{e^{{2\eps}/{c_ic_j}}}/{(1+e^{{2\eps}/{c_ic_j}})}
\sim 1$ in
the $g_{ij}$ direction and $\delta c_p\sim {\eps}/{I_0(g_{ij})}$
in the $c_p$ direction. Since $I_0(g_{ij})\in[-(\ln 2)/2,0]$,
we have $\eps/I_0(g_{ij})>\eps$. It is thus enough to sample
on a grid of size $\eps$. Similar analysis following ~\eqref{EQ:e Grad
  g}, ~\eqref{EQ:t Grad g}, ~\eqref{EQ:e Grad c} and ~\eqref{EQ:t Grad
  c} shows that, to achieve an accuracy $\eps$ on the $\E$ and $\T$
constraints, we need to sample on grids of size at most on the order
of $\eps$ in both the $g_{ij}$ and $c_p$ directions. The total computational complexity is thus $\sim (1/\eps)^{N+N^2}$ in terms of function evaluations.

To run the SQP algorithm, for each $(\E,\T)$ constraint, we start from
a collection of $L_c^N\times L_g^{N^2}$ initial guesses. These initial
guesses are generated on the uniform grid of $L_c$ intervals in each
$c_p$ direction and $L_g$ intervals in each $g_{ij}$ direction. They
are then rescaled linearly (if necessary) to satisfy the $\E$ and $\T$
constraints. The results of the algorithm after convergence are
collected to be compared with the sampling algorithm. We observe that in most cases, the algorithm converges to identical results starting from different initial guesses.

Finally we note that there are not many parameters we need to tune to
get the results that we need. We observe that the SQP algorithm is
very robust using the general default algorithmic parameters. The only
parameters that we can adjust are $L_c$ and $L_g$ which control how
many initial guesses we want to run. Our calculation shows that
$L_c=L_g=10$ is enough for all the cases we studied. When we increased
$L_c$ and $L_g$ to get more initial guesses, we did not gain any new
minimizers.

\subsection{Benchmarking the Computations}

Before using the computational algorithms to explore the regions of the phase space that we plan to explore, we first benchmark our codes by reproducing some theoretically known results. 

\begin{figure}[ht]
\centering
\includegraphics[angle=0,width=0.23\textwidth]{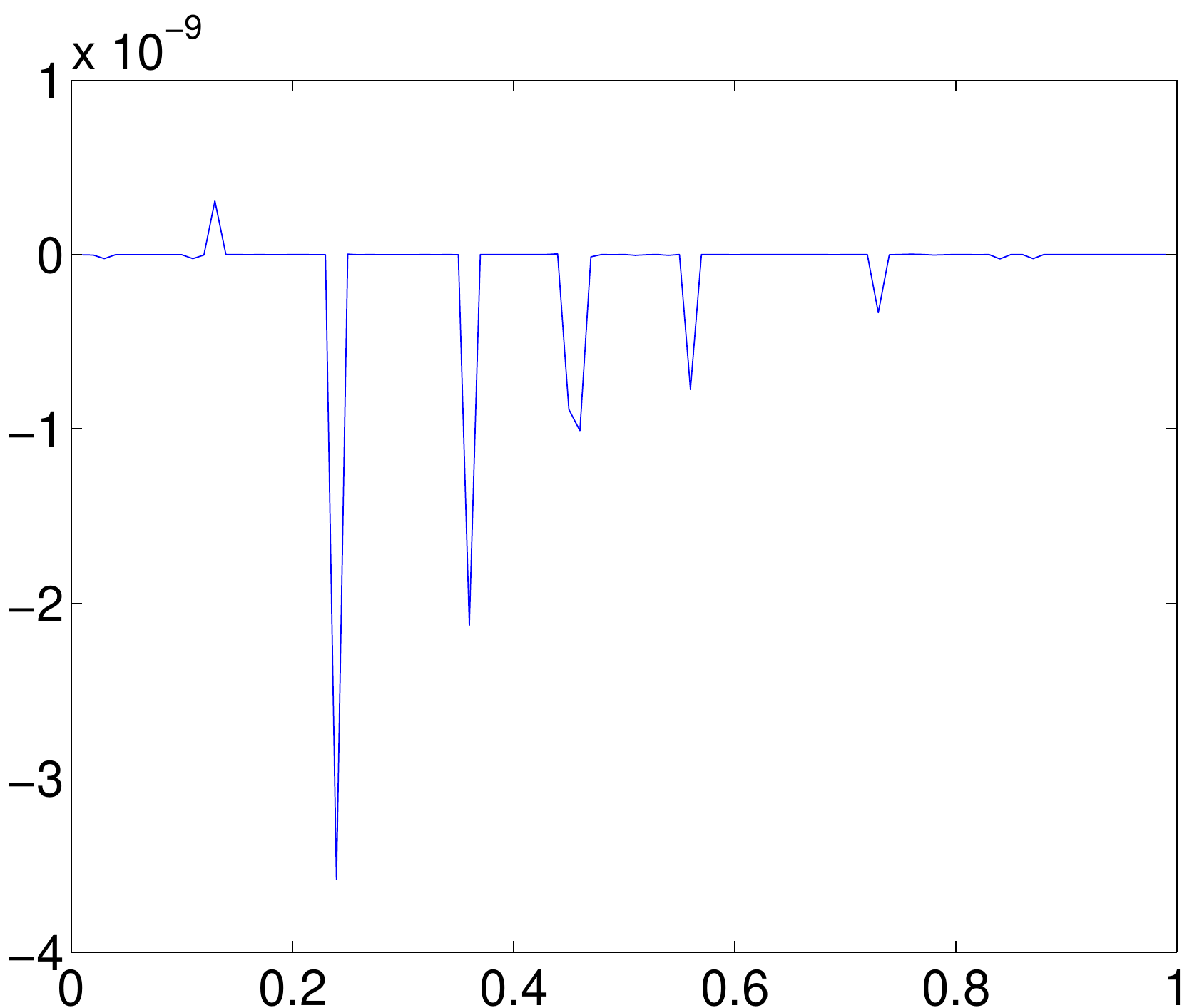} 
\includegraphics[angle=0,width=0.23\textwidth]{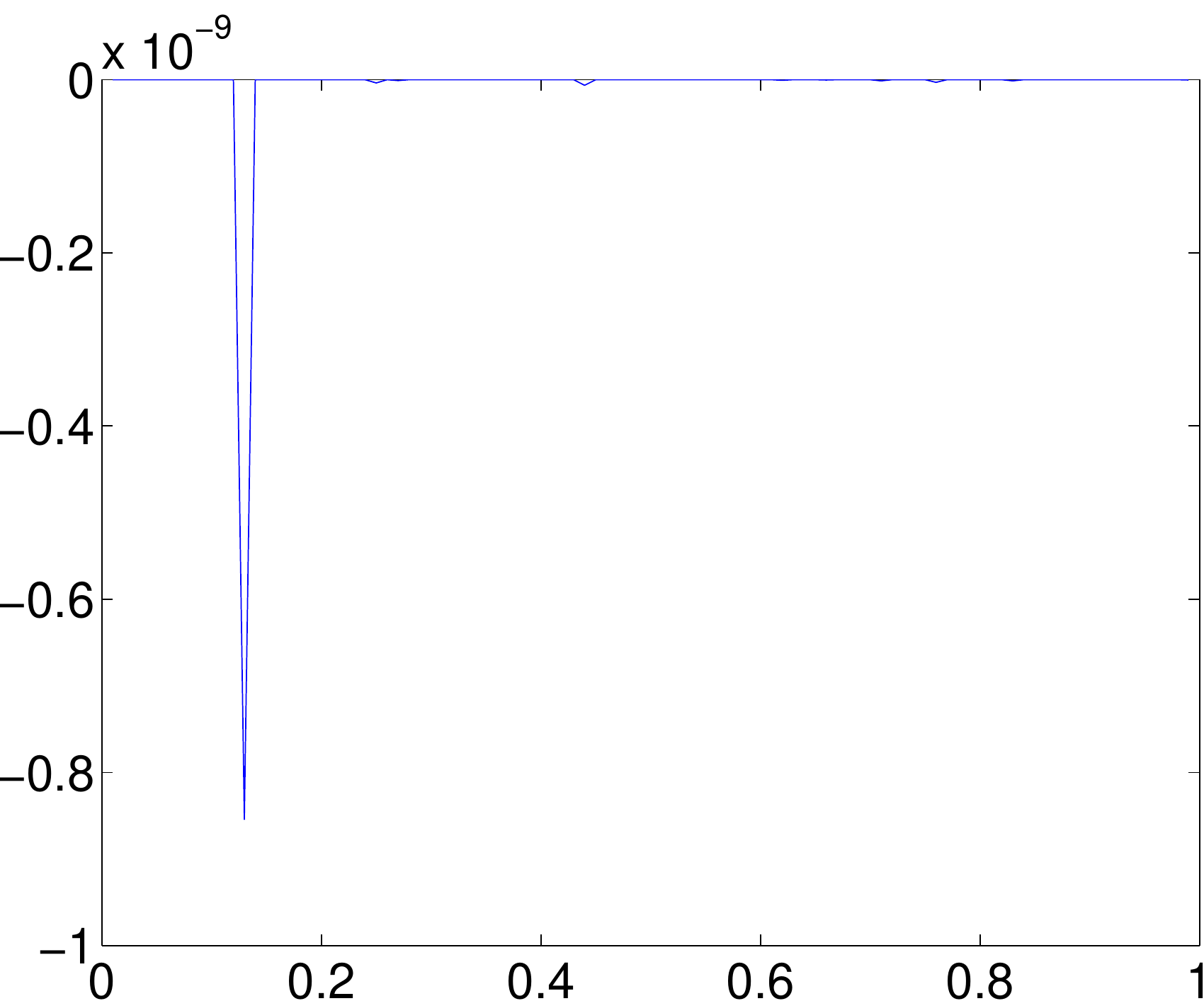} 
\includegraphics[angle=0,width=0.23\textwidth]{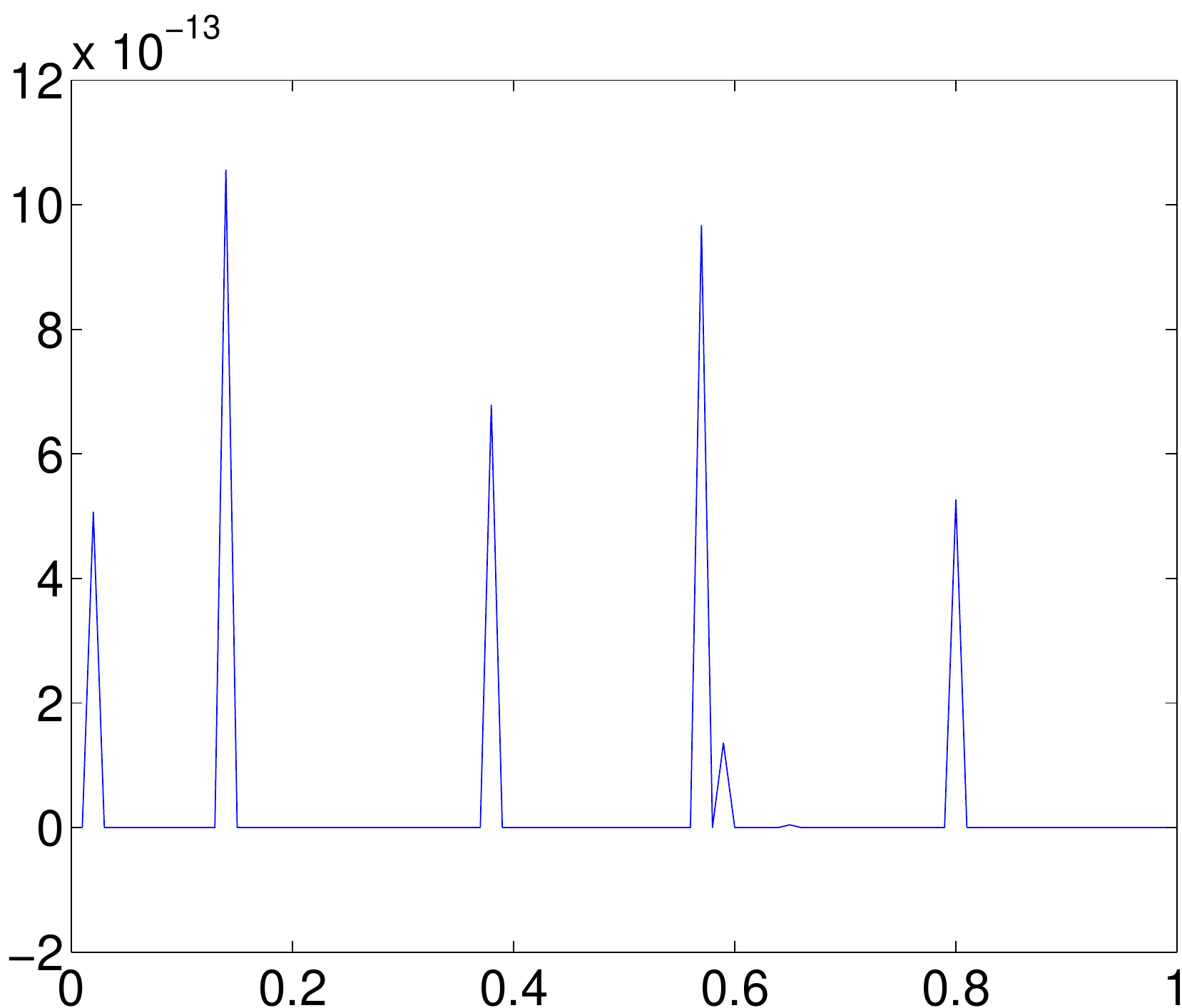} 
\includegraphics[angle=0,width=0.23\textwidth]{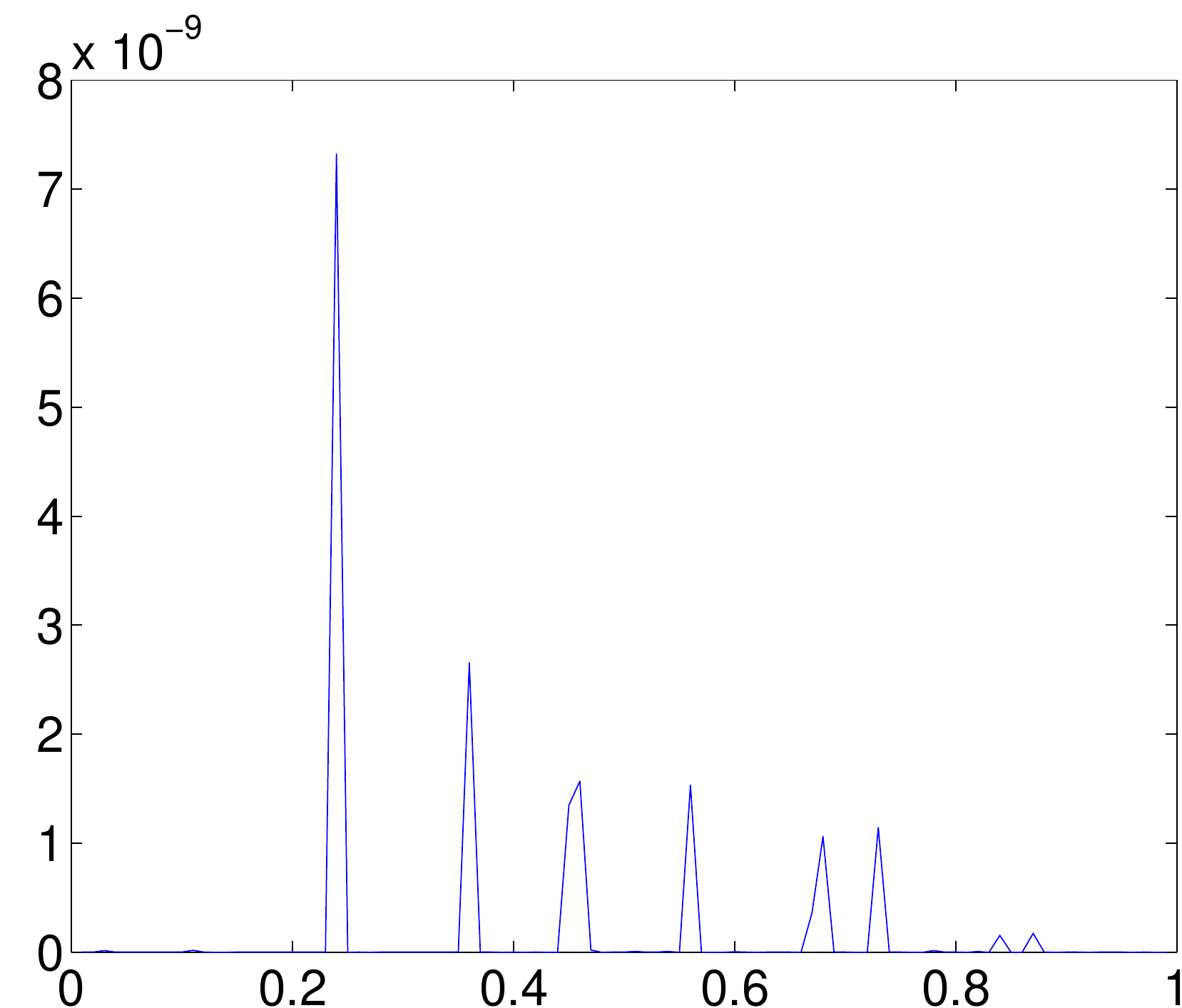}\\ 
\vspace*{-0.1cm}\hspace*{-0.0\textwidth} {$\epsilon$}\hspace*{0.22\textwidth} {$\epsilon$}\hspace*{0.22\textwidth} {$\epsilon$}\hspace*{0.22\textwidth} {$\epsilon$}
\caption{Difference between numerically computed (with superscript $num$)  and true (with superscript $true$) minimizing
  graphons on the curve $\T=\E^{\frac{3}{2}}$. From left to right:
  $c^{true}-c^{num}$, $g_{11}^{true}-g_{11}^{num}$, $g_{22}^{true}-g_{22}^{num}$ and $g_{12}^{true}-g_{12}^{num}$.}
\label{FIG:Benchmark-Upper}
\end{figure}

\paragraph{On the upper boundary.} We first reproduce minimizing
graphons on the curve $\T=\E^{\frac{3}{2}}$. It is known~\cite{RS2} that the minimizing graphons are equivalent to the bi\partite graphon with $c=\sqrt{\E}$, $g_{11}=1$ and $g_{12}=g_{21}=g_{22}=0$. The minimum value of the rate function is $I=0$. 
In Fig.~\ref{FIG:Benchmark-Upper} we show the difference between the simulated minimizing graphons and the true minimizing graphons given by the theory. We observe that the difference is always well below $10^{-6}$ for all the results on the curve.

\paragraph{On the segment $(\E,\T)=(0,0.5)\times \{0\}$.} It is known that on the segment $(\E,\T)=(0,0.5)\times \{0\}$ the minimizing graphons are symmetric bi\partite with $g_{11}=g_{22}=0$ and $g_{12}=g_{21}=2\E$. 

\begin{figure}[!ht]
\centering
\includegraphics[angle=0,width=0.23\textwidth]{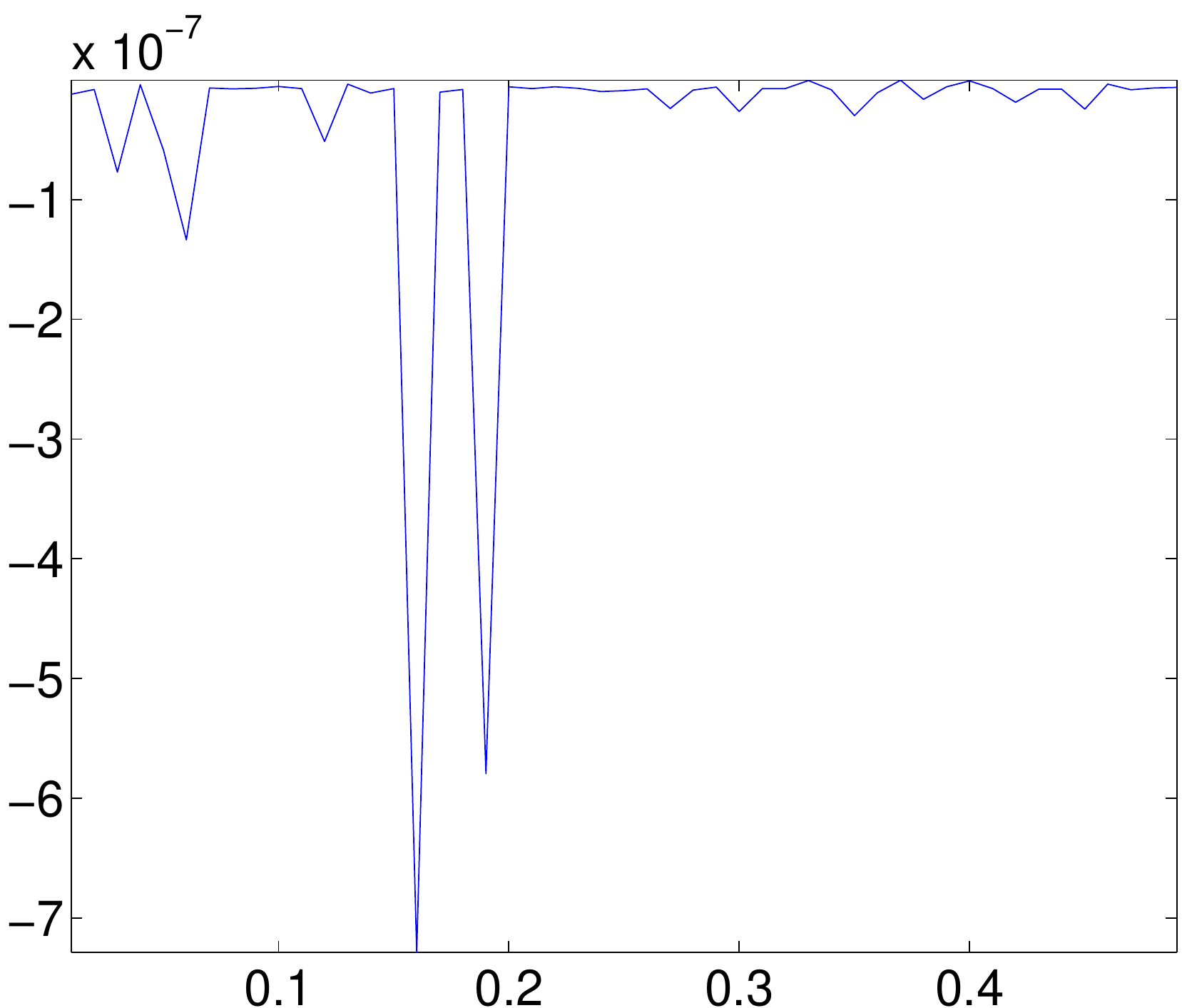}
\includegraphics[angle=0,width=0.23\textwidth]{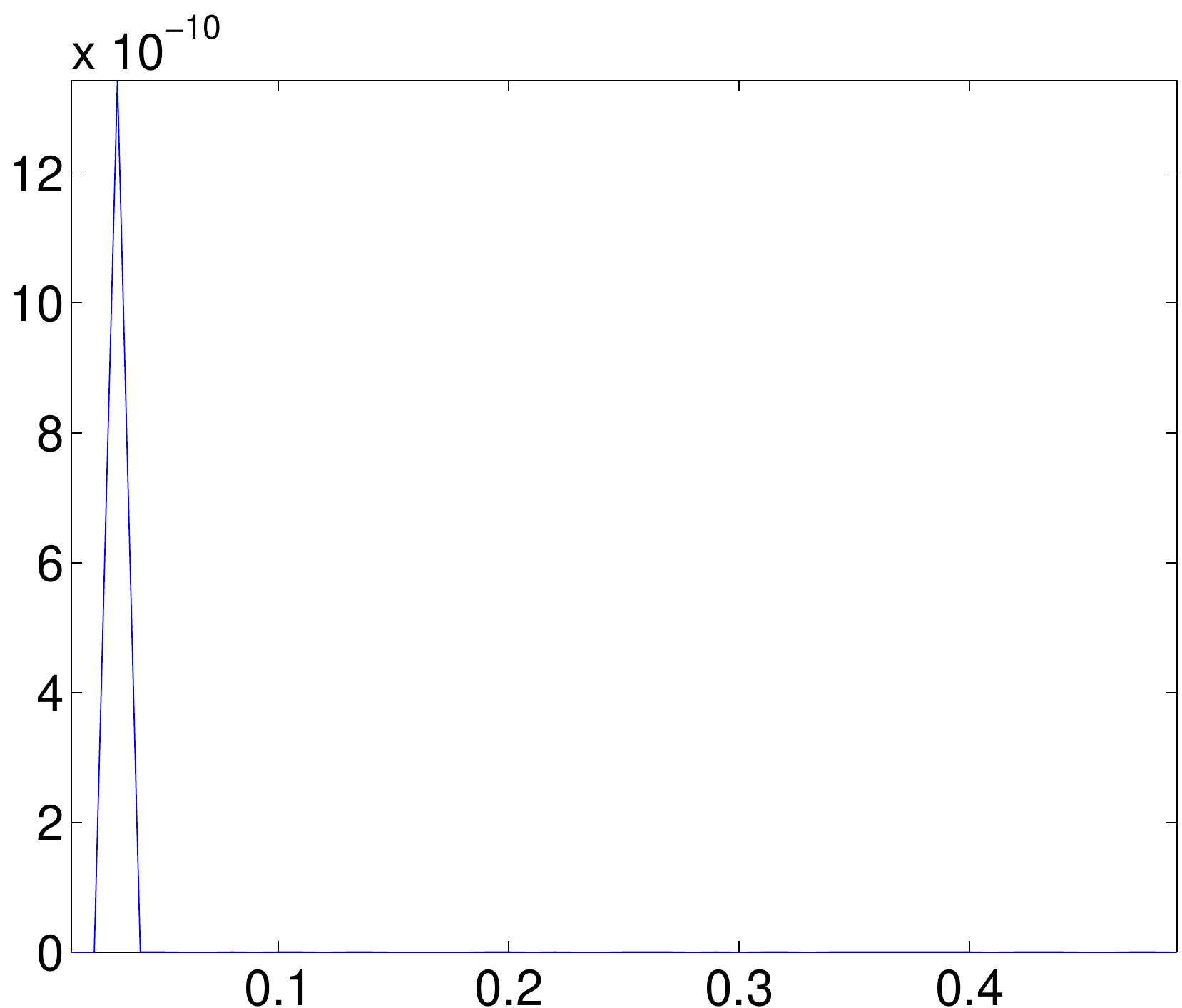} 
\includegraphics[angle=0,width=0.23\textwidth]{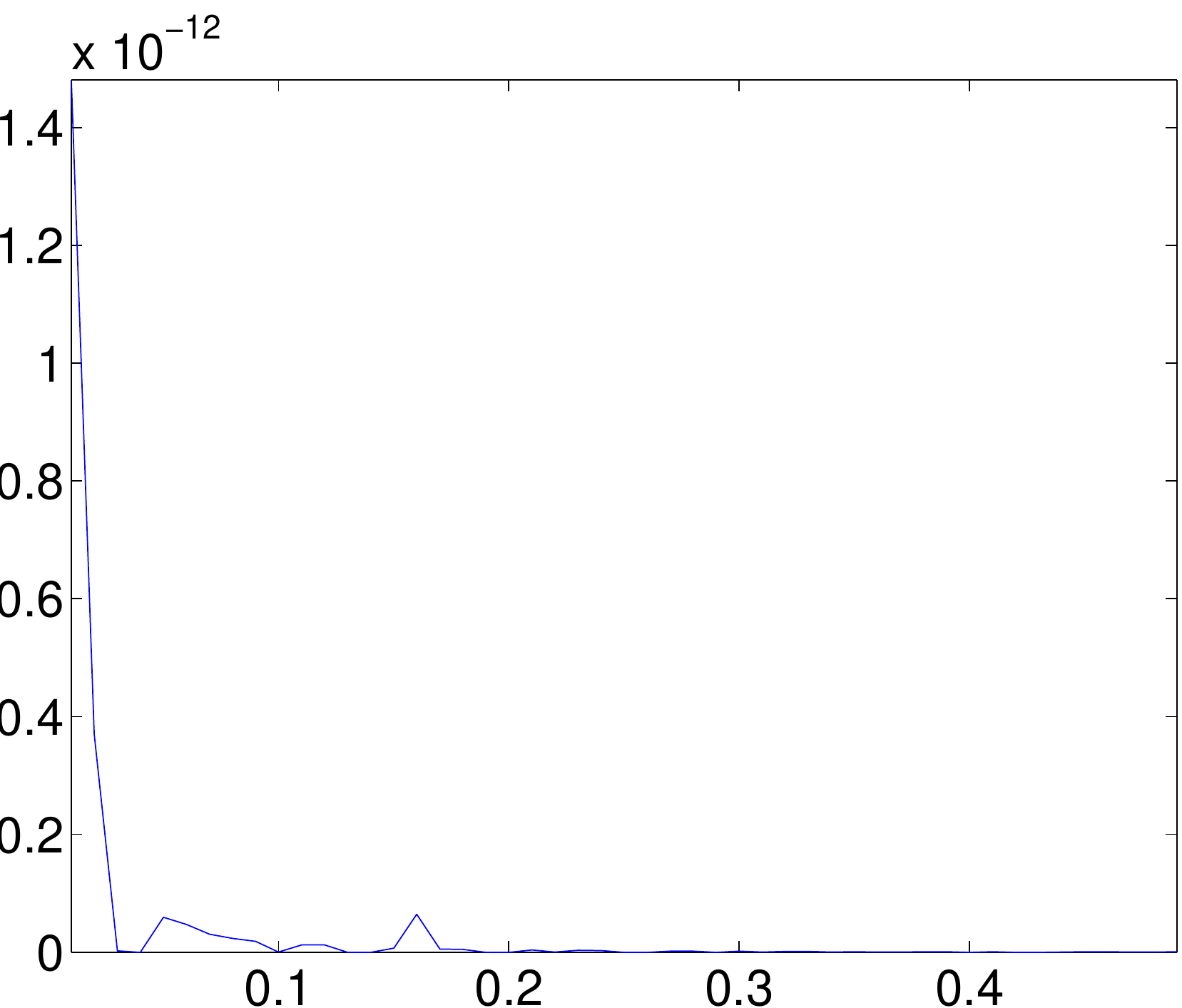} 
\includegraphics[angle=0,width=0.23\textwidth]{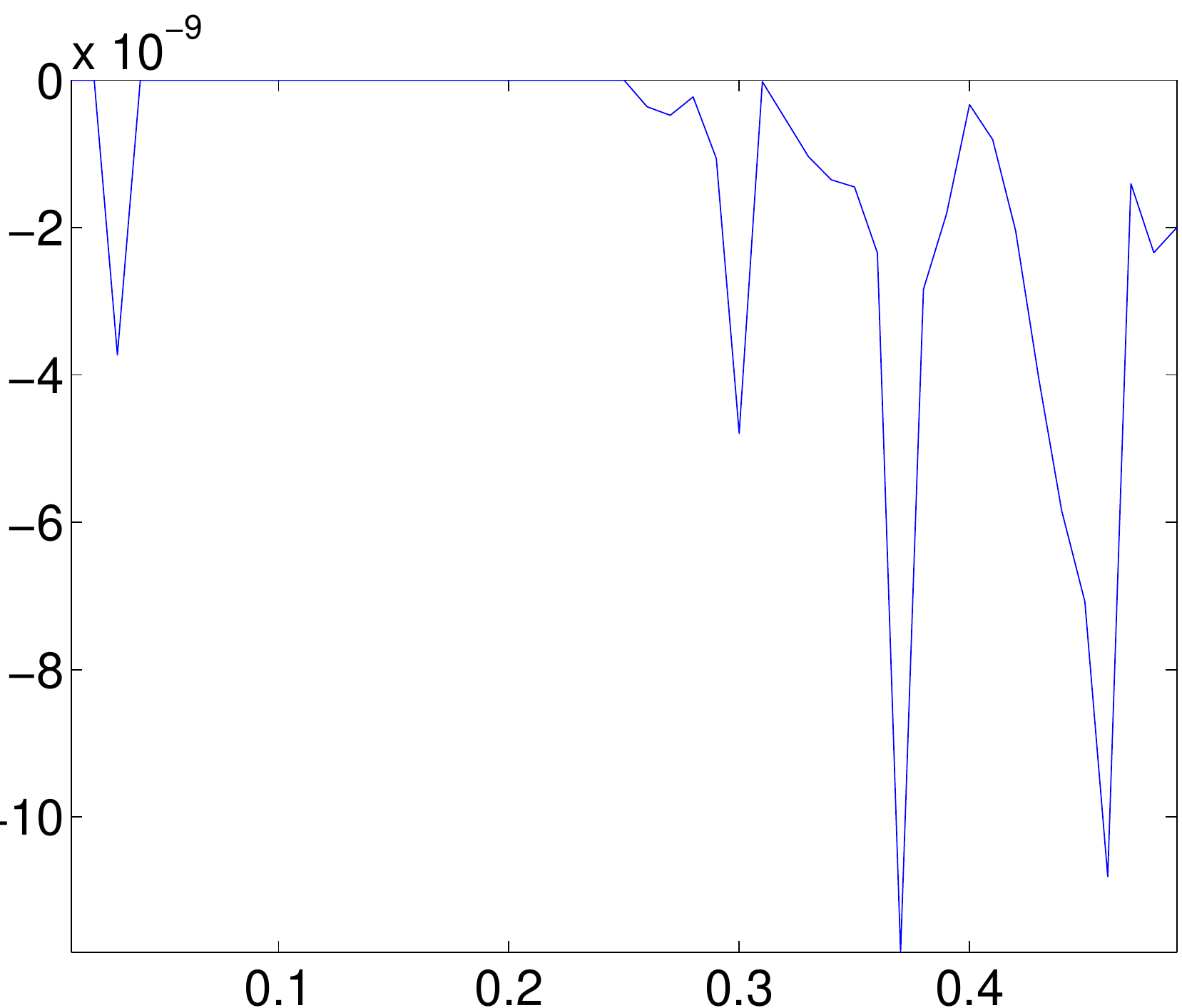}\\ 
\vspace*{-0.1cm}\hspace*{-0.0\textwidth} {$\epsilon$}\hspace*{0.22\textwidth} {$\epsilon$}\hspace*{0.22\textwidth} {$\epsilon$}\hspace*{0.22\textwidth} {$\epsilon$}
\caption{Difference between numerically computed (with superscript $num$) and true (with superscript $true$) minimizing graphons on the line segment $(\E,\T)=(0,0.5)\times \{0\}$. From left to right: $c^{true}-c^{num}$, $g_{11}^{true}-g_{11}^{num}$, $g_{22}^{true}-g_{22}^{num}$ and $g_{12}^{true}-g_{12}^{num}$.}
\label{FIG:Benchmark-Scallop1}
\end{figure}

In Fig.~\ref{FIG:Benchmark-Scallop1} we show the difference between the simulated minimizing graphons and the true minimizing graphons given by the theory~\cite{RS2}. The differences are again very small, $<10^{-6}$. This shows again that our numerical computations are fairly accurate.

\begin{figure}[!ht]
\centering
\includegraphics[angle=0,width=0.23\textwidth]{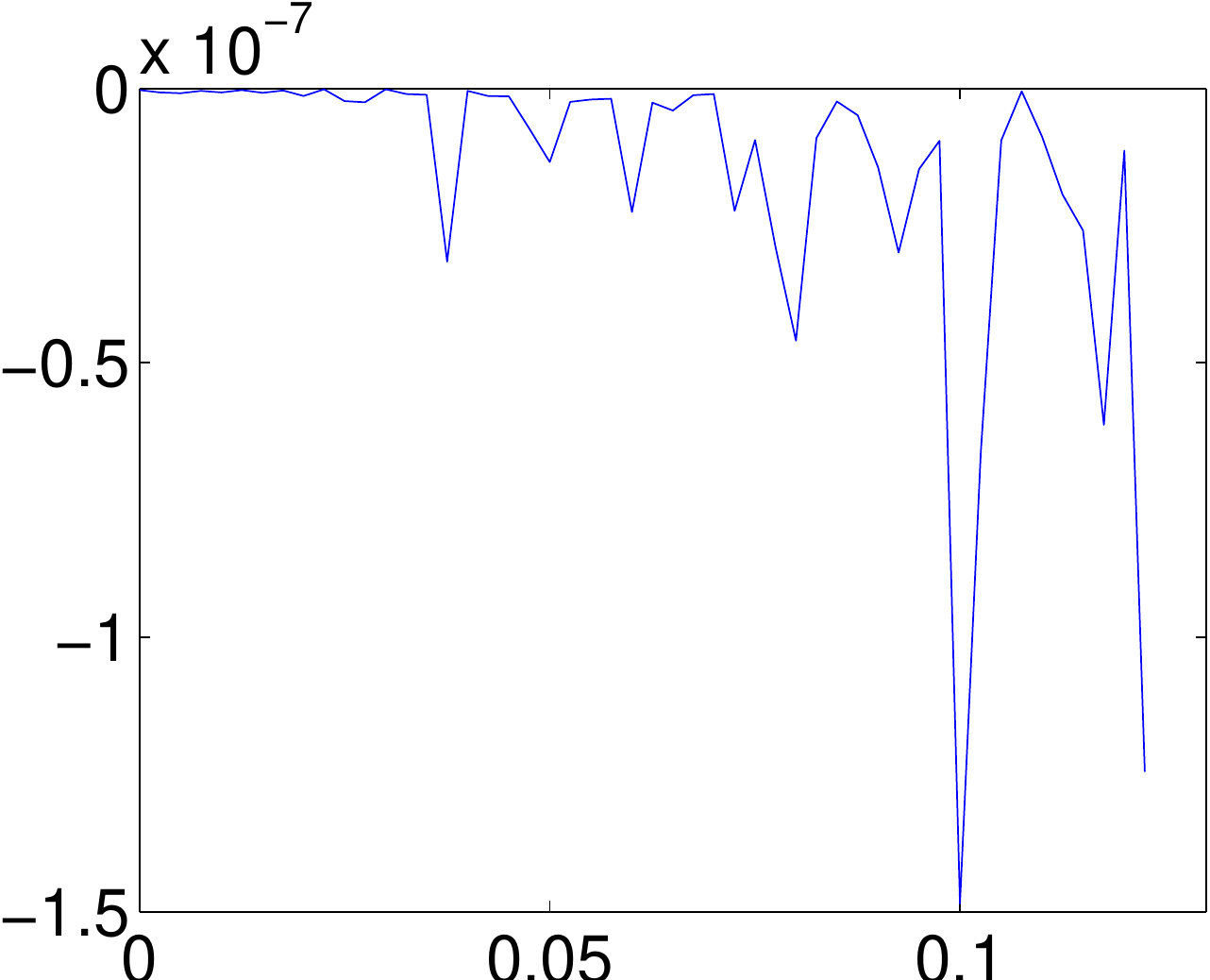} 
\includegraphics[angle=0,width=0.23\textwidth]{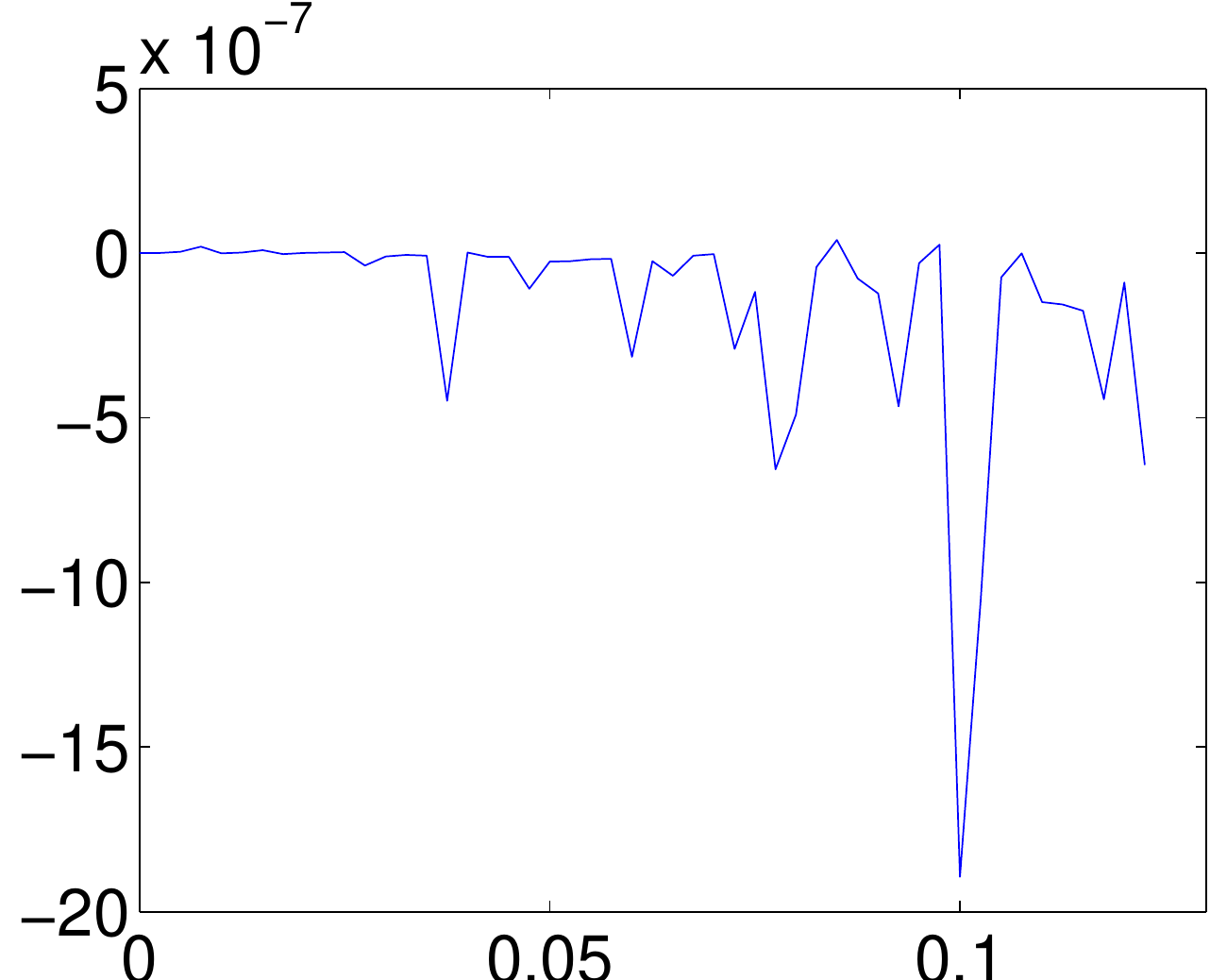} 
\includegraphics[angle=0,width=0.23\textwidth]{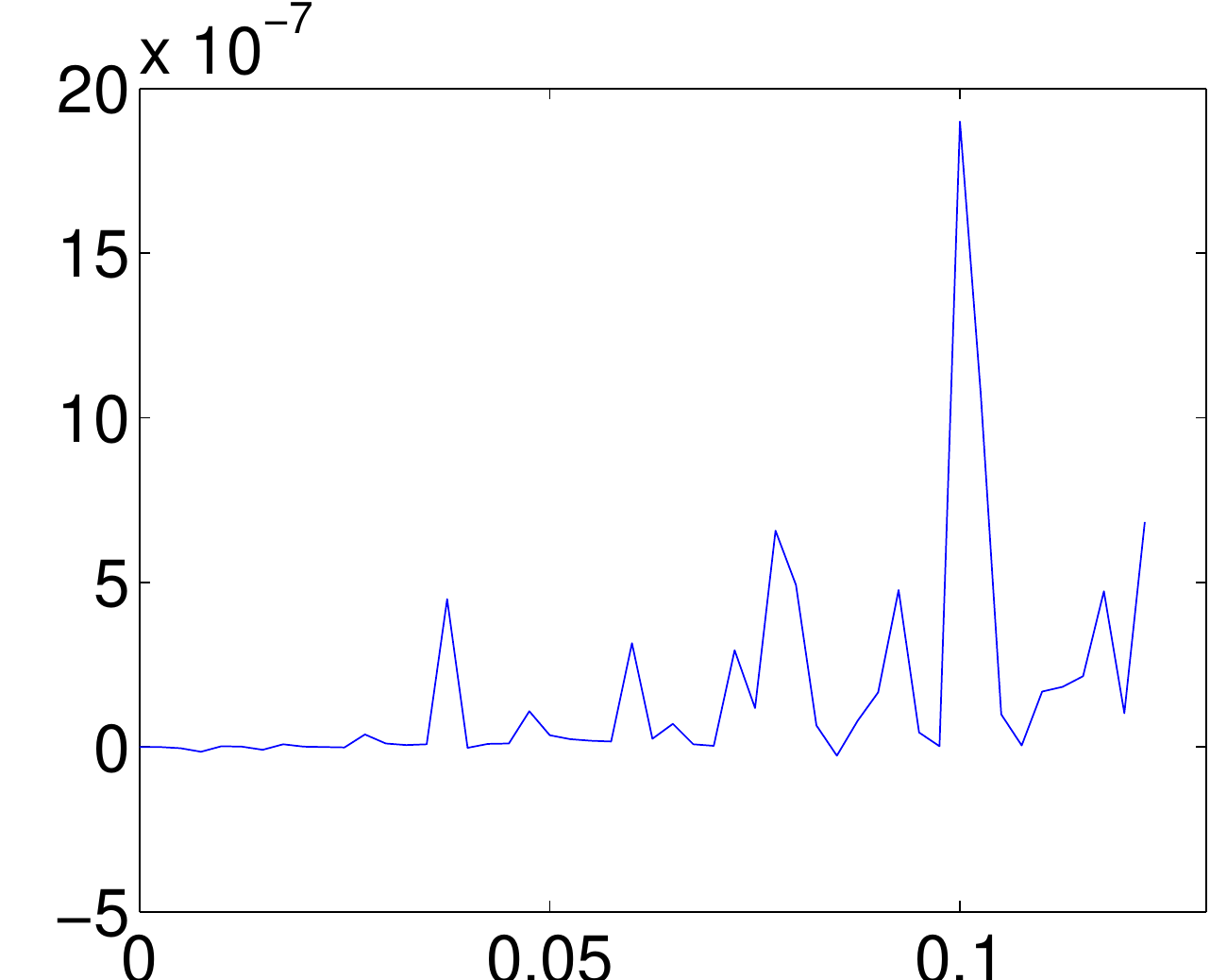} 
\includegraphics[angle=0,width=0.23\textwidth]{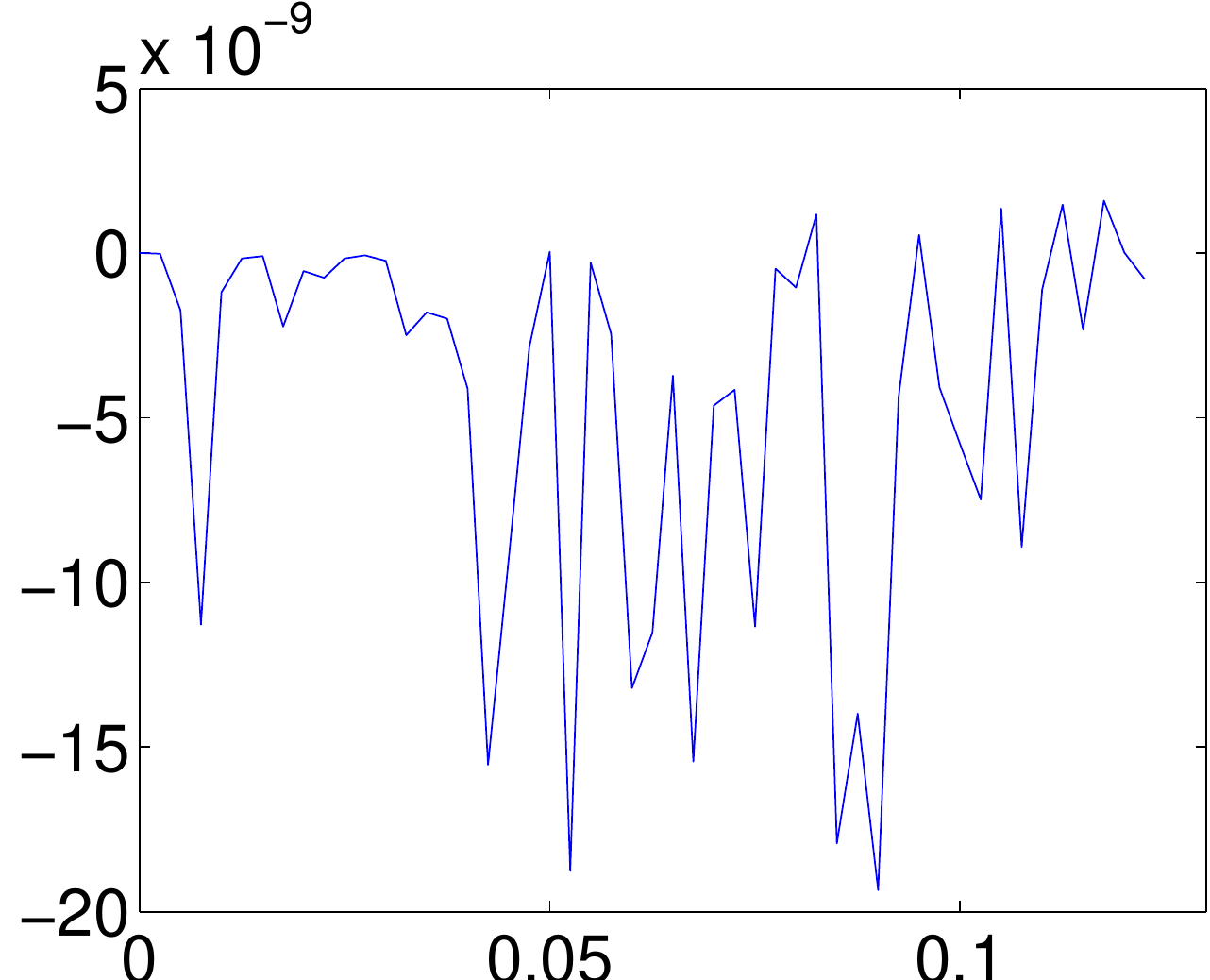}\\ 
\vspace*{-0.1cm}\hspace*{-0.0\textwidth} {$\tau$}\hspace*{0.22\textwidth} {$\tau$}\hspace*{0.22\textwidth} {$\tau$}\hspace*{0.22\textwidth} {$\tau$} 
\caption{Difference between numerically computed (with superscript $num$) and true (with superscript $true$) minimizing graphons on the line segment $(\E,\T)=\{0.5\}\times (0,0.5^3)$. From left to right: $c^{true}-c^{num}$, $g_{11}^{true}-g_{11}^{num}$, $g_{22}^{true}-g_{22}^{num}$ and $g_{12}^{true}-g_{12}^{num}$.}
\label{FIG:Benchmark-e=1/2}
\end{figure}

\paragraph{On the segment $(\E,\T)=\{0.5\}\times (0,0.5^3)$.} It is
again known from the theory in~\cite{RS2} that on this segment
the minimizing graphons are symmetric bi\partite  with
$g_{11}=g_{22}=0.5-(0.5^3-\T)$ and $g_{12}=g_{21}=0.5+(0.5^3-\T)$. In 
Fig.~\ref{FIG:Benchmark-e=1/2} we show the difference between the simulated minimizing graphons and the true minimizing graphons.
The last point on this segment is $(0.5,0.125)$. This is the point
where the global minimum of the rate function is achieved. Our
algorithms produce a minimizer that gives the value $I=-0.34657359$
which is less than $10^{-8}$ away from the true minimum value of $-(\ln 2)/2$.

\paragraph{Symmetric bi\partite graphons.} In Section \ref{SUBSEC:Phases} we find a formula for the optimizing graphons for phase II, namely the following symmetric bi\partite graphons:

\begin{equation}
g(x,y) = \begin{cases} \E - (\E^3-\T)^{1/3} & x,y < 1/2 \hbox{ or } x,y > 1/2 \cr
\E + (\E^3-\T)^{1/3} & x< \frac{1}{2} < y \hbox{ or } y < \frac{1}{2} < x
\end{cases}
\end{equation}

Our computational algorithms also find these minimizing symmetric bi\partite graphons; see discussions in the next section.

\subsection{Numerical Experiments}
\label{SEC:Num}
\begin{figure}[ht]
\centering
\rotatebox{90}{\hspace*{4.1cm}{\Large $\tau$}}
\includegraphics[angle=0,width=0.6\textwidth]{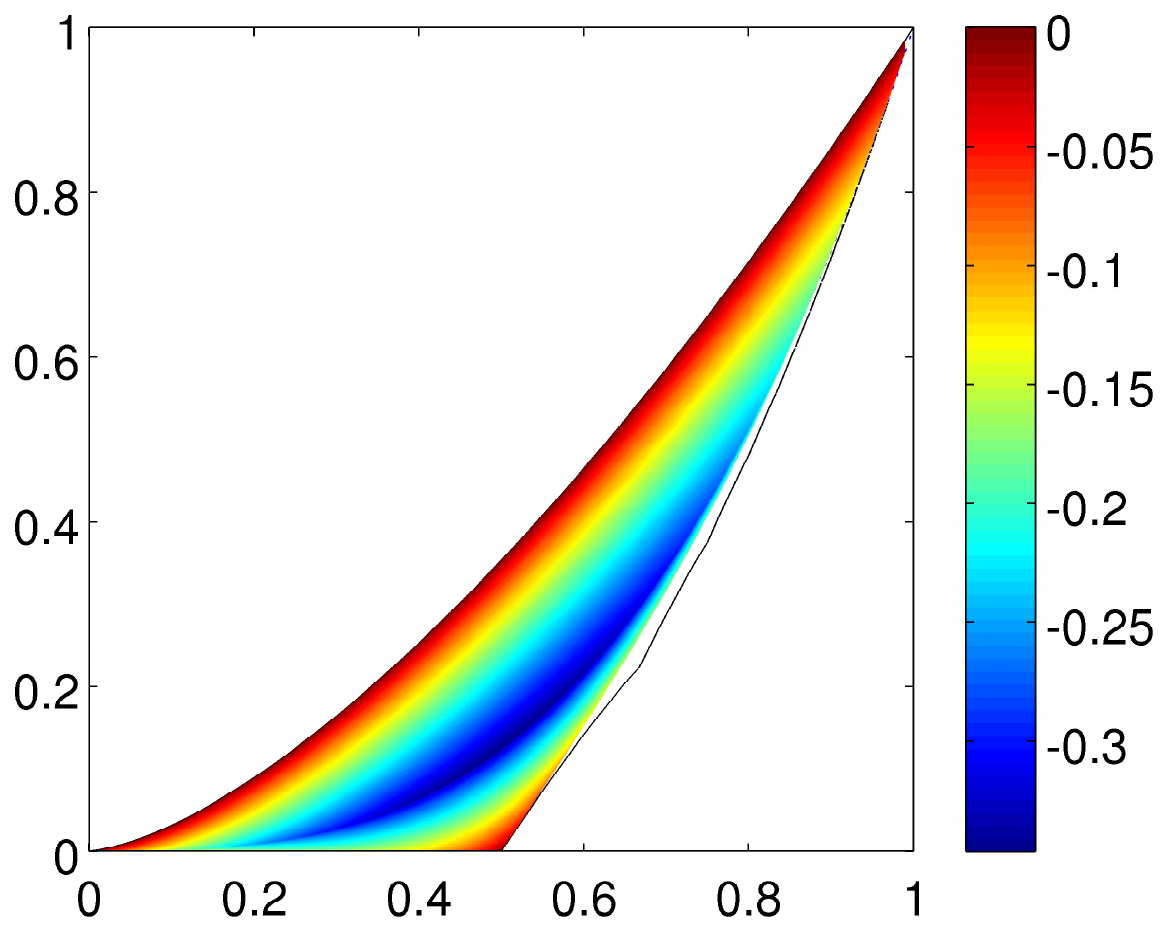}\\ 
\vspace*{-0.2cm}\hspace*{-1.1cm} {\Large $\epsilon$}
\caption{The minimal value of $I$ at different $(\E,\T)$ values.}
\label{FIG:Phase-I-Value}
\end{figure}
We now present numerical explorations of the minimizing graphons in some subregions of the phase space. We focus on two main subregions. The first is the subregion above the ER curve and below the upper boundary $\T=\E^{3/2}$. This is the phase I in Fig.~\ref{FIG:Phase-Boundary}. The second subregion is the region below the ER curve and above the lower boundary of the region where bi\partite graphons exist. We further split this region into two separate phases, II and III; see Fig.~\ref{FIG:Phase-Boundary}. Phase II is the region where the minimizing graphons are symmetric bi\partite while in phase III asymmetric bi\partite graphons are the minimizers.
\begin{figure}[ht]
\centering
\includegraphics[angle=0,width=0.19\textwidth]{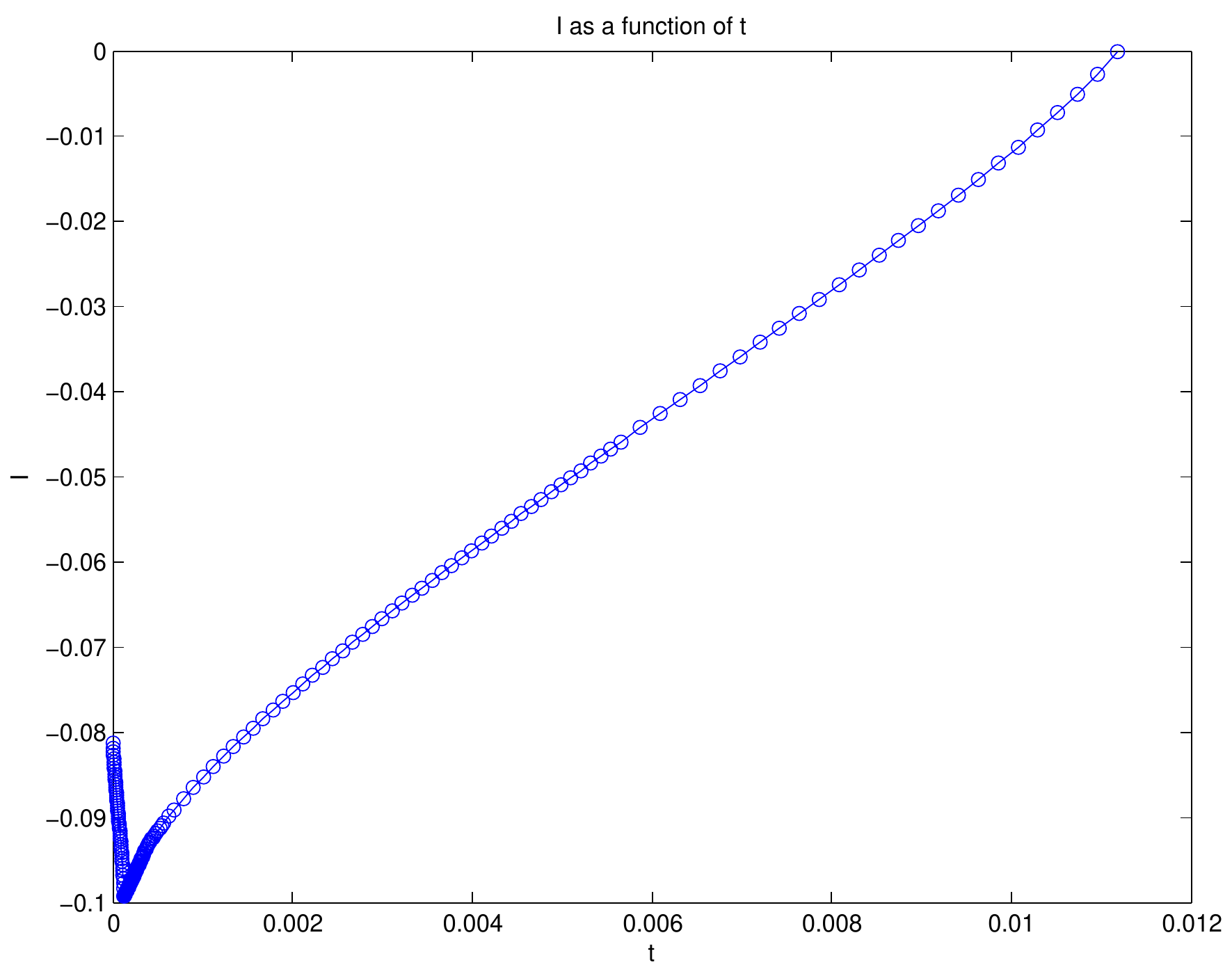} 
\includegraphics[angle=0,width=0.19\textwidth]{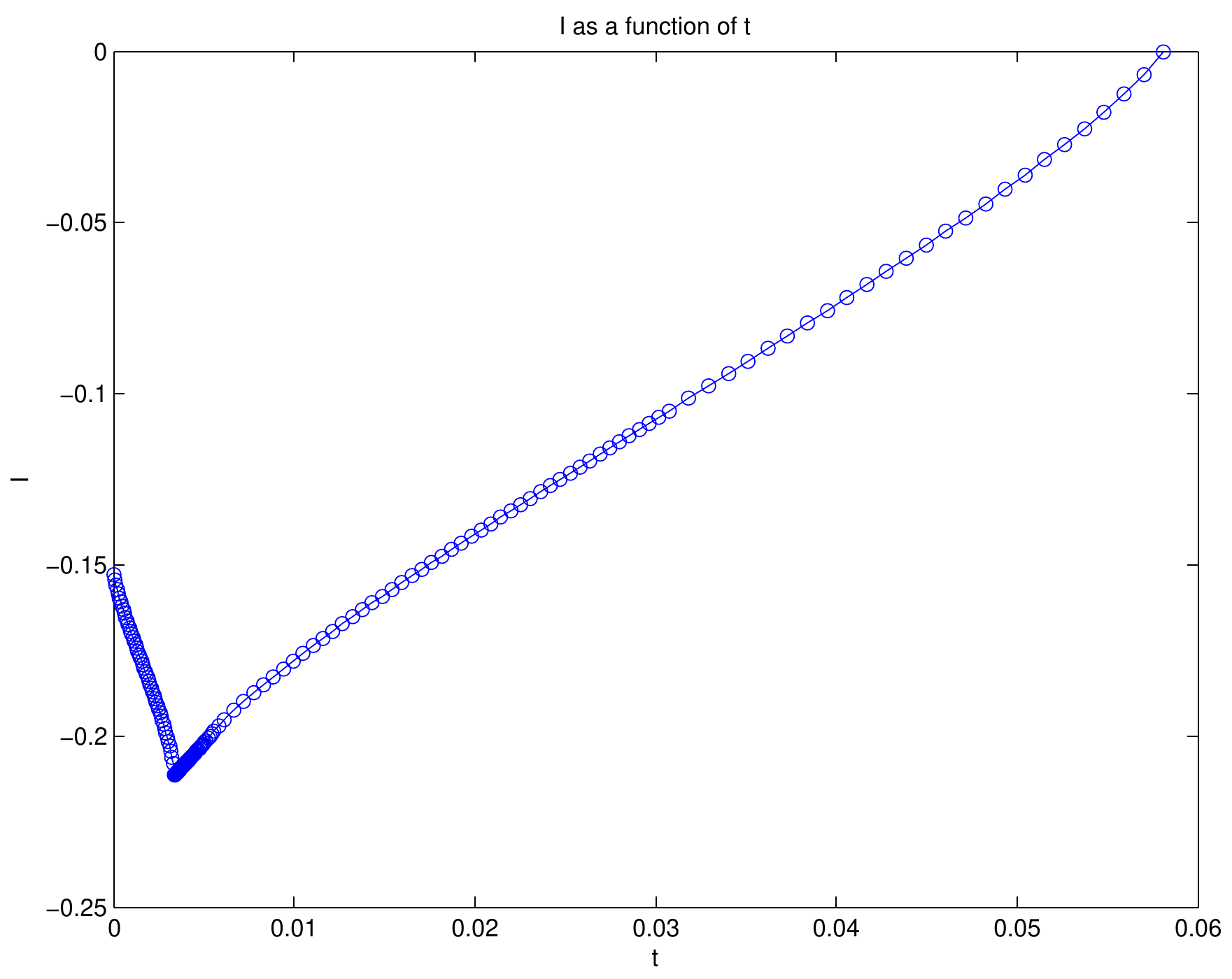} 
\includegraphics[angle=0,width=0.19\textwidth]{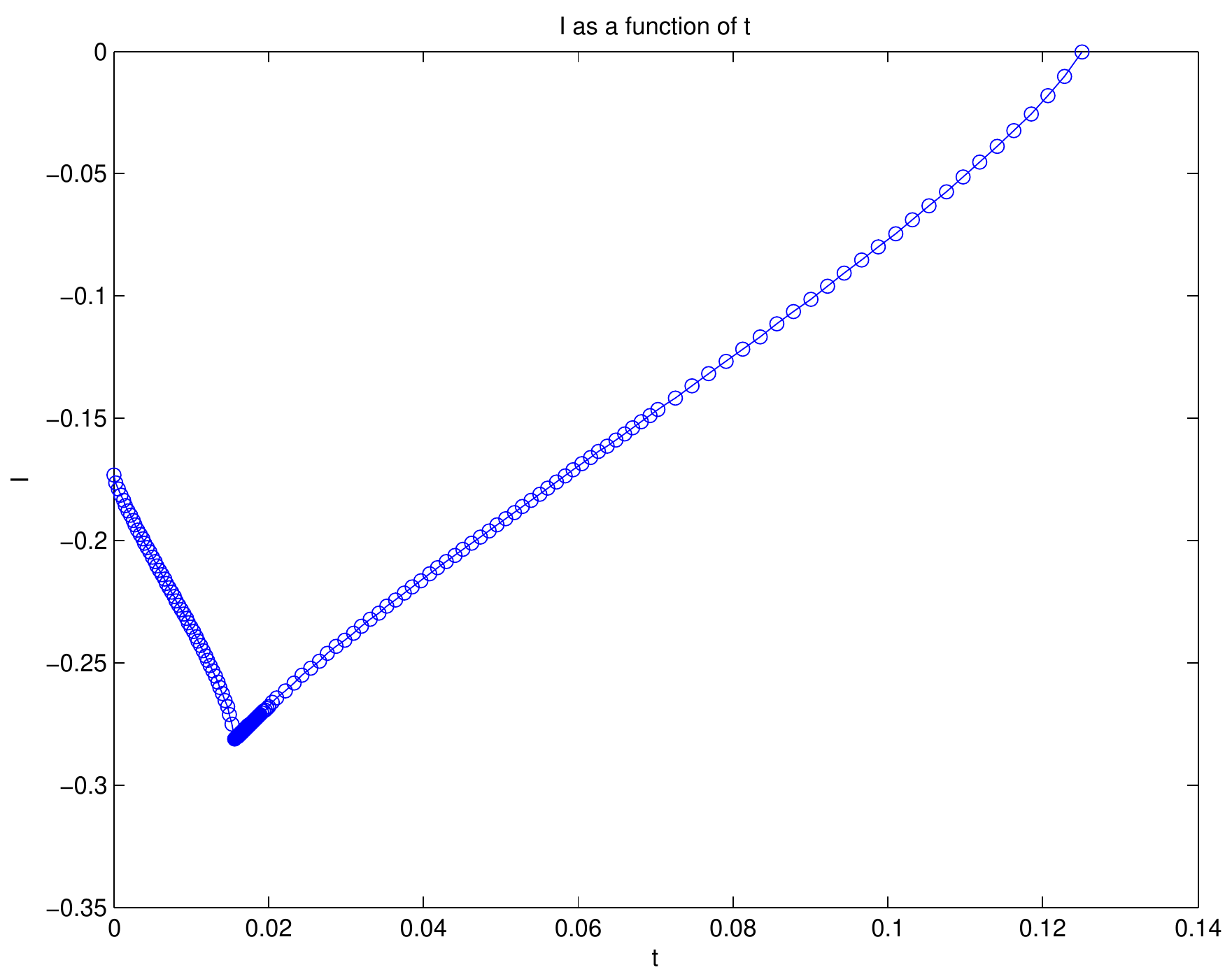} 
\includegraphics[angle=0,width=0.19\textwidth]{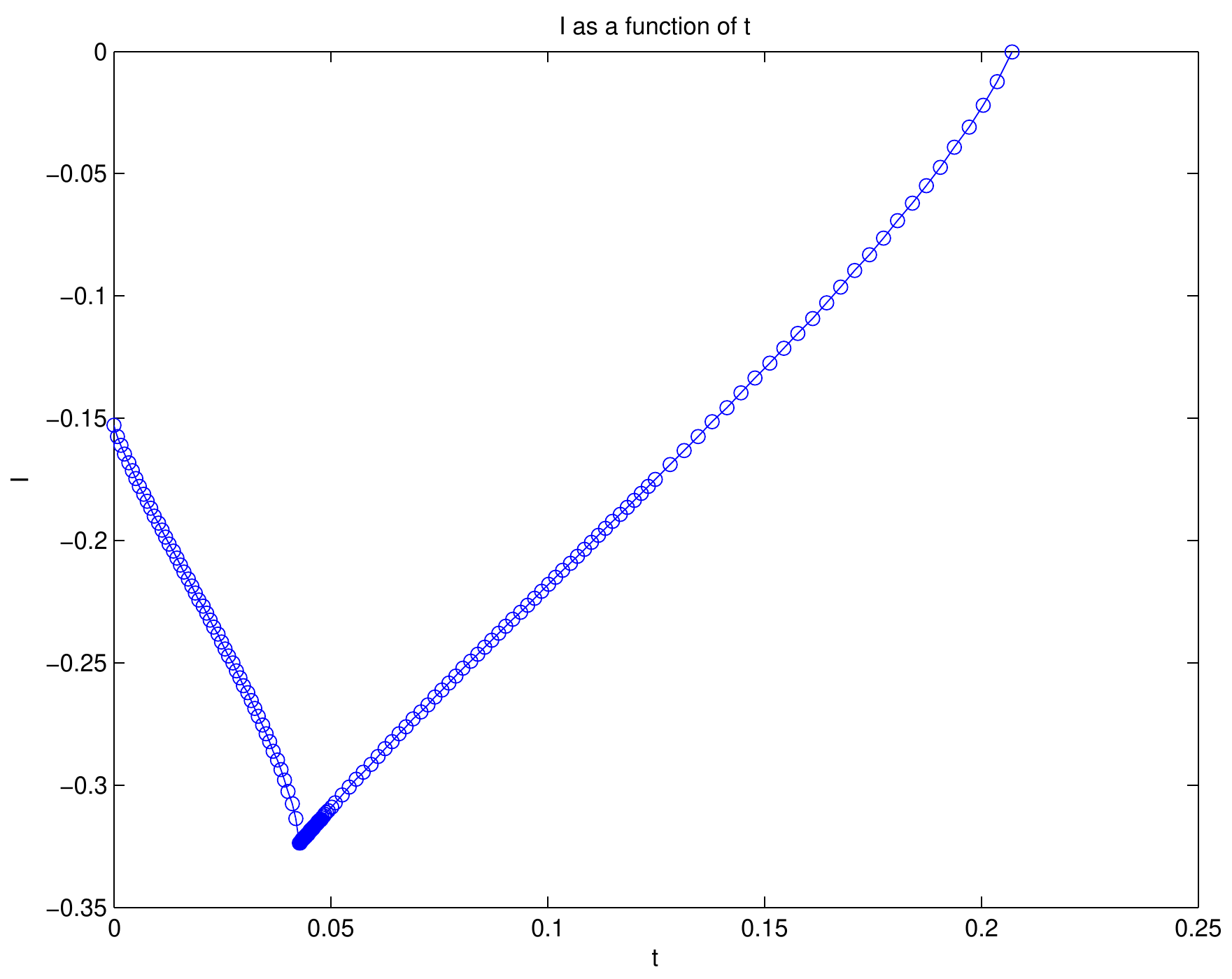} 
\includegraphics[angle=0,width=0.19\textwidth]{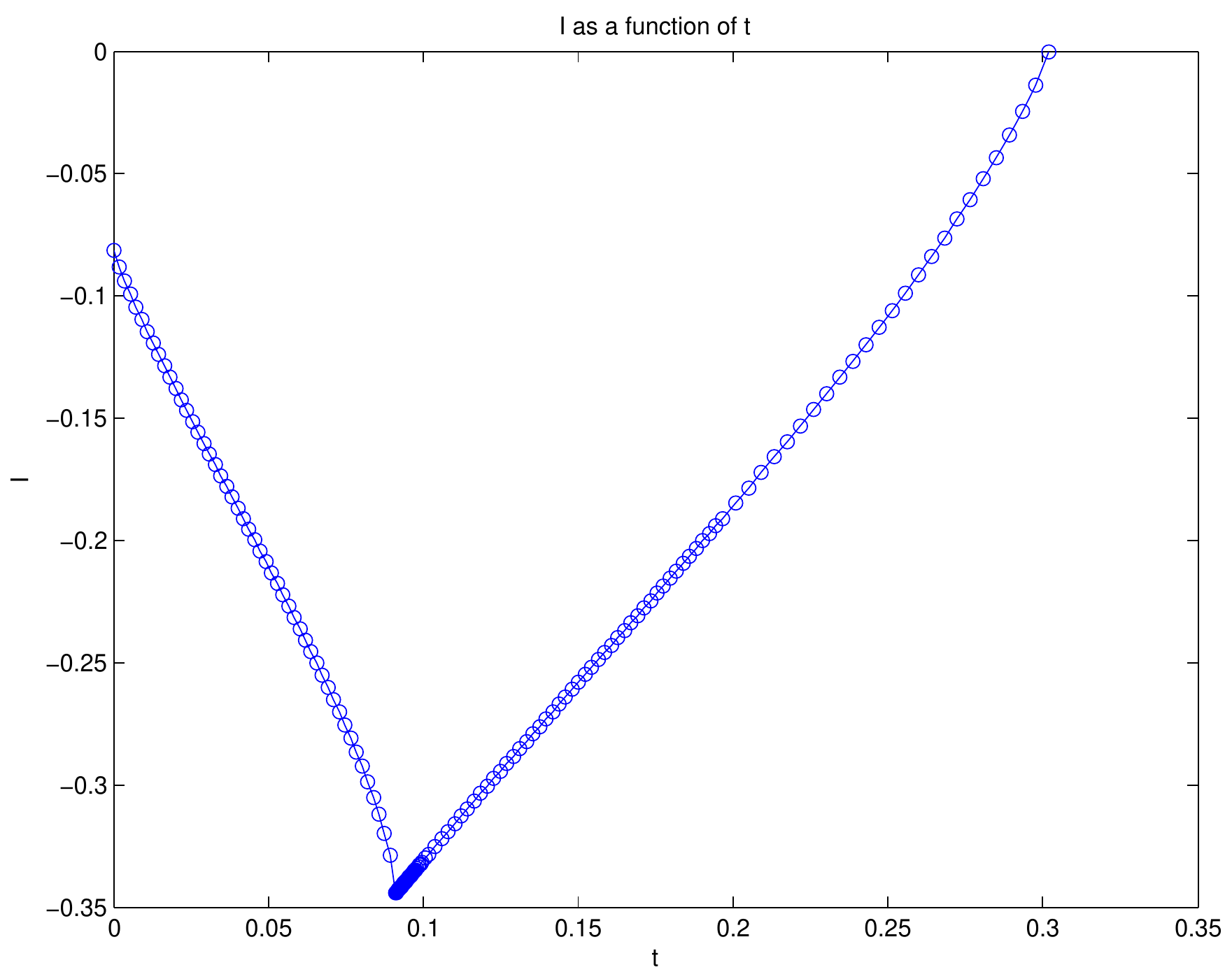}\\ 
\includegraphics[angle=0,width=0.19\textwidth]{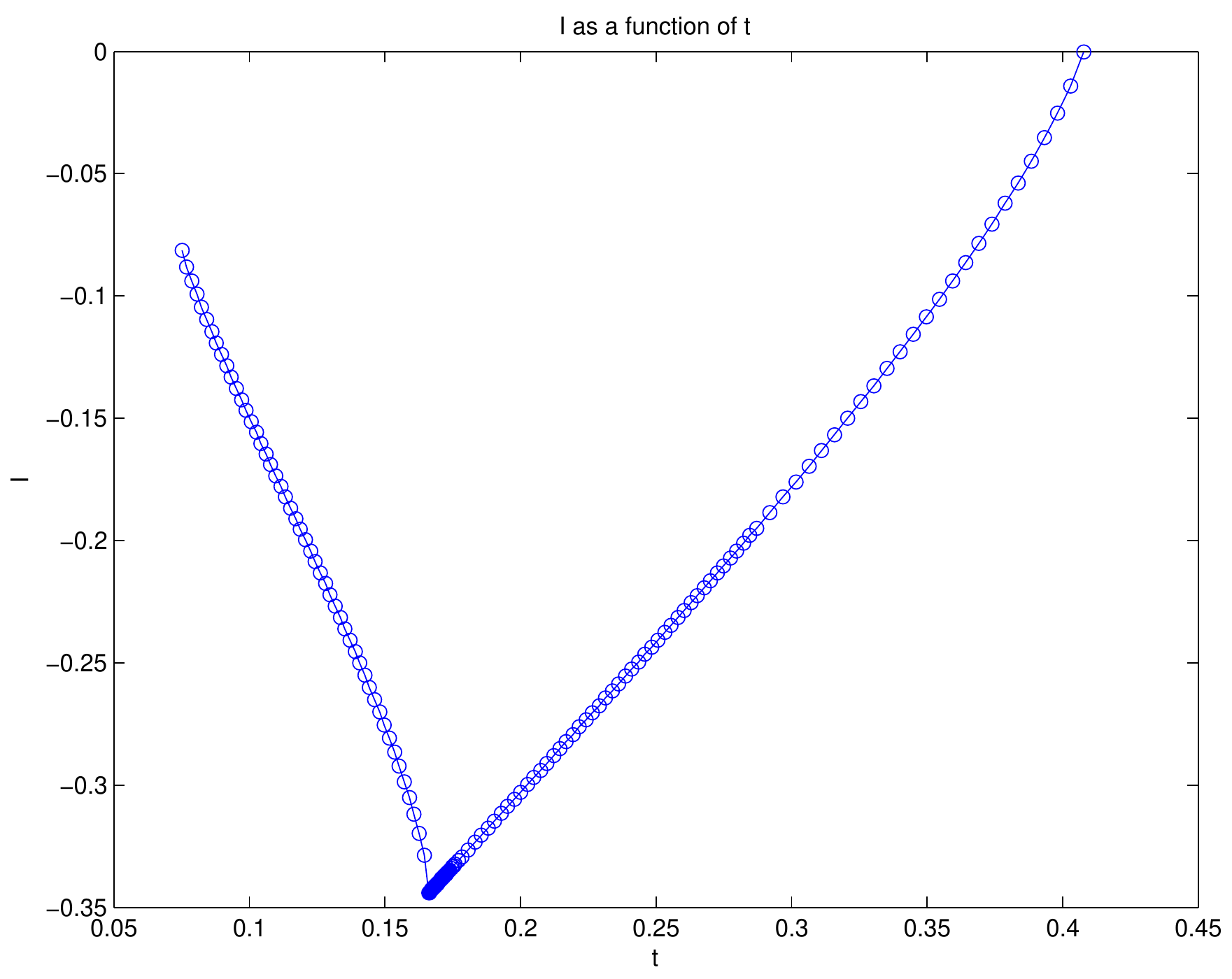} 
\includegraphics[angle=0,width=0.19\textwidth]{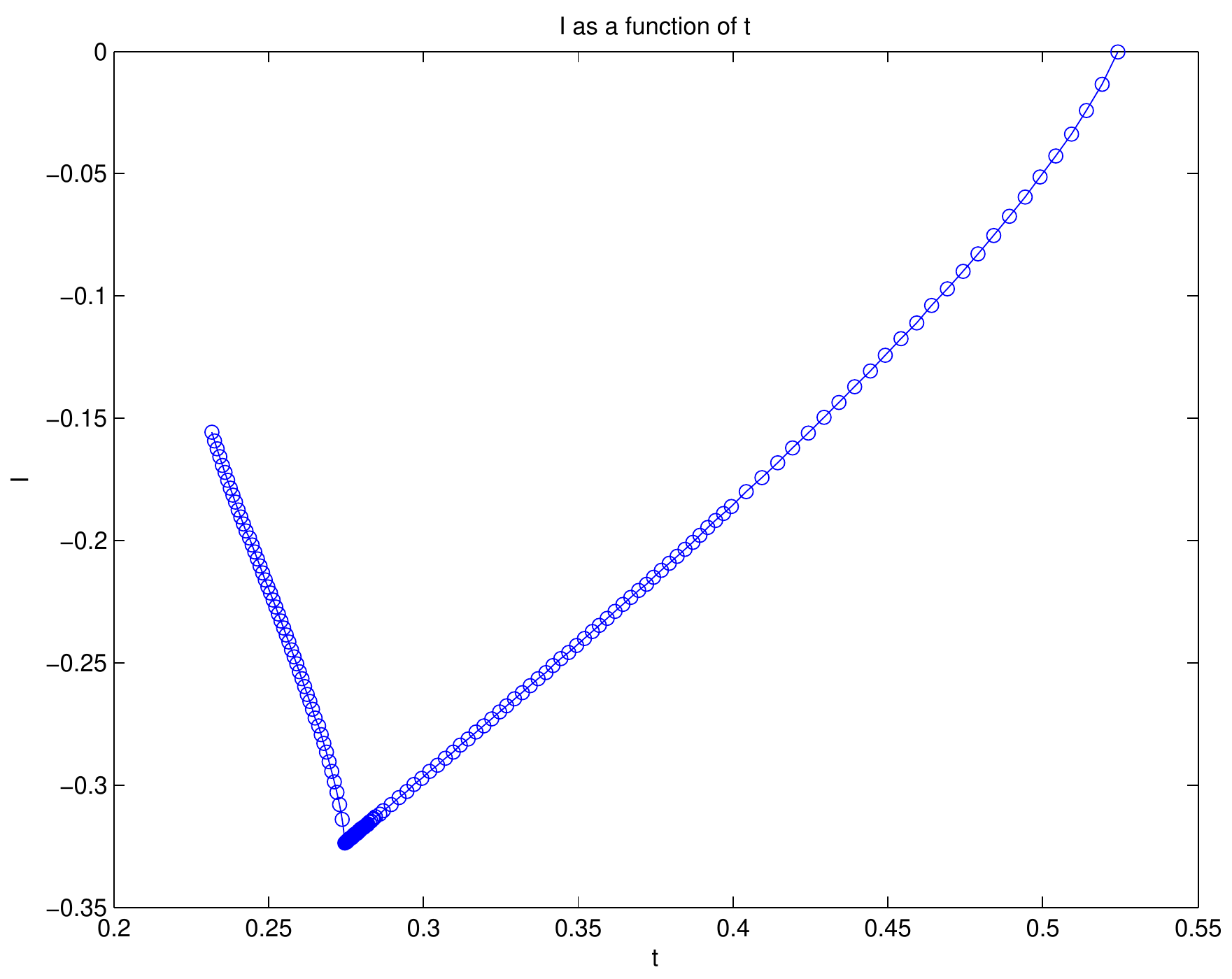} 
\includegraphics[angle=0,width=0.19\textwidth]{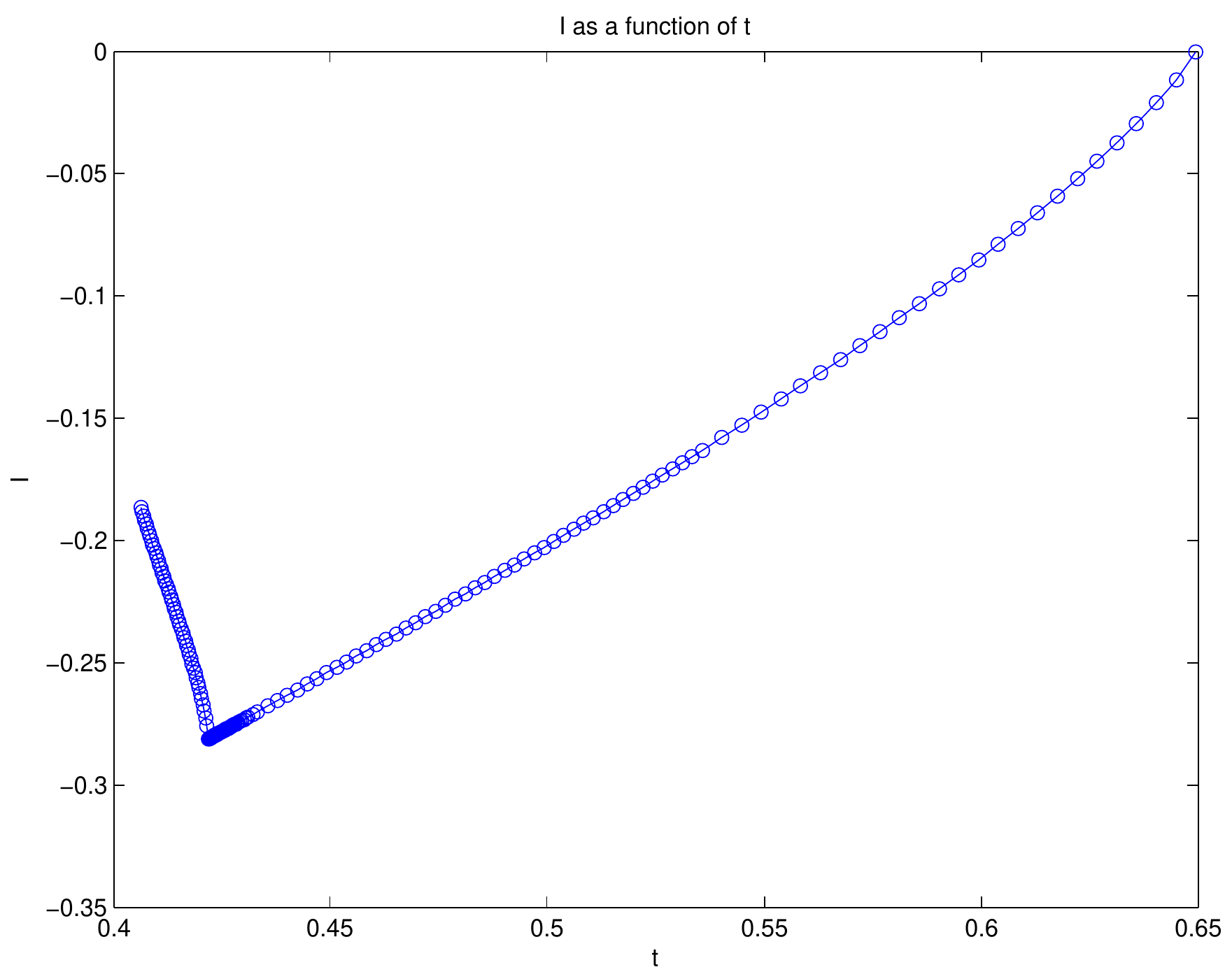} 
\includegraphics[angle=0,width=0.19\textwidth]{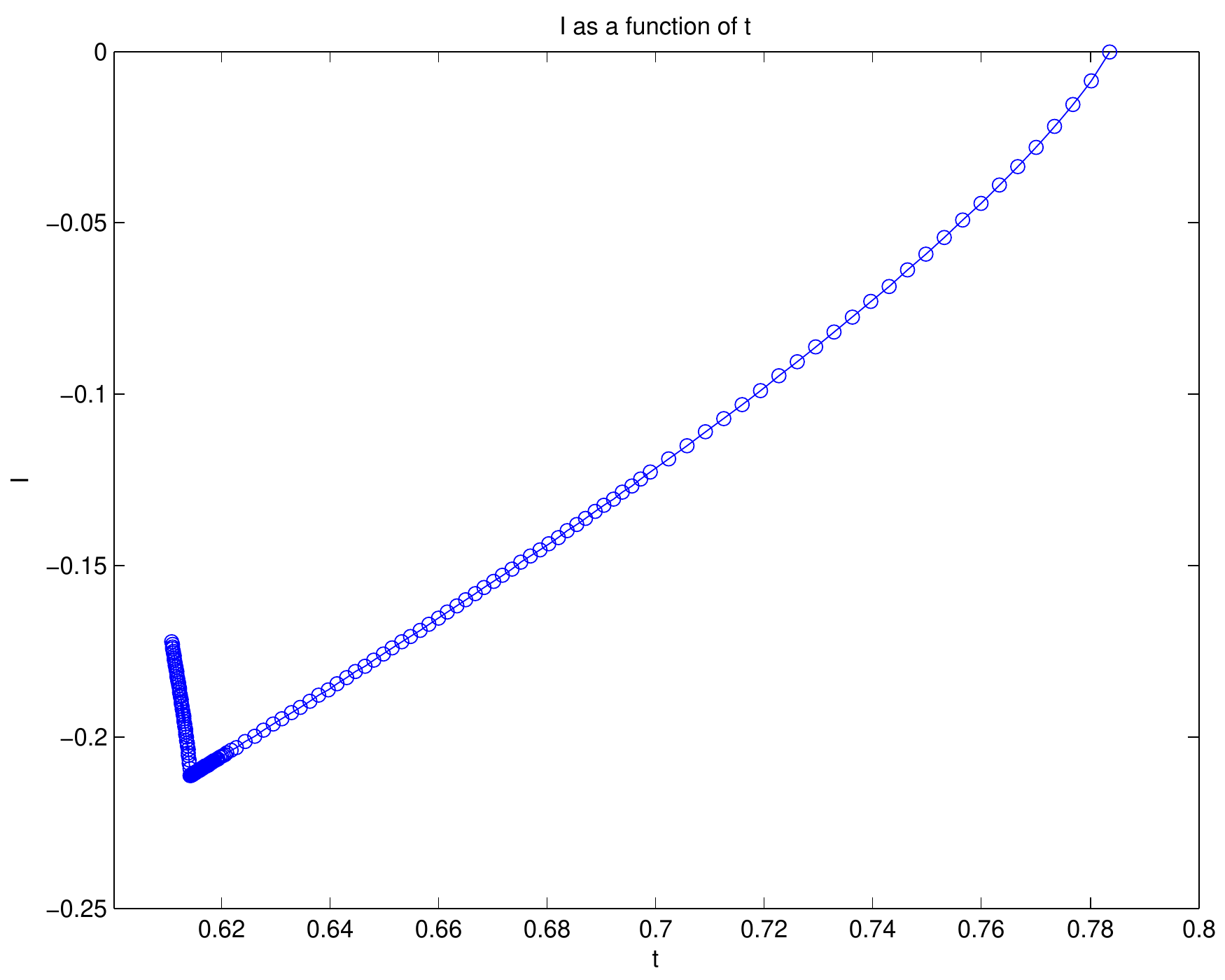} 
\includegraphics[angle=0,width=0.19\textwidth]{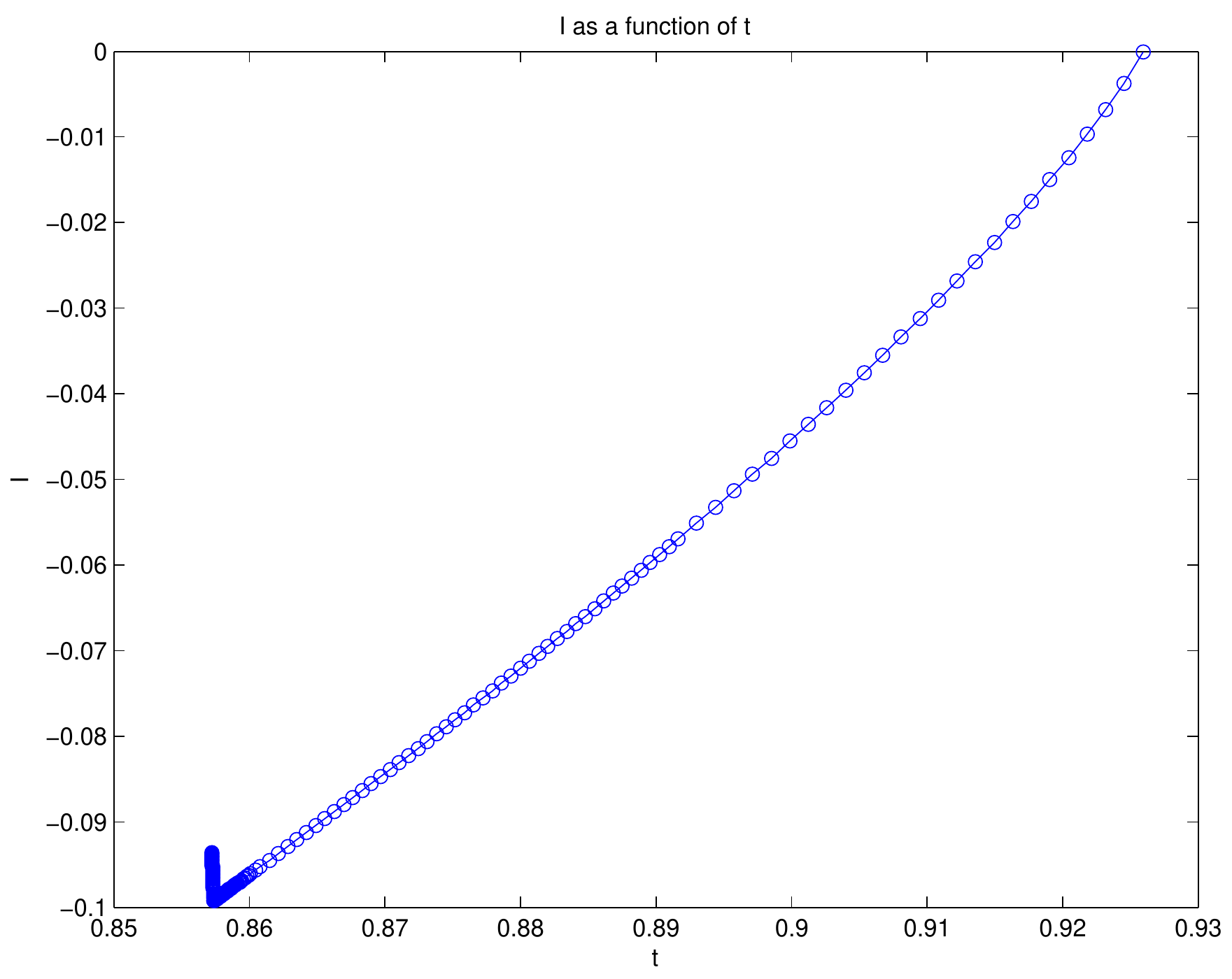} 
\caption{Cross-sections of the minimal value of the rate function $I$ as a function of $\T$ along lines $\E=a_k$ ($k=1,\cdots, 10$).}
\label{FIG:Phase-I-CS}
\end{figure}
We first show in Fig.~\ref{FIG:Phase-I-Value} the minimal values of
$I$ in the whole region in which we are interested: ${\rm I\cup II\cup
  III}$. The minimum of $I$ for fixed $\E$ is achieved on the ER curve. The cross-sections in Fig.~\ref{FIG:Phase-I-CS} along the lines $\E=a_k=(k-1)*0.1+0.05$ give a better visualization of the landscape of the $I$.

\begin{figure}[!ht]
\centering
\includegraphics[angle=0,width=0.23\textwidth]{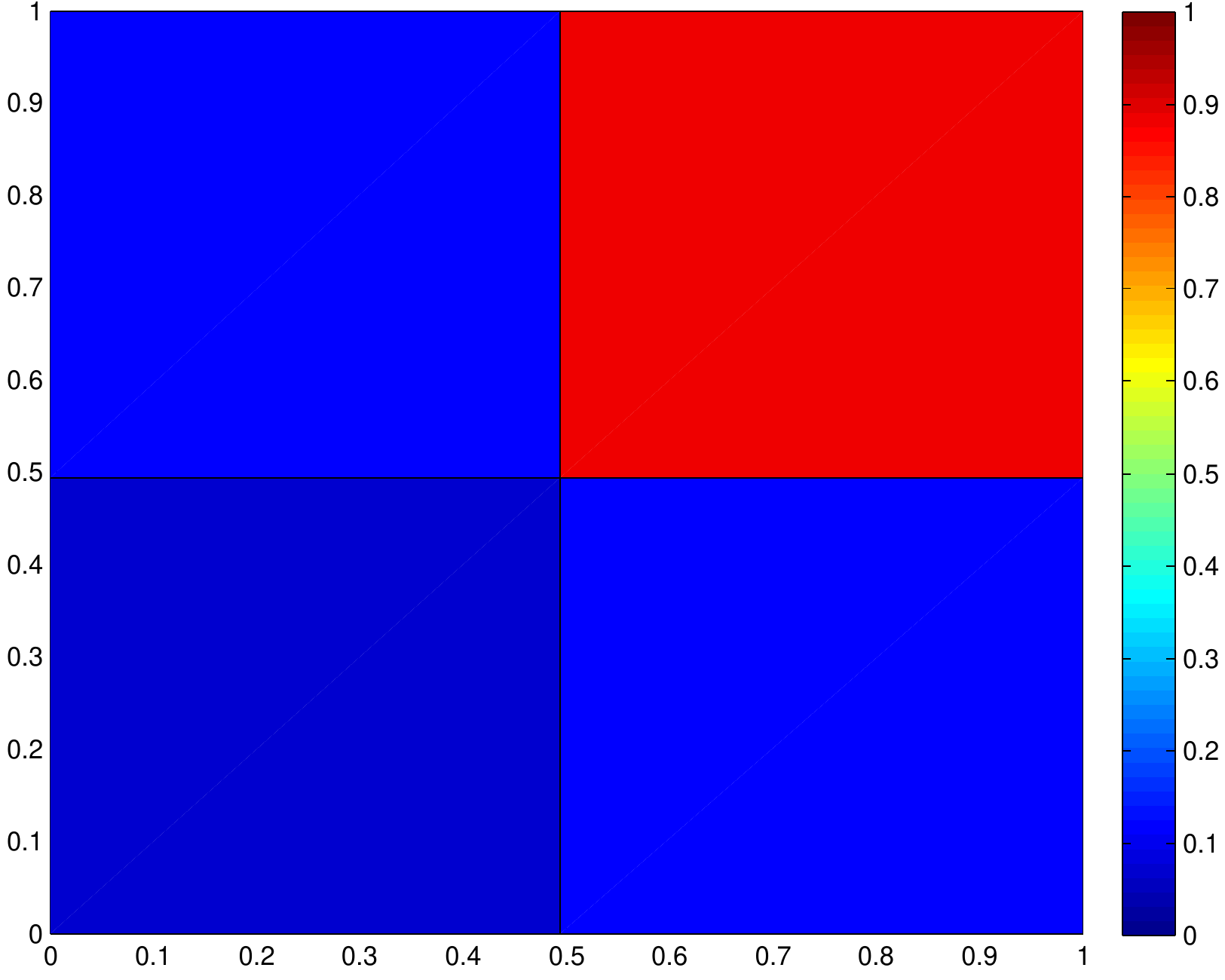} 
\includegraphics[angle=0,width=0.23\textwidth]{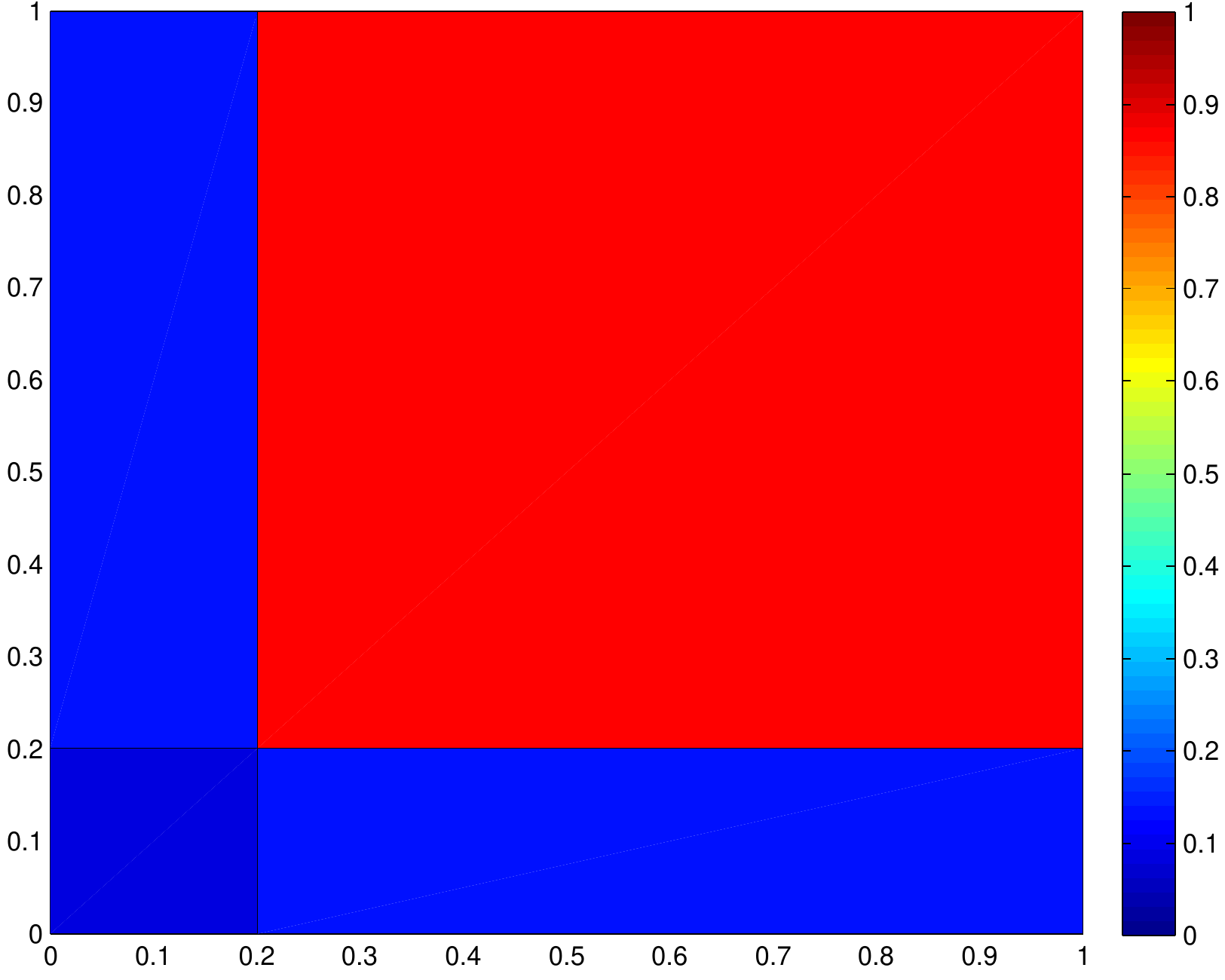} 
\includegraphics[angle=0,width=0.23\textwidth]{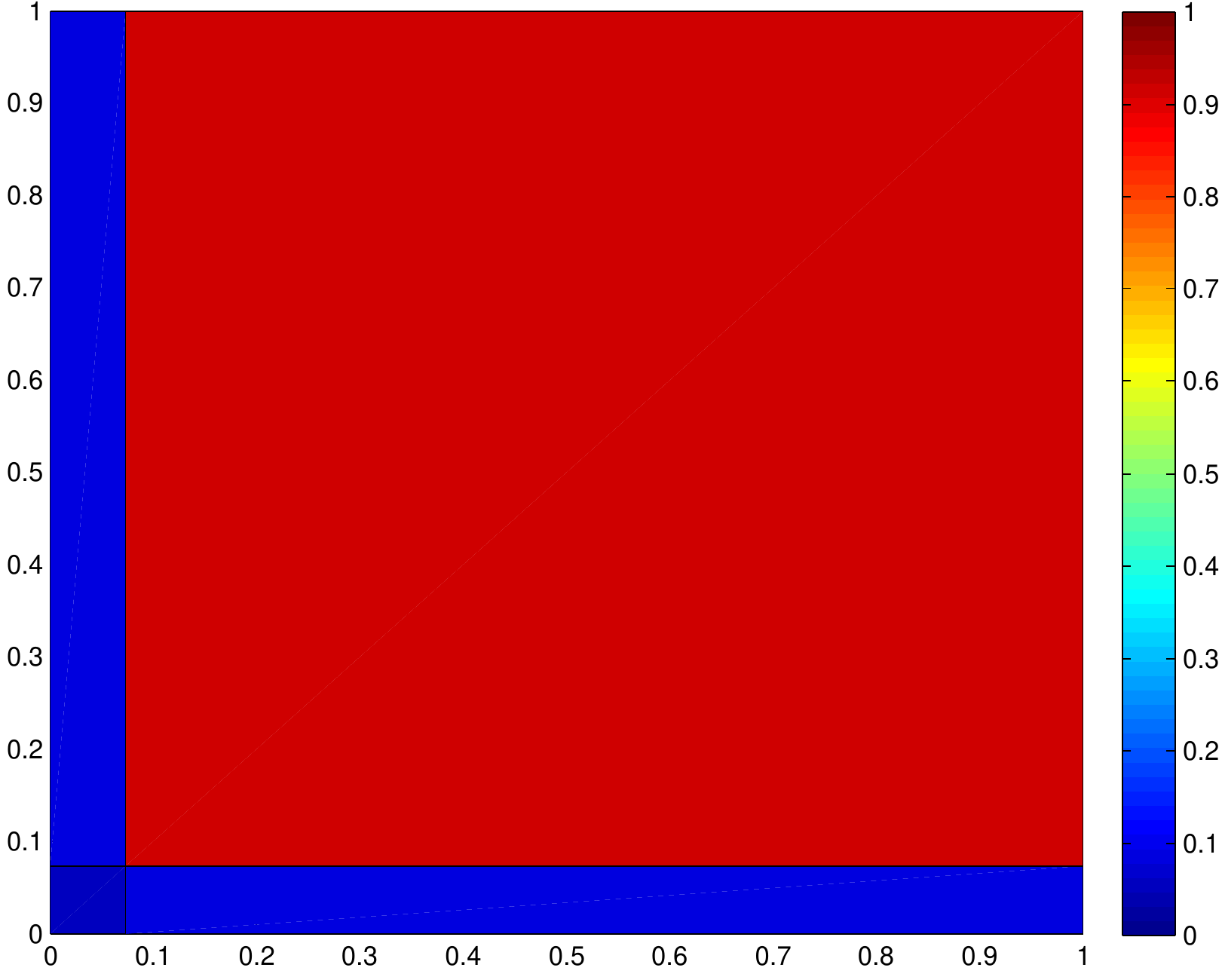}\\
\includegraphics[angle=0,width=0.23\textwidth]{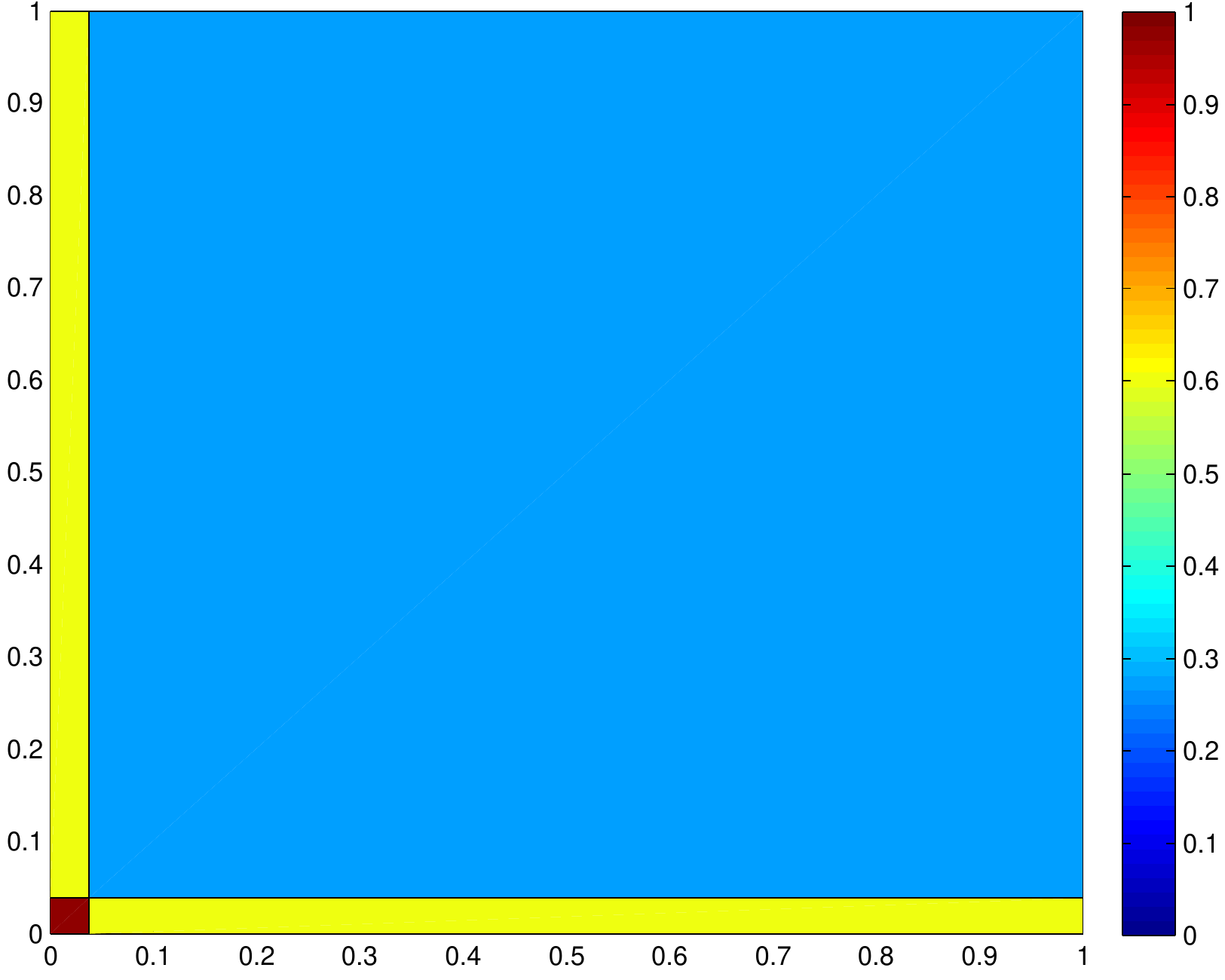} 
\includegraphics[angle=0,width=0.23\textwidth]{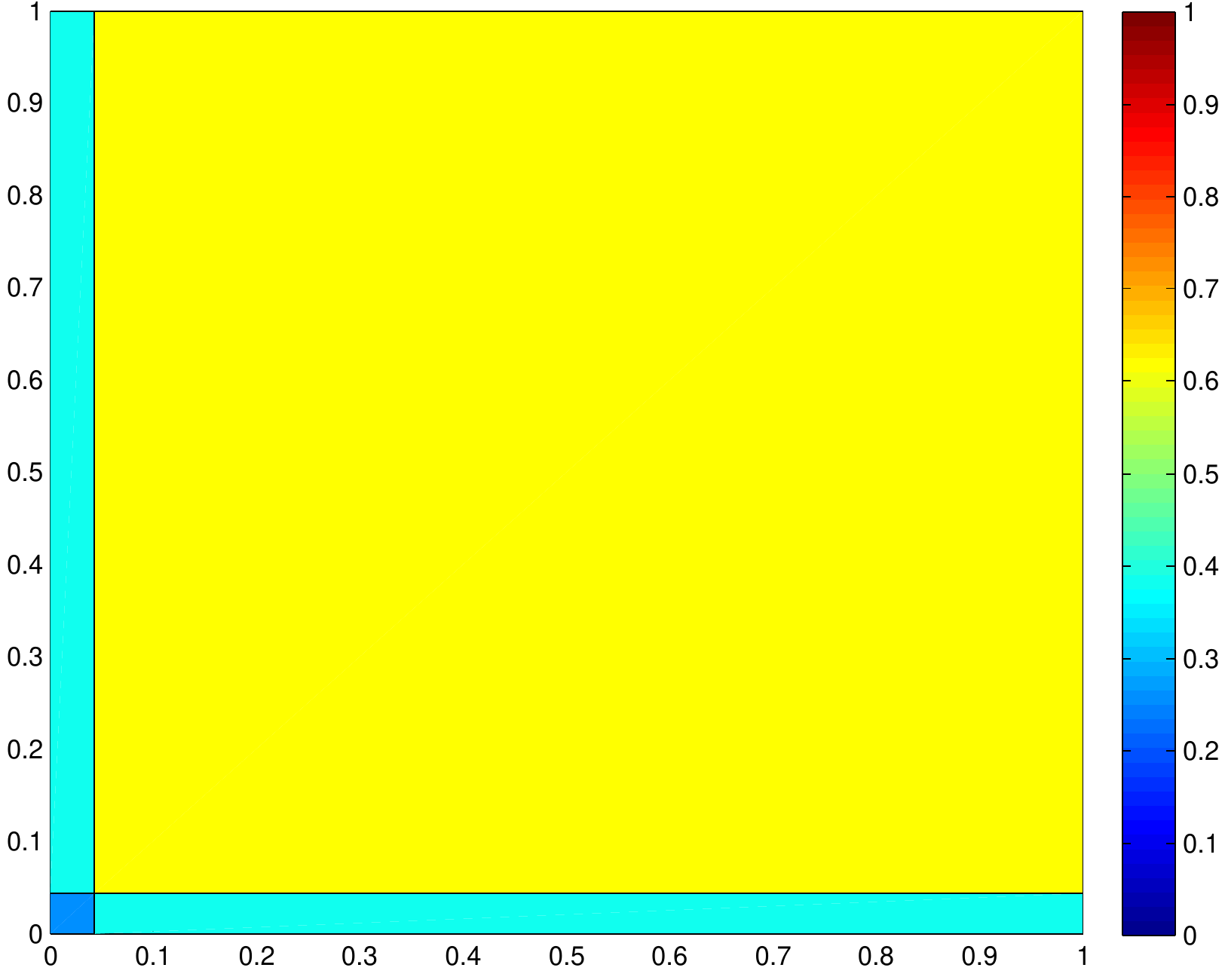} 
\includegraphics[angle=0,width=0.23\textwidth]{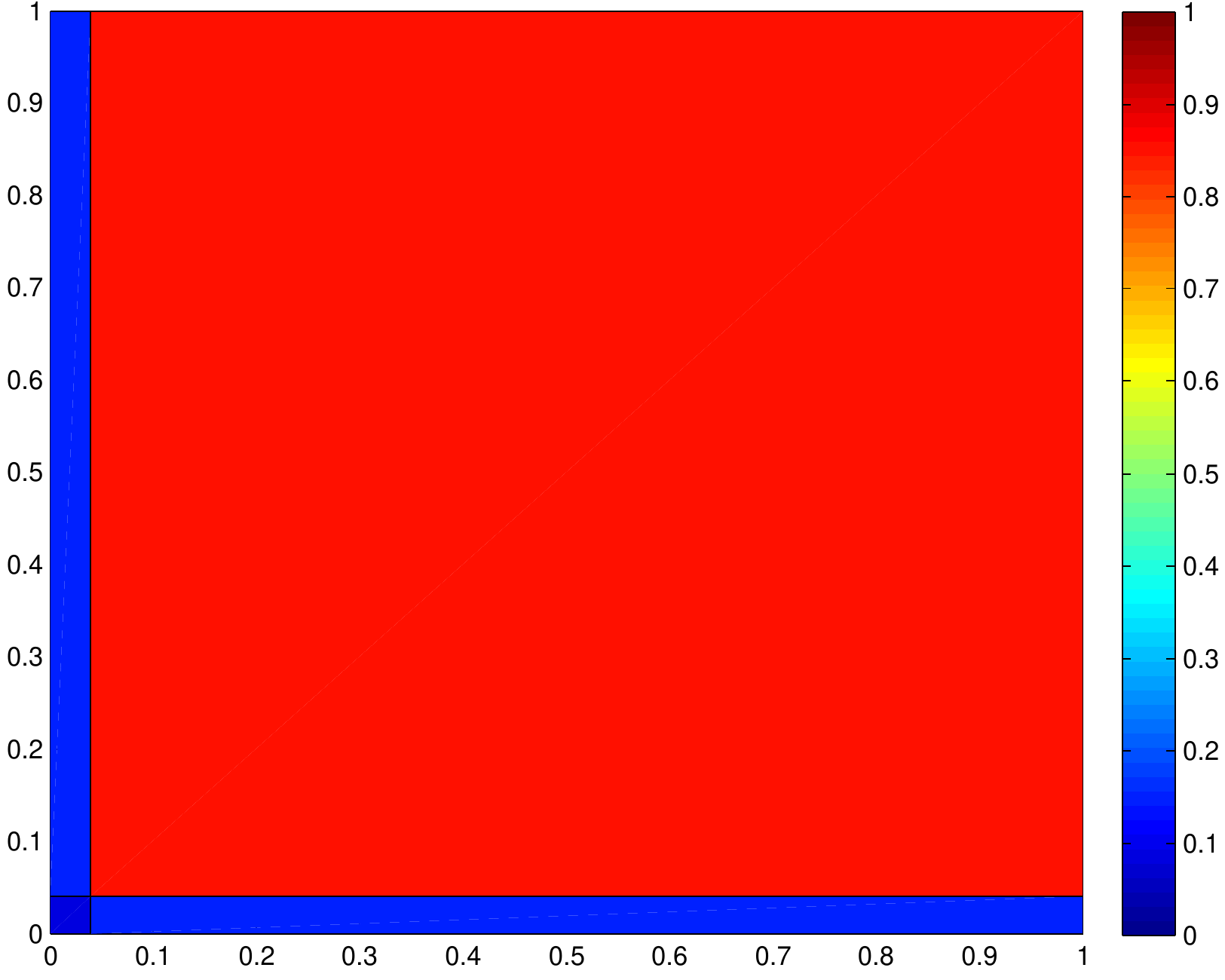}\\
\includegraphics[angle=0,width=0.23\textwidth]{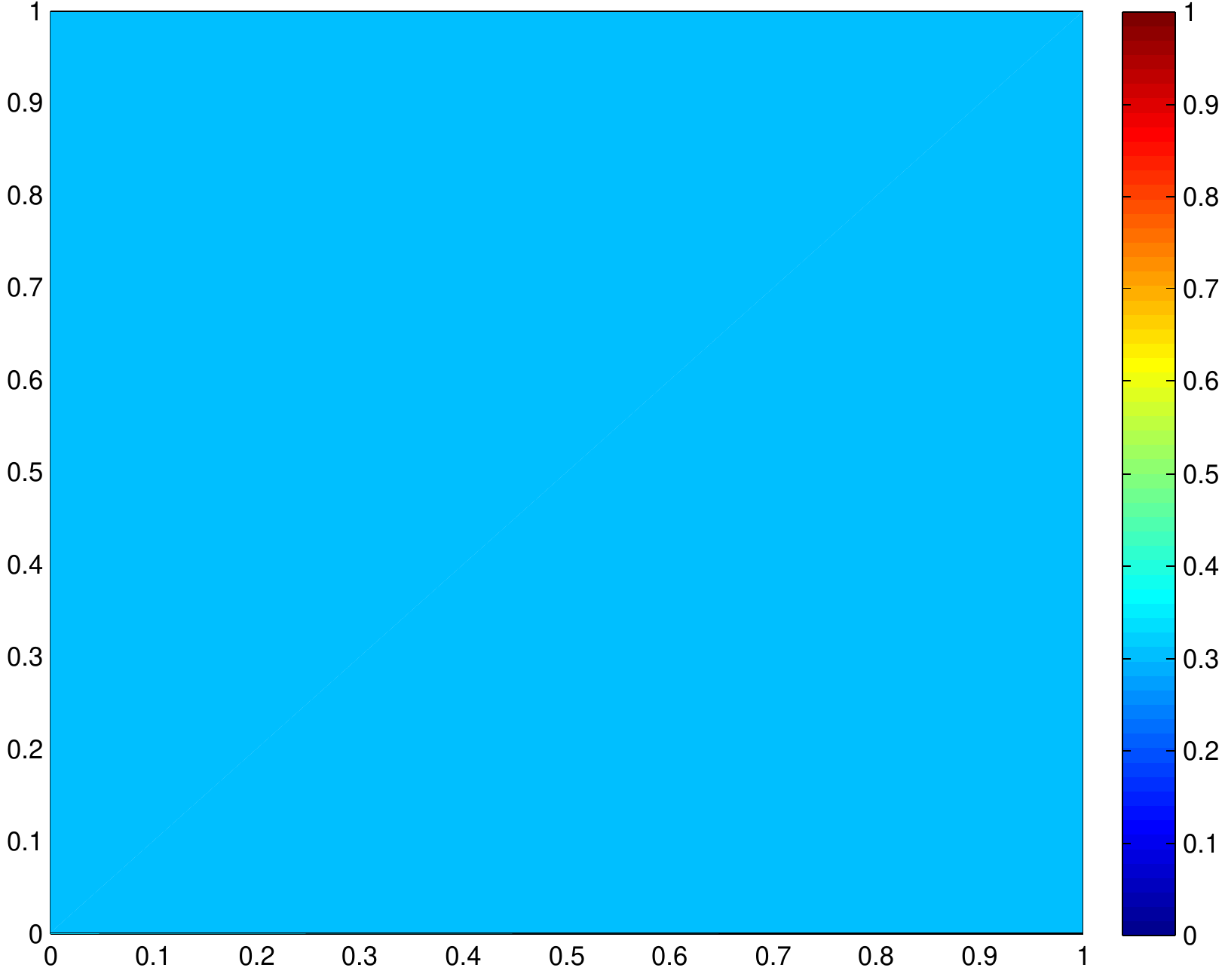} 
\includegraphics[angle=0,width=0.23\textwidth]{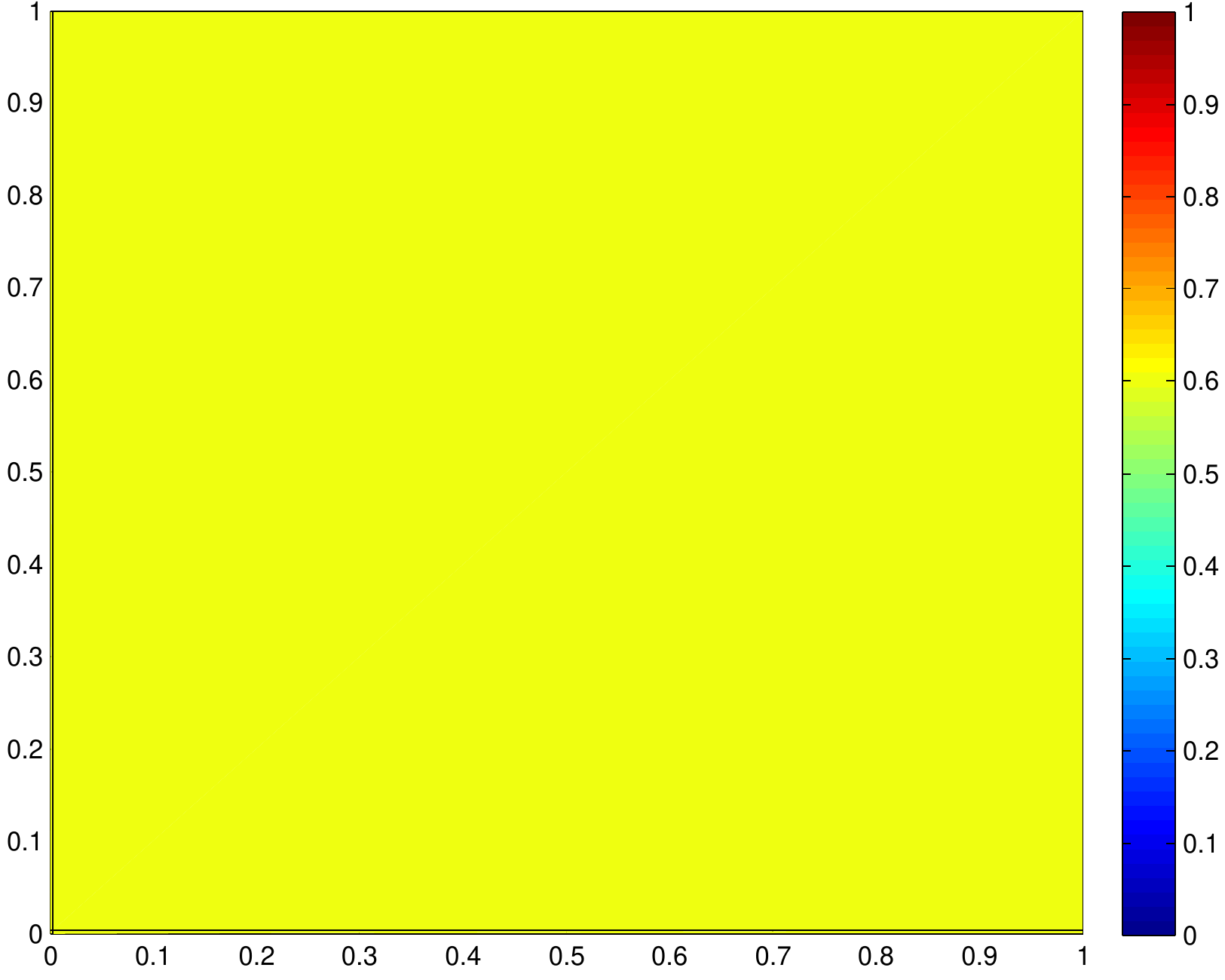} 
\includegraphics[angle=0,width=0.23\textwidth]{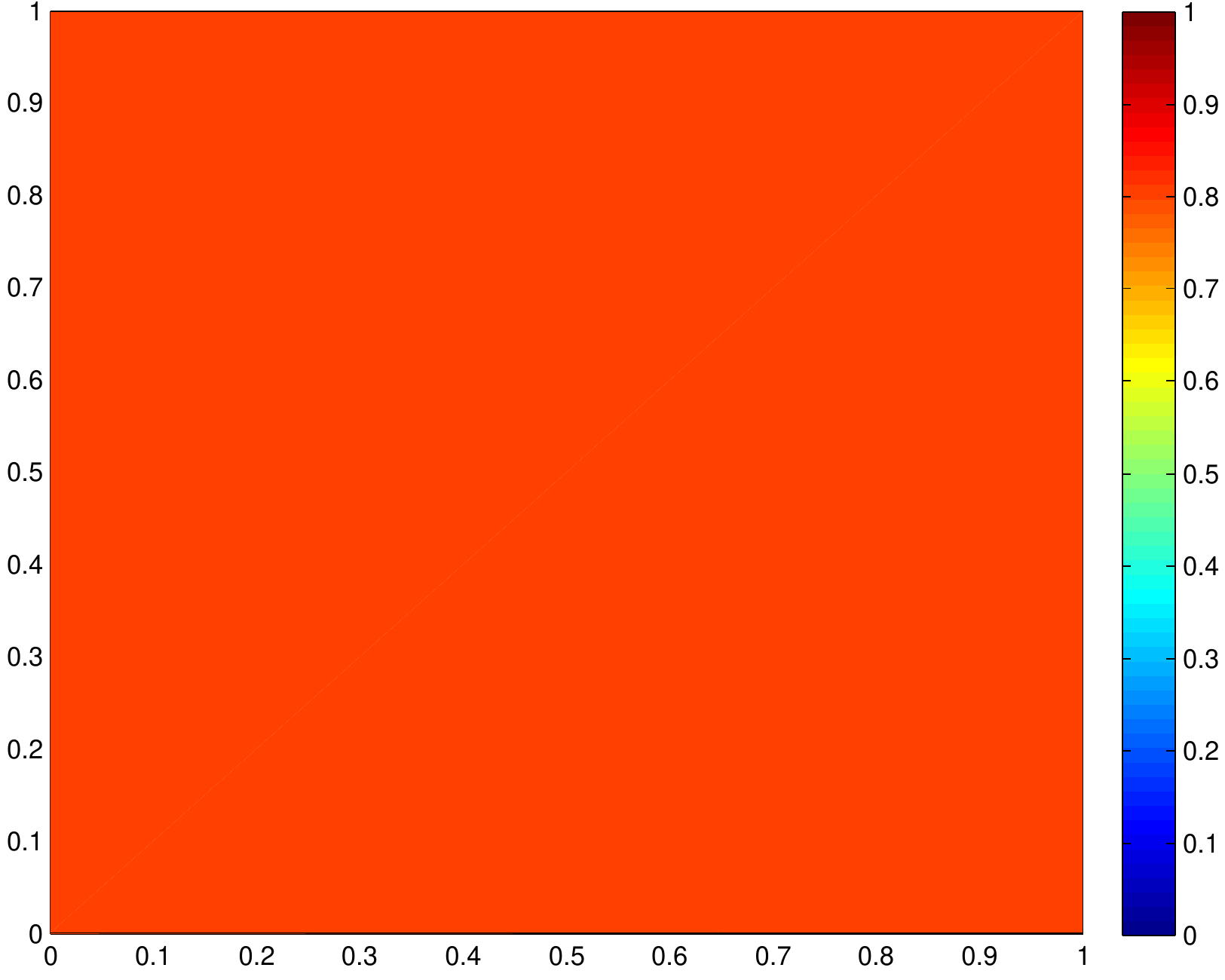}\\ 
\includegraphics[angle=0,width=0.23\textwidth]{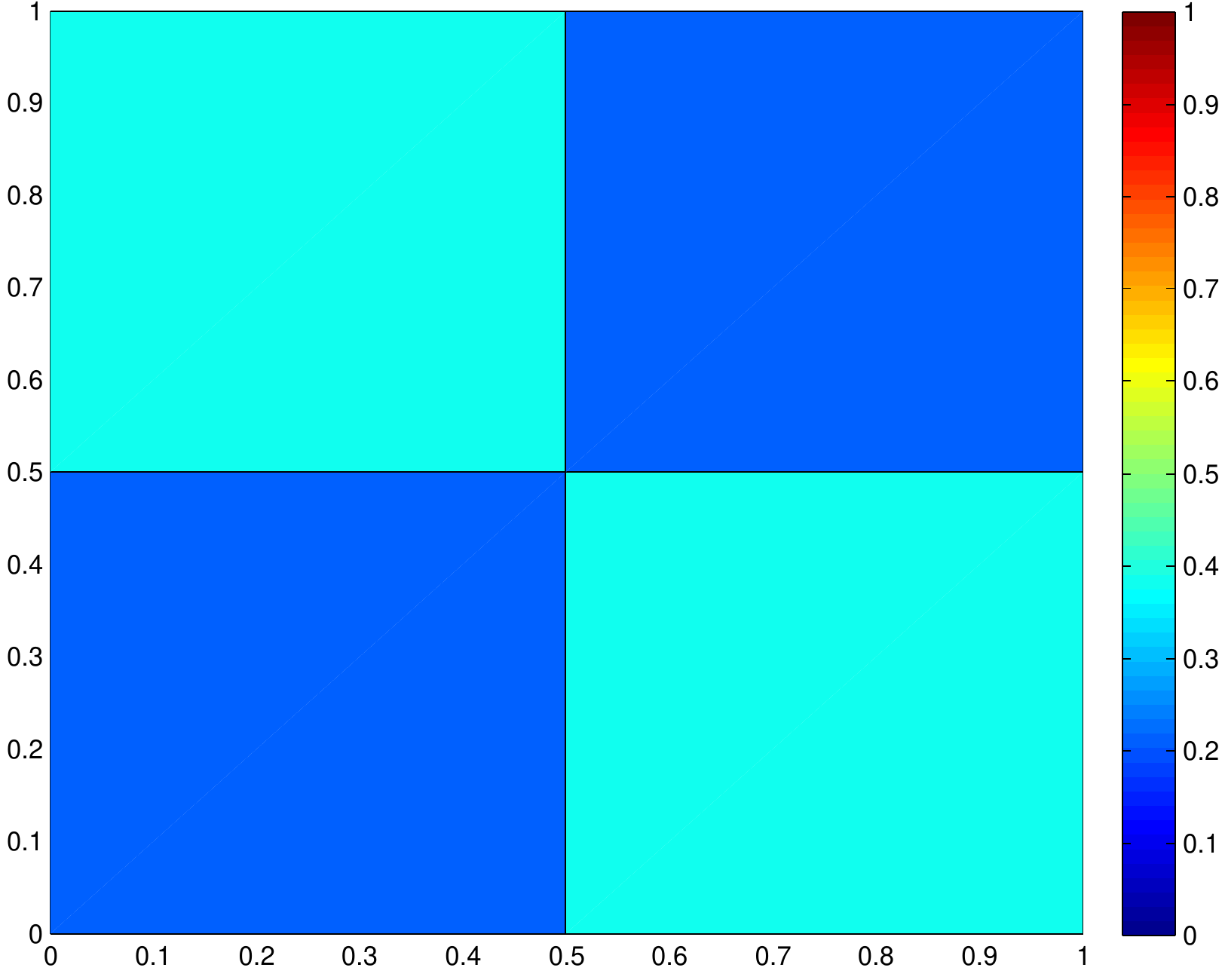} 
\includegraphics[angle=0,width=0.23\textwidth]{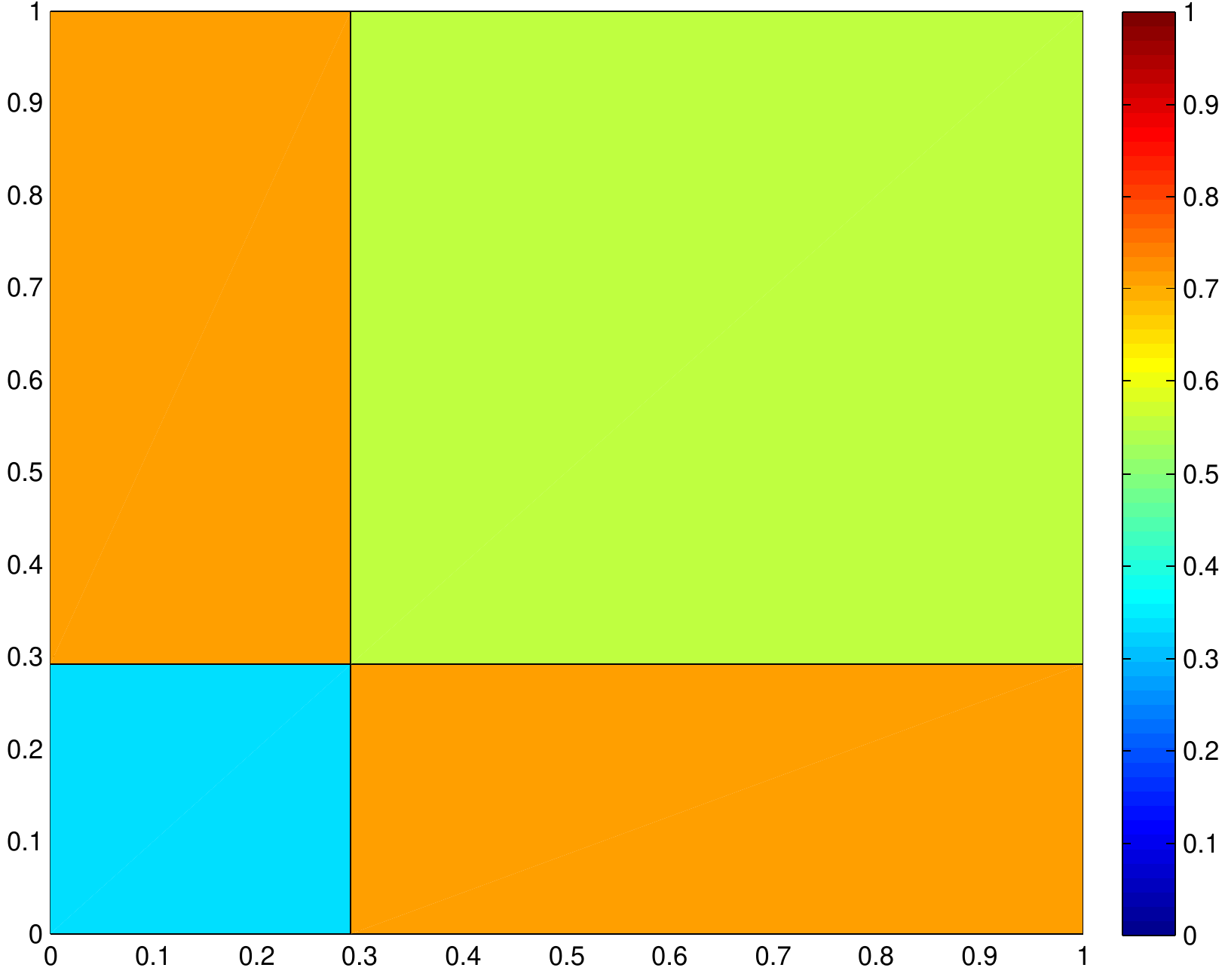} 
\includegraphics[angle=0,width=0.23\textwidth]{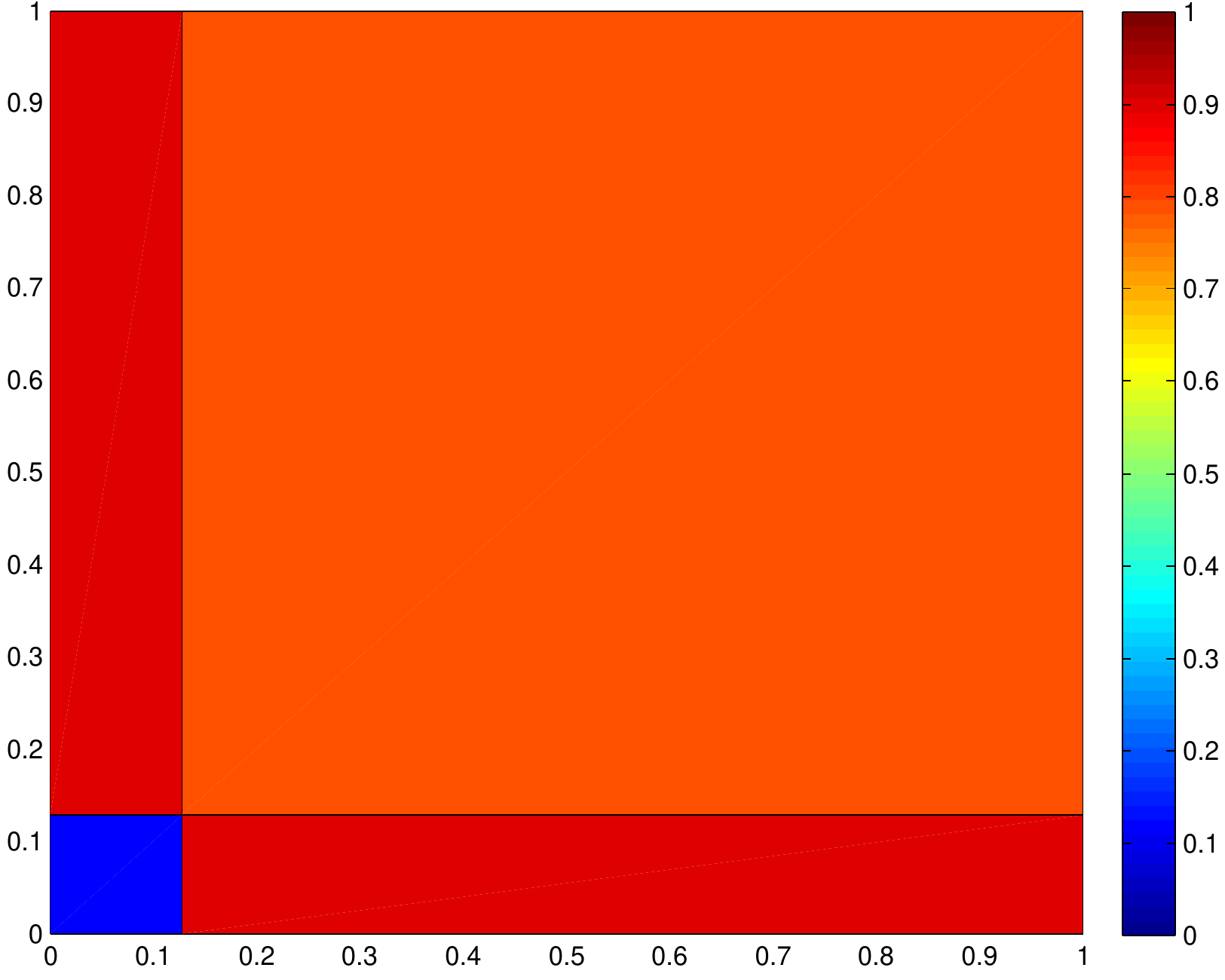}\\ 
\includegraphics[angle=0,width=0.23\textwidth]{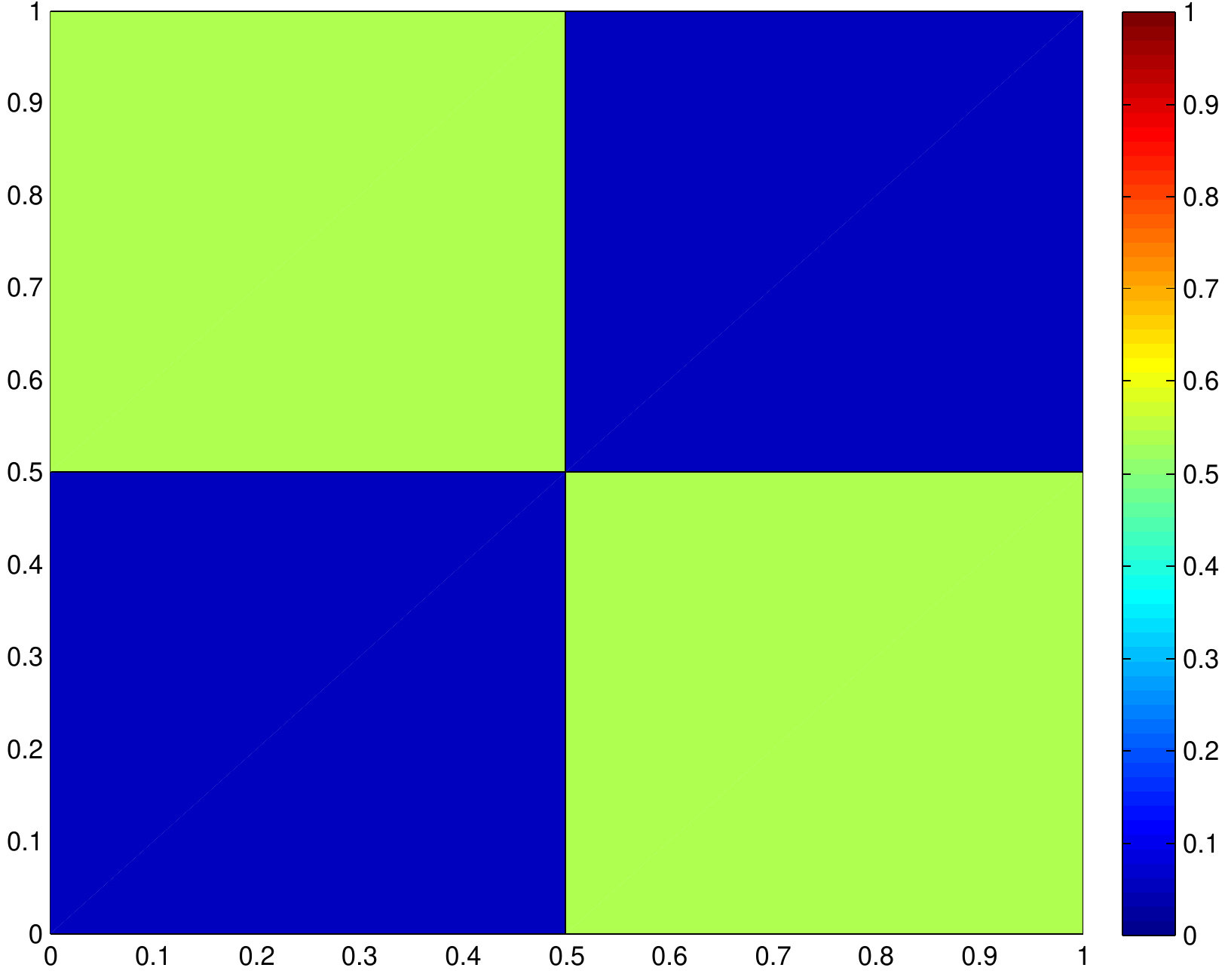} 
\includegraphics[angle=0,width=0.23\textwidth]{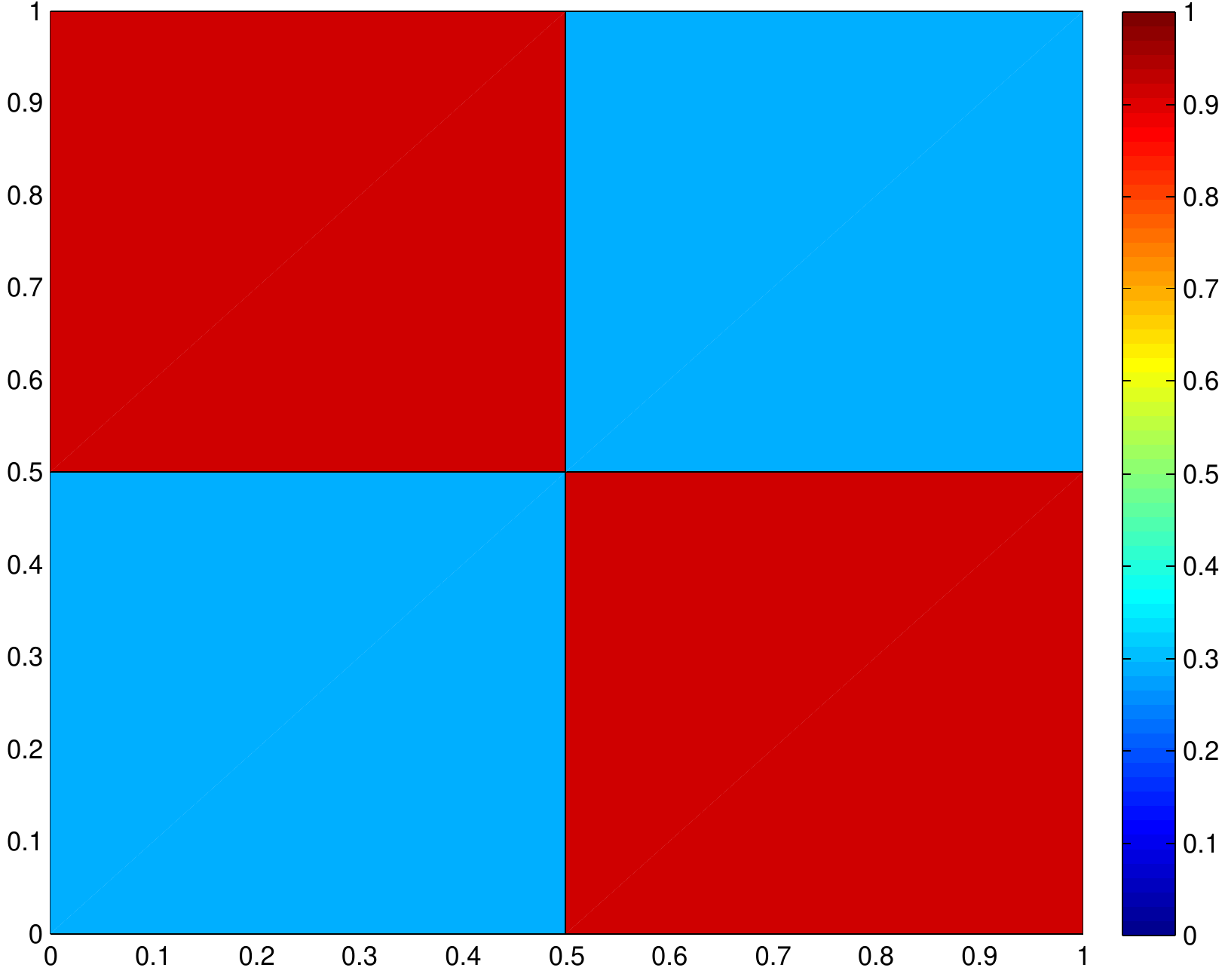} 
\includegraphics[angle=0,width=0.23\textwidth]{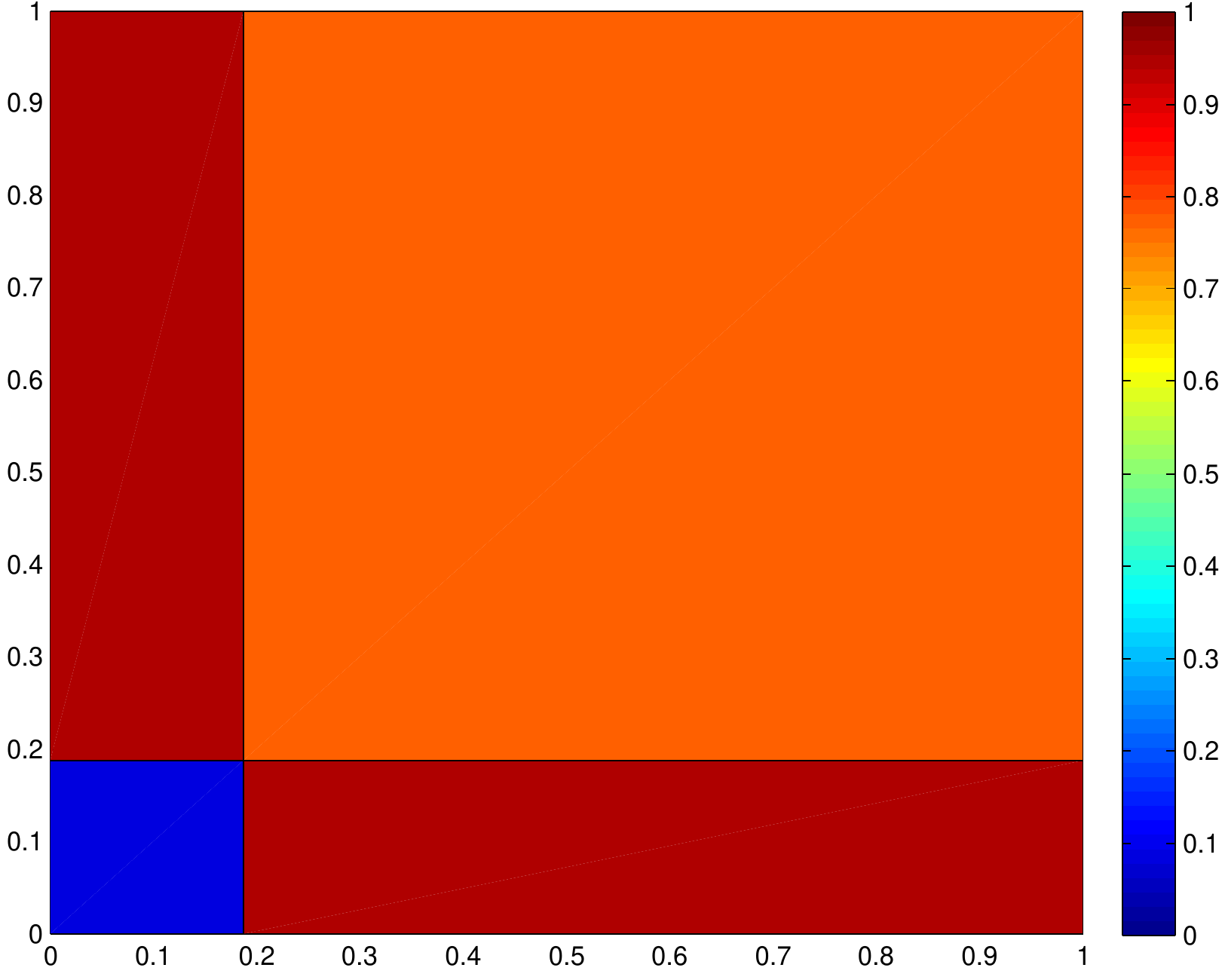}
\caption{Minimizing graphons at $\E=0.3$ (left column),
  $\E=0.6$ (middle column) and $\E=0.8$ (right column). For each
  column $\T$ values decrease from top to bottom. The second row
  represents points just above the ER curve, the third row points on
  the ER curve, and the fourth row points just below the ER curve. The exact $\T$ values are
  given in Section~\ref{SEC:Num}.}
\label{FIG:Graphons}
\end{figure}

The main feature of the computational results is that the minimal values of the rate function $I$ in the region we show in  Fig.~\ref{FIG:Phase-I-Value} are achieved by bi\partite graphons. In other words, the minimizing graphons in the whole subregion are bi\partite\!\!. On the ER curve, the graphons are constant. Thus the partition does not matter anymore. We consider them as bi\partite purely from the continuity perspective. We show in Fig.~\ref{FIG:Graphons} minimizing graphons at some typical points in the phase space. 
Shown are minimizing graphons at $\E=0.3$ (left column), $\E=0.6$ (middle column) and $\E=0.8$ (right column) respectively. For the $\E=0.3$ (resp. $\E=0.6$ and $\E=0.8$) column, the $\T$ values for the corresponding graphons from top to bottom are  $\T=0.09565838$ (resp. $\T=0.34037900$ and $\T=0.61784171$) which is in the middle of phase I,  $\T=0.03070755$ (resp. $\T=0.22010451$ and $\T=0.55677919$) which is just above the ER curve, $\T=0.02700000$ (resp. $\T=0.21600000$ and $\T=0.51200000$) which is on the ER curve, $\T=0.02646000$ (resp. $\T=0.21472000$ and $\T=0.51104000$) which is just below the ER curve, and $\T=0.01296000$ (resp. $\T=0.18400000$ and $\T=0.50800000$) which is in the middle of phase II (resp. phase II and phase III) respectively.
\begin{figure}[!ht]
\centering
\rotatebox{90}{\hspace*{4.1cm}{\Large $\tau$}}
\includegraphics[angle=0,width=0.6\textwidth]{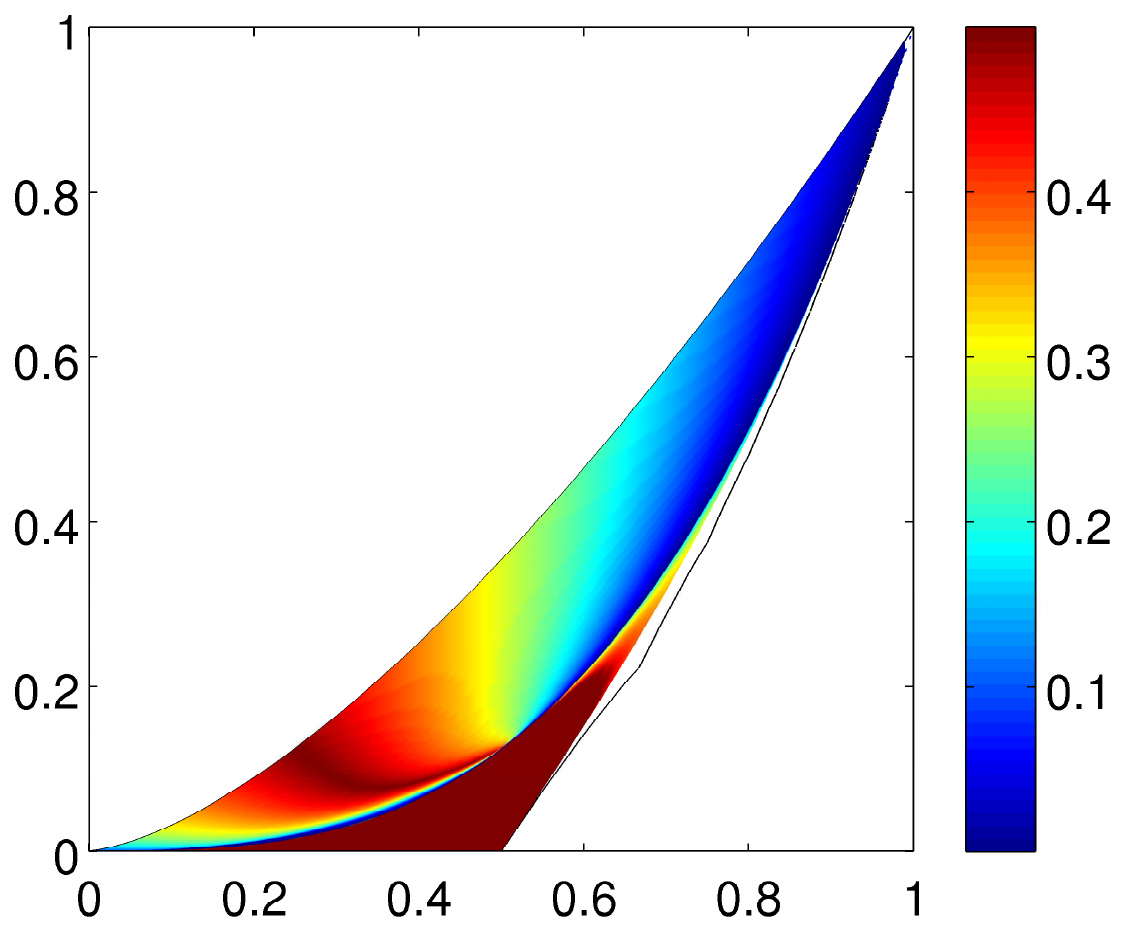}\\ 
\vspace*{-0.2cm}\hspace*{-0.7cm} {\Large $\epsilon$}
\caption{The minimal $c$ values of the minimizing bi\partite graphons as a function of $(\E,\T)$.}
\label{FIG:Phase-C-Value}
\end{figure}

\begin{figure}[!ht]
\centering
\includegraphics[angle=0,width=0.19\textwidth]{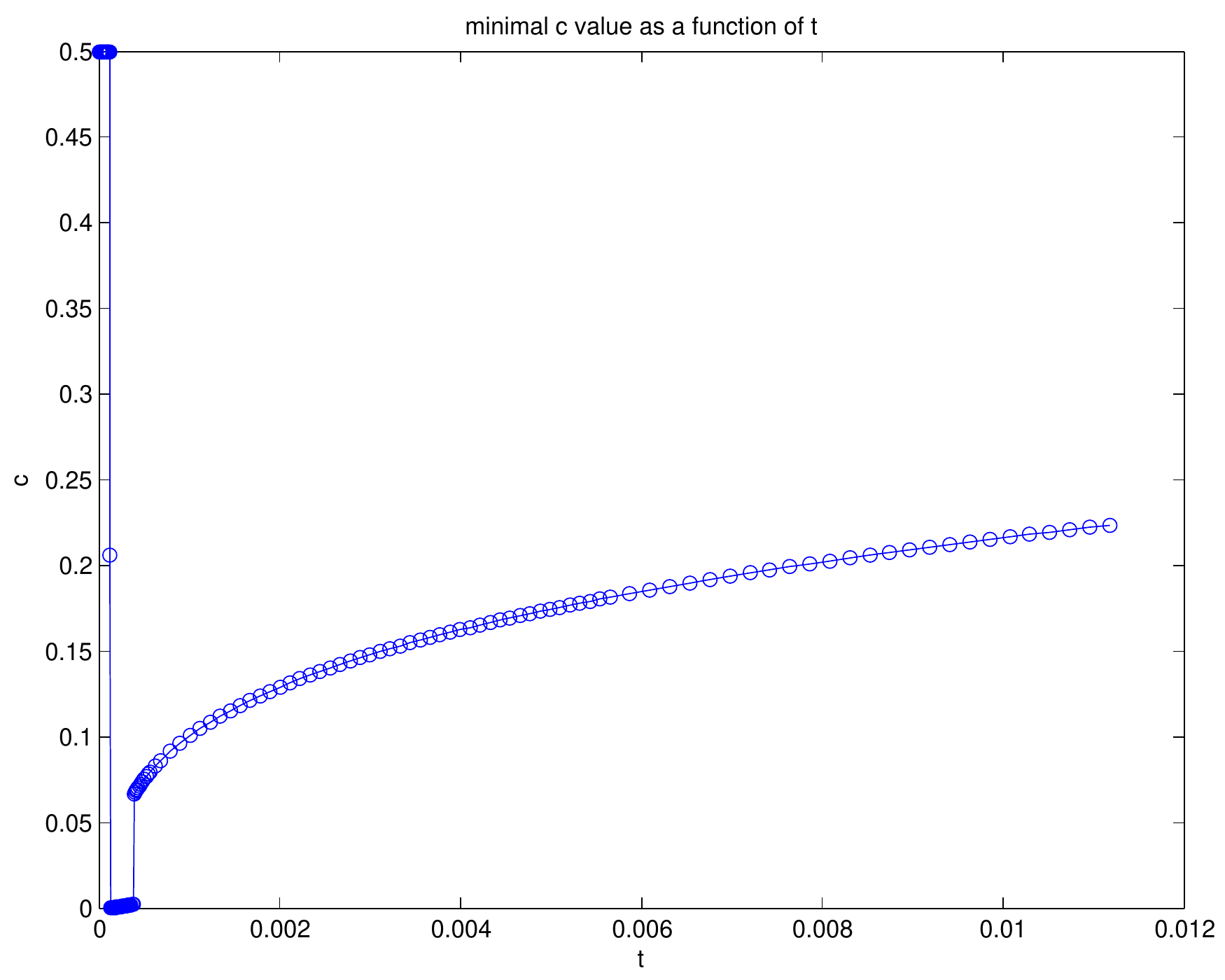} 
\includegraphics[angle=0,width=0.19\textwidth]{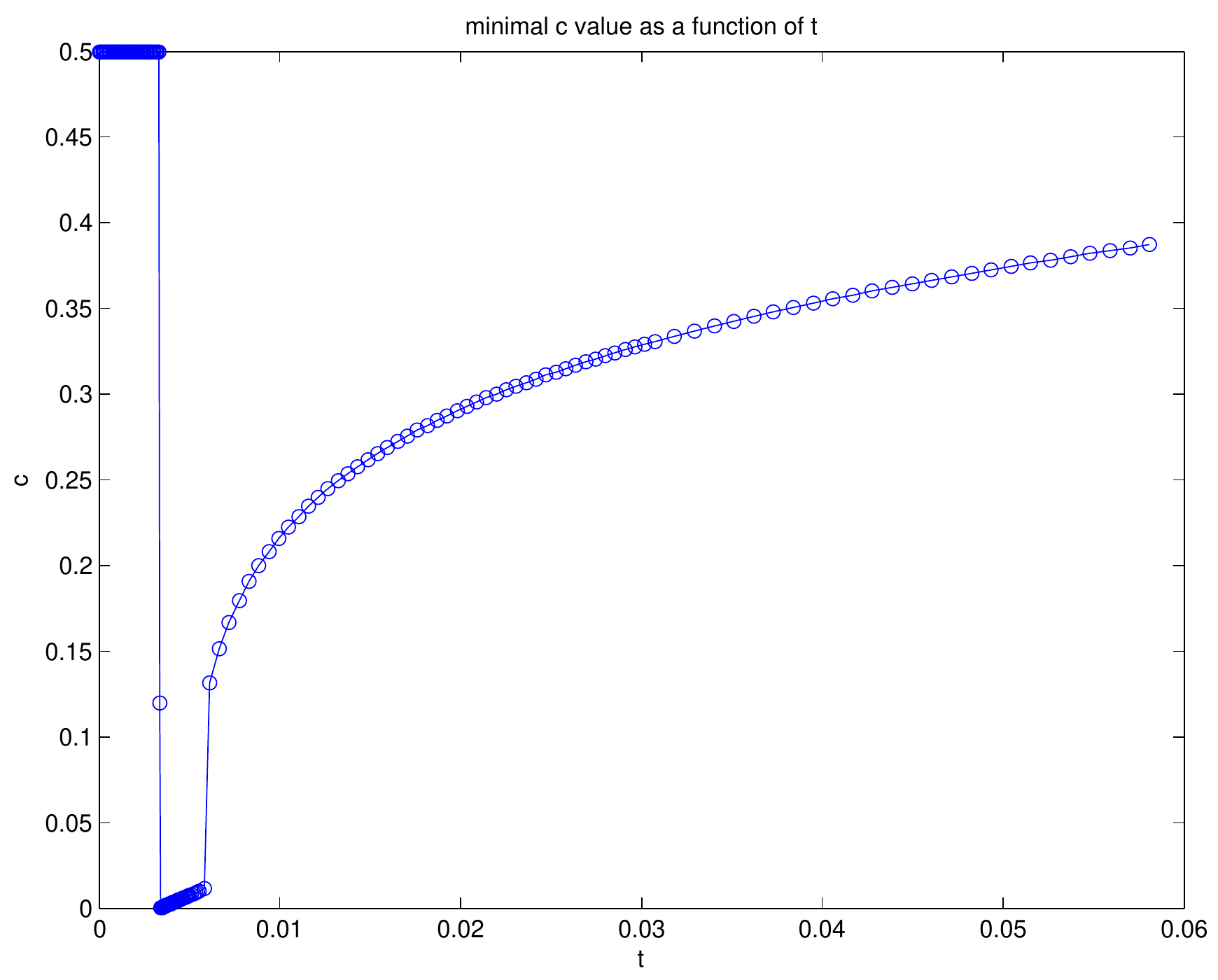} 
\includegraphics[angle=0,width=0.19\textwidth]{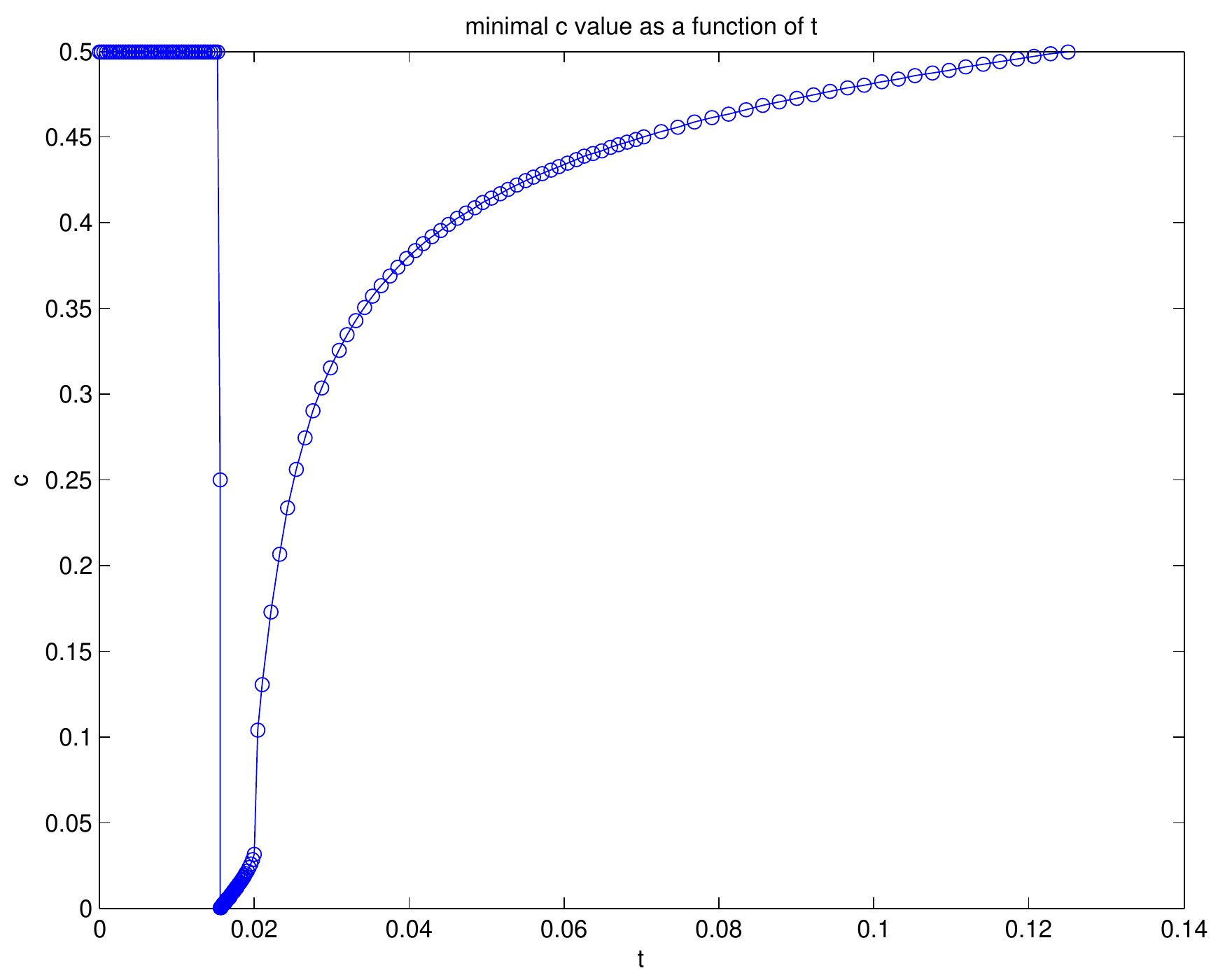} 
\includegraphics[angle=0,width=0.19\textwidth]{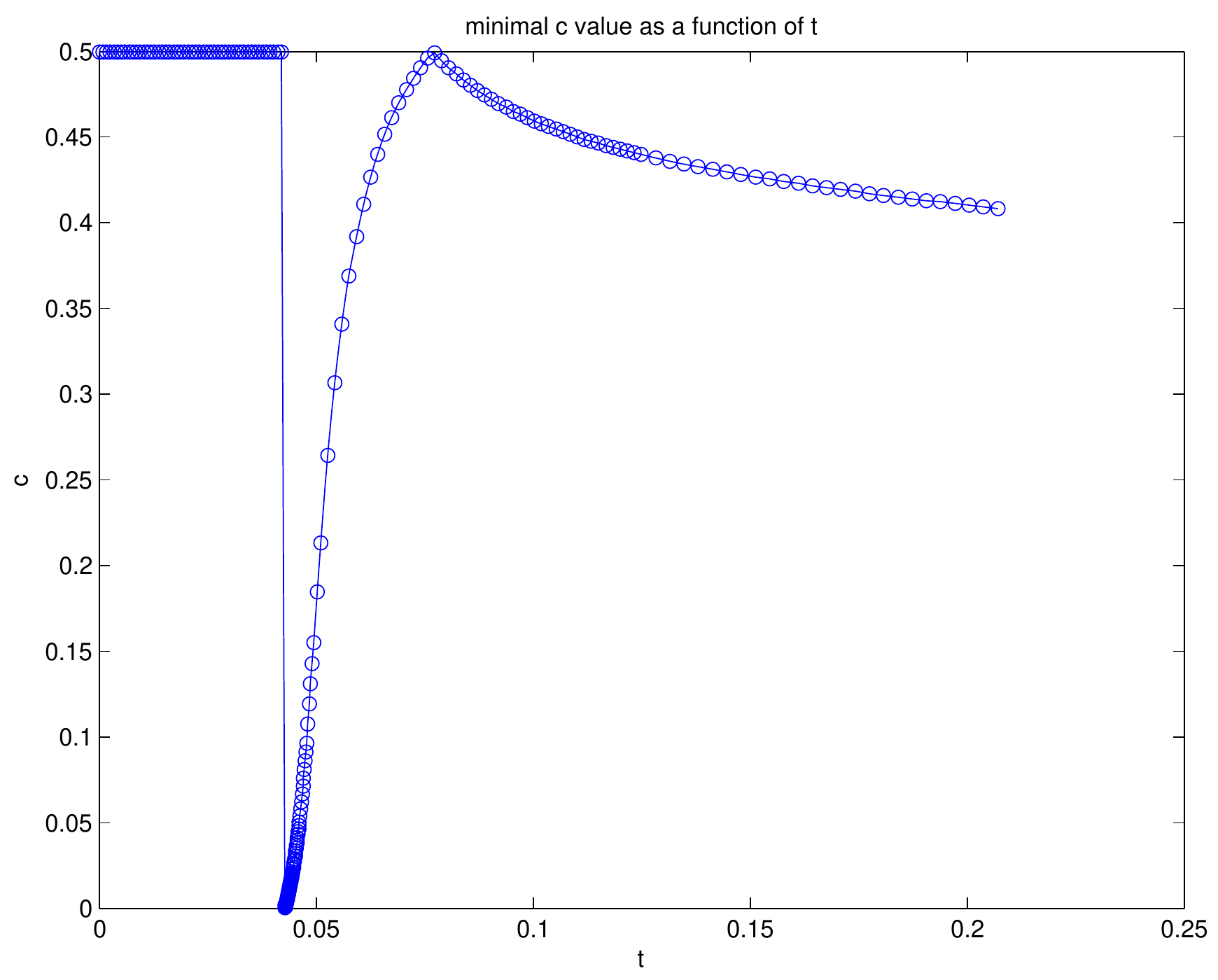} 
\includegraphics[angle=0,width=0.19\textwidth]{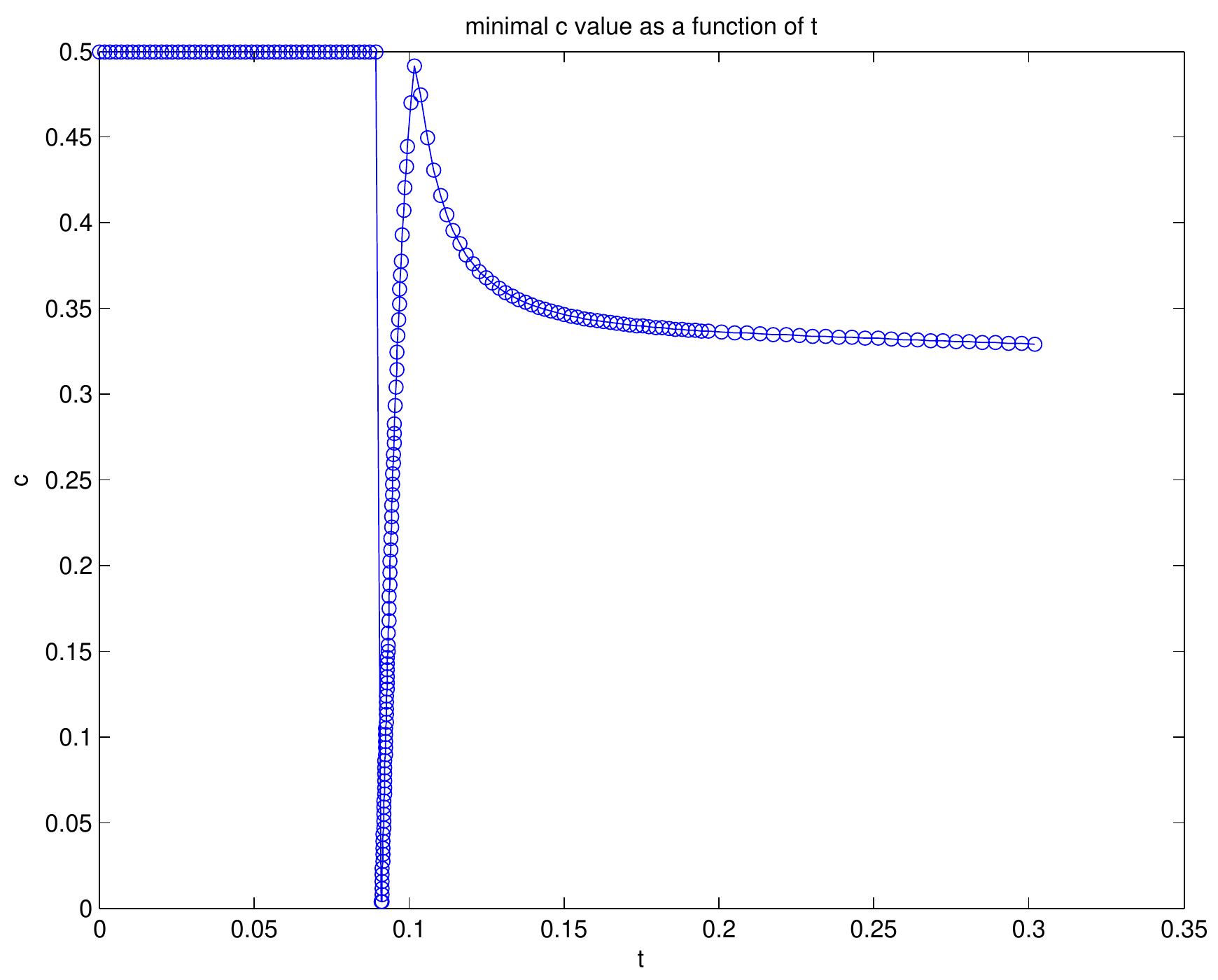}\\ 
\includegraphics[angle=0,width=0.19\textwidth]{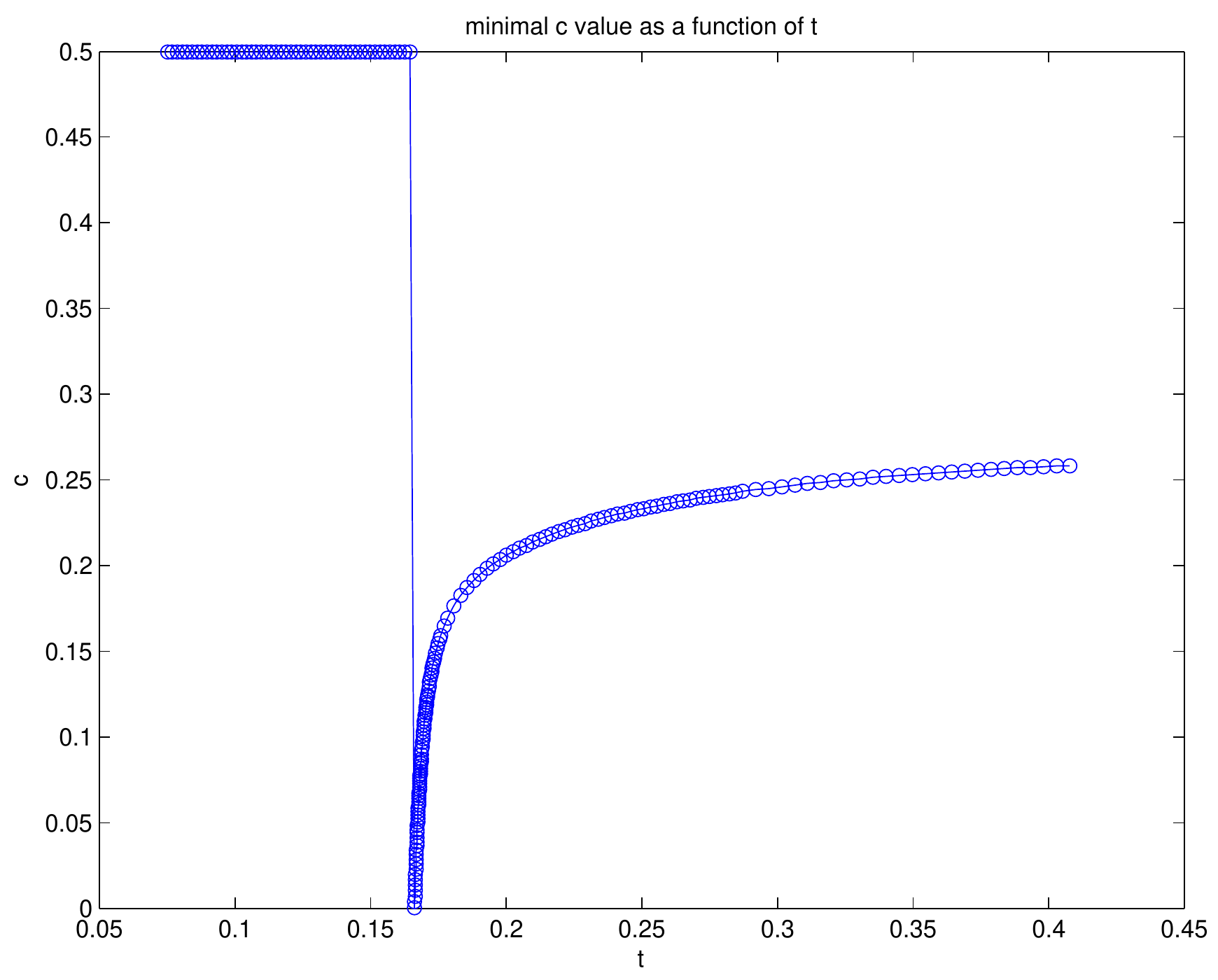} 
\includegraphics[angle=0,width=0.19\textwidth]{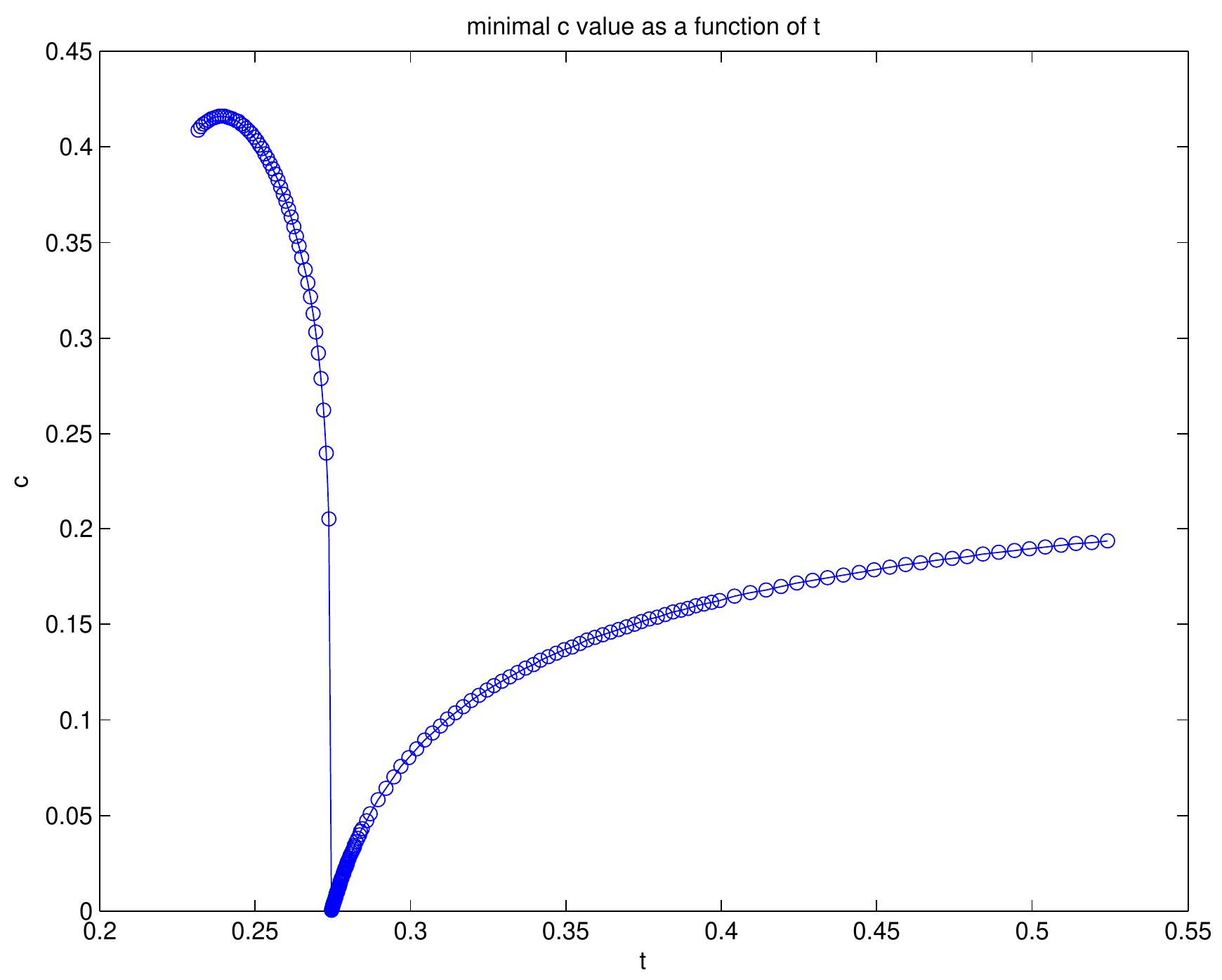} 
\includegraphics[angle=0,width=0.19\textwidth]{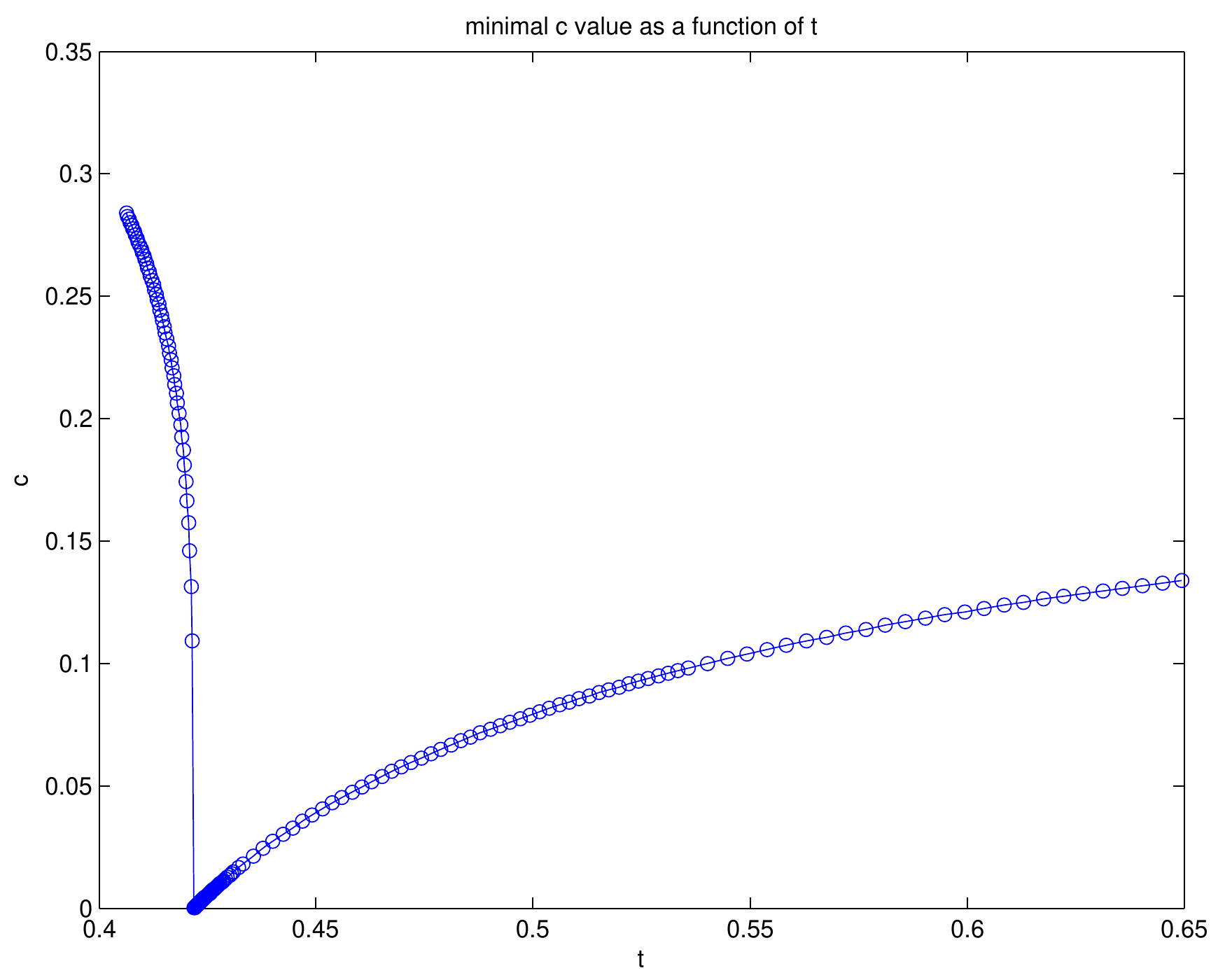} 
\includegraphics[angle=0,width=0.19\textwidth]{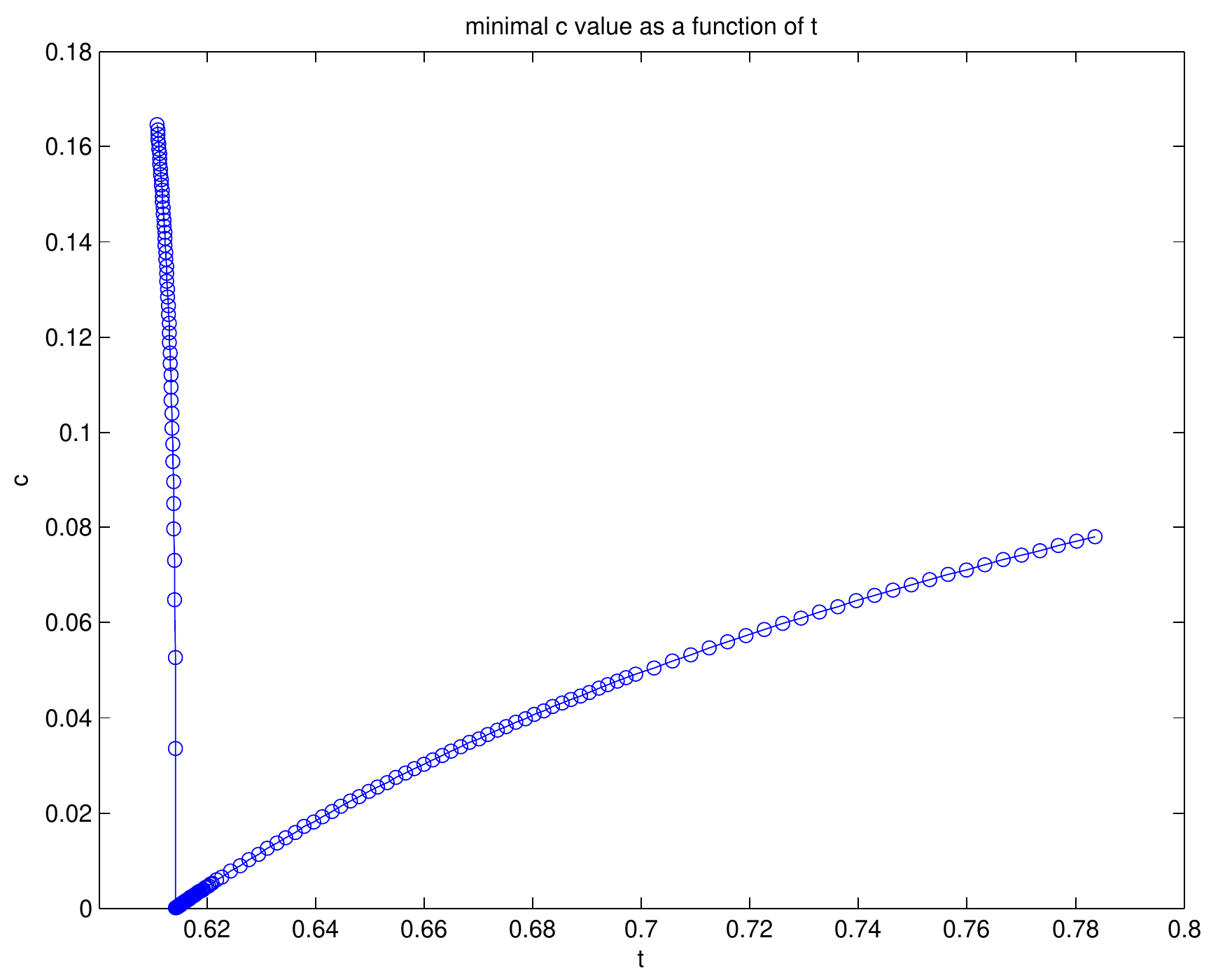} 
\includegraphics[angle=0,width=0.19\textwidth]{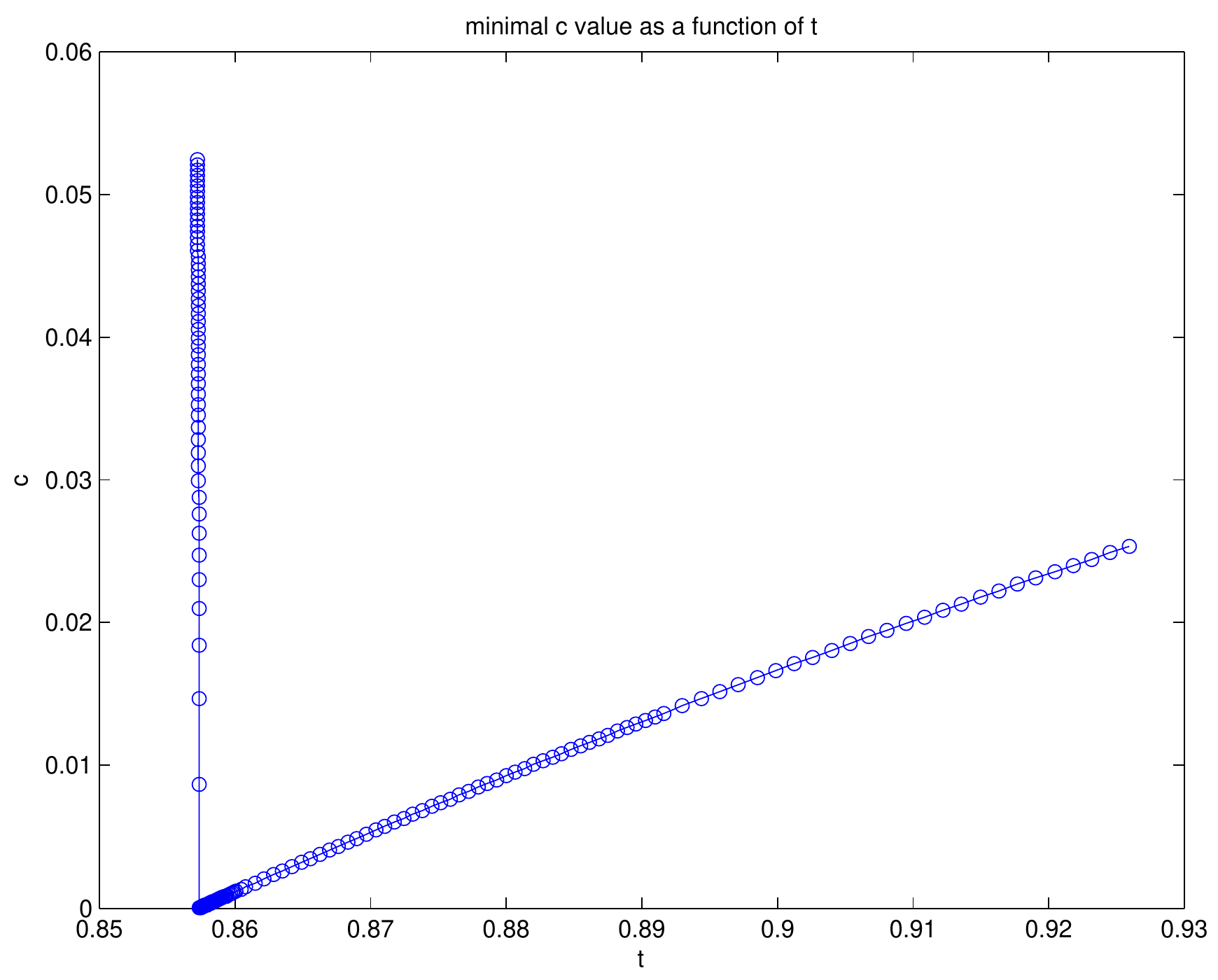} 
\caption{Cross-sections of the function $c(\E,\T)$ along lines of
  $\E=0.05$, $\E=0.15$, $\cdots$, $\E=0.95$.  The 4th and 5th graphs
  display coordinate singularities as the size of a cluster passes
  through 0.5, since we have chosen $c$ to always be 0.5 or less. The
  right-most kinks in those graphs are not phase transitions.}
\label{FIG:Phase-C-CS}
\end{figure}

A glance at the graphons in Fig.~\ref{FIG:Graphons} gives the impression that the minimizing graphons are very different at different values of $(\E,\T)$. A prominent feature is that the sizes of
the two blocks vary dramatically. Since there is only one parameter
that controls the sizes of the blocks, that is, if the first block has
size $c_1=c$ then the second block has size $c_2=1-c$, we can easily
visualize the change of the sizes of the two blocks in the minimizing
graphons in the phase space through $c$, the size of the smaller block as a function of $(\E,\T)$. We show in Fig.~\ref{FIG:Phase-C-Value} the minimal c values associated with the minimizing bi\partite graphons. The cross-sections along the lines $\E=a_k$ ($k=1,\cdots, 10$) are shown in Fig.~\ref{FIG:Phase-C-CS}.

\begin{figure}[ht]
\centering
\rotatebox{90}{\hspace*{4.1cm}{\Large $\tau$}}
\includegraphics[angle=0,width=0.6\textwidth]{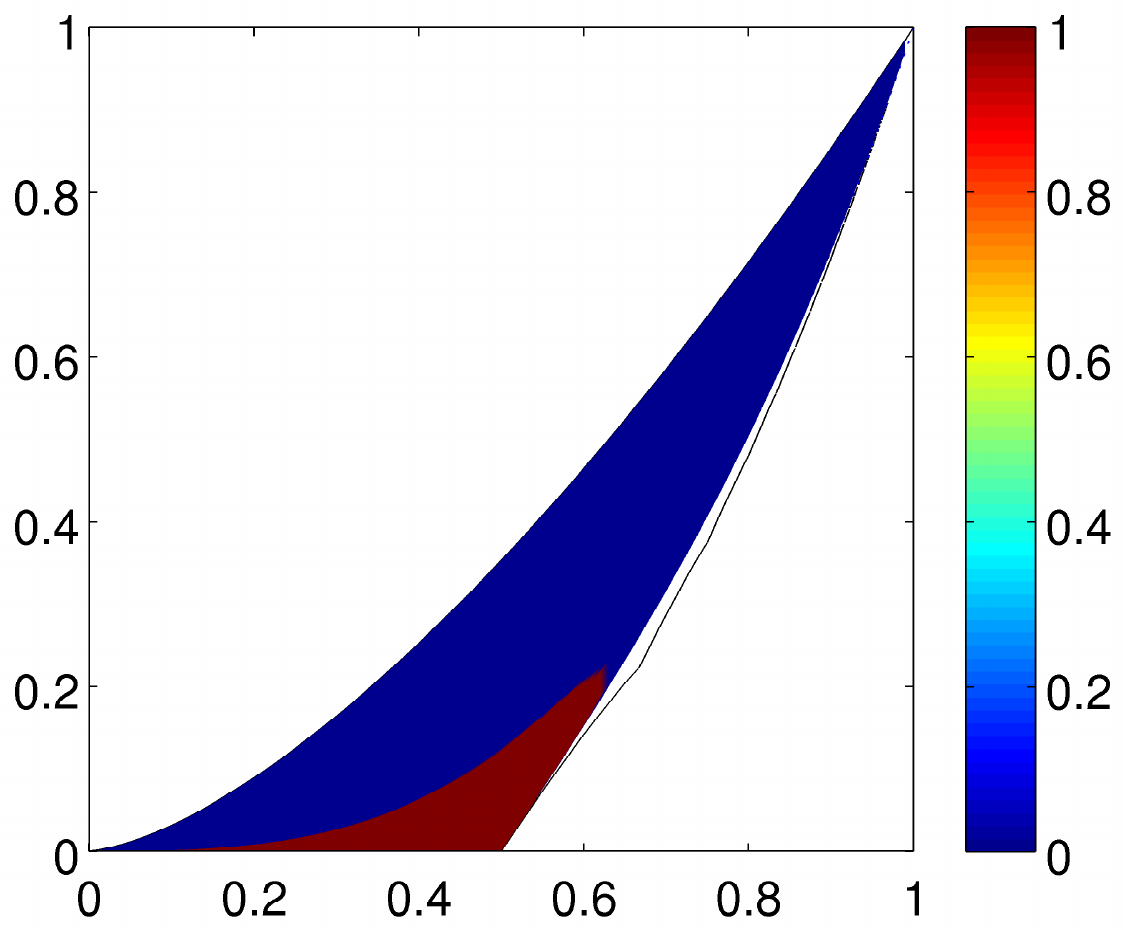}\\ 
\vspace*{-0.1cm}\hspace*{-0.60cm} {\Large $\epsilon$} 
\caption{The regions of symmetric bi\partite and asymmetric bi\partite optima.}
\label{FIG:Phase-Symmetric-Asymmetric}
\end{figure}

Among the minimizing bi\partite graphons some are symmetric and the
others are asymmetric. We show in
Fig.~\ref{FIG:Phase-Symmetric-Asymmetric} the region where the
minimizing bi\partite graphons are symmetric (brown) and the region
where they are asymmetric (blue). This is done by checking the
conditions $c_1=c_2=0.5$ and $g_{11}=g_{22}$ (both with accuracy up to
$10^{-7}$) for symmetric bi\partite graphons. Our numerical
computation here agrees with the theoretical results we obtain later
in Section~\ref{SUBSEC:Phases}.

\section{Local Analysis}
\label{SEC:Analysis}

In this section we do a perturbative analysis of all three phases near
the ER curve. This gives a qualitative explanation for why these three
phases appear in the form they do, and we compute {exactly} the
boundary between phases II and III, assuming the multipodal structure
of the phases.

\subsection{General Considerations}

We treat the graphon $g(x,y)$ as the integral kernel of an operator on
$L^2([0,1])$. Let $|\phi_1\rangle \in L^2([0,1])$ be the constant
function $\phi_1(x)=1$. Then the edge density is $e(g)=\E = \langle \phi_1 |
g | \phi_1 \rangle$ and the triangle density is $t(g)=\T = Tr(g^3)$.  On the
ER curve the optimal graphon is $g_0 = \E |\phi_1 \rangle \langle \phi_1
|$. Near the ER curve we take $g = g_0 + \delta g$, where $\langle \phi_1 |
\delta g | \phi_1 \rangle = 0$. For fixed $\E$, the rate function is
minimized at the ER curve, where $\T=\E^3$. For nearby values of $t(g)=\T =
\E^3 + \delta \T$, minimizing the rate function involves solving two
simpler optimization problems:
\begin{enumerate}
\begin{item}We want to minimize the size of $\delta g$, as measured in the $L^2$ norm, for a given $\delta \T$. \end{item}
\begin{item}We want to choose the form of $\delta g$ that minimizes $\delta I$ for a given $\|\delta g\|^2_{L^2}=
\iint \delta g(x,y)^2\, dxdy$. \end{item}
\end{enumerate}
If we could solve both optimization problems simultaneously, we would
have a rigorously derived optimum graphon. This is indeed what happens
when $\E=1/2$ and $\T<\E^3$ \cite{RS2}. For other values of $(\E,\T)$, 
the solutions to the two optimizations disagree somewhat. 
In phases I and III, the actual best graphon appears to be a compromise
between the two optima, while in phase II it appears to be a solution to the
first optimization problem, but of a form suggested by the second
problem. 

By playing the two problems off of one another, we derive candidates 
for the optimal graphon and gain insight into the numerical results. 
By measuring the extent to which these candidates
fail to solve each problem, we can derive rigorous estimates on the
behavior of the entropy. However, we do not claim to prove that our
bipodal ansatz is correct. 

A simple expansion of $g^3=(g_0 + \delta g)^3$ shows that 
\begin{equation} \delta \T = 3\E \langle \phi_1 | \delta g^2 | \phi_1 \rangle 
+ Tr(\delta g^3). \end{equation}
The first term is positive-definite, while the second is indefinite. 

For $\delta \T$ negative, the solution to the first optimization
problem is to have $\delta g = \nu |\psi \rangle \langle \psi |$,
where $\psi$ is an arbitrary normalized vector in $L^2([0,1])$ such
that $\langle \phi_1 | \psi \rangle = 0$. This eliminates the
positive-definite term and makes the indefinite term as negative as 
possible, namely equal to $\nu^3$.

For $\delta \T$ positive, the solution to the first optimization
problem is of the form
\begin{equation}\label{ansatz1}
\delta g = \frac{\mu}{\sqrt{2}} | \psi \rangle \langle \phi_1| + \frac{\mu}{\sqrt{2}} | \phi_1 \rangle \langle \psi|
+ \nu |\psi \rangle \langle \psi |,
\end{equation}
with $\langle \phi_1 | \psi \rangle = 0$. 
We then have 
\begin{equation}
\|\delta g\|_{L^2}^2 = \mu^2 + \nu^2; \qquad \delta \T = \nu^3 + \frac{3 \mu^2}{2}(\E+\nu).
\end{equation}
Maximizing $\delta \T$ for fixed $\mu^2+\nu^2$ is then a Lagrange
multiplier problem in two variables with one constraint. There are
minima for $\delta \T$ at $\mu=0$ and maxima at $\nu =
\mu^2/(2\E)$.

Note that the form of the vector $\psi$ did not enter into this
calculation.  The form of $\psi$ is determined {\bf entirely} from
considerations of $\delta I$ {\it versus} $\|\delta g\|_{L^2}$. We now 
turn to this part of the problem. 

The rate function $I_0(u)$ is
an even function of $u-{1}/{2}$. Furthermore, all even derivatives of
this function at $u={1}/{2}$ are positive, so the function can be written as a 
power series $I_0(u) = -{(\ln 2)}/{2} + \sum_{n=1}^\infty c_n 
\left (u-{1}/{2} \right )^{2n}$, where the coefficients $c_n$ are all
positive. As a function of $(u -{1}/{2})^2$, $I_0(u)$ is then concave up. 
Thus the way to minimize $I(g) = \iint I[g(x,y)]\, dxdy$ for fixed 
$\iint \left (g(x,y) - \frac{1}{2}\right)^2\, dxdy$ is to have 
$\left [g(x,y) - {1}/{2}\right]^2$
constant; in other words to have $g(x,y)$ only take on two values, and for
those values to sum to 1.

This is exactly our second minimization problem. We have a fixed value of 
$\|\delta g\|_{L^2}^2$, and hence a fixed value of 
\begin{eqnarray}
\iint \left ( g(x,y) - \frac{1}{2} \right)^2\, dxdy & =& 
\iint \left (\E - \frac{1}{2} + \delta g(x,y)\right )^2\, dxdy \cr
&=&  
\left ( \E - \frac{1}{2}\right)^2 + 2\left (\E -\frac{1}{2}\right) 
\iint \delta g(x,y)\, dxdy + \iint \delta g(x,y)^2\, dxdy \cr  
&=& \left(\E-\frac{1}{2}\right)^2 + \|\delta g\|_{L^2}^2. 
\end{eqnarray}
In order for $g(x,y)$ to only take on two values, we need $\psi$ to
only take on two values. After doing a measure-preserving transformation
of $[0,1]$, we can assume that 
\begin{equation} \label{psic}
\psi(x) = \psi_c(x) := \begin{cases} \sqrt{\frac{1-c}{c}} & 
x < c \cr & \cr -\sqrt{\frac{\phantom{1}c\phantom{1}}{1-c}} & x > c \end{cases}
\end{equation}
for some constant $c\le \frac{1}{2}$. 

We henceforth restrict our attention to the ansatz (\ref{ansatz1}) 
with the function (\ref{psic}), and vary the parameters $\mu$, $\nu$
and $c$ to minimize $I(g)$ while preserving the fixed values of $\E$ and $\T$. 

\subsection{Phase I: $\T>\E^3$}

Suppose that $\T$ is slightly greater than $\E^3$.  By our previous analysis,
we want $\delta g$ to be small in an $L^2$ sense. Maximizing $\delta \T$
for fixed $\|\delta g\|_{L^2}$ means taking $\nu = {\mu^2}/{(2\E)} \ll \mu$,
so
\begin{equation}
\delta \T = \frac{3\mu^2\E}{2} + \frac{3\mu^4}{4} + \frac{\mu^6}{8\E^3};
\qquad \|\delta g\|_{L^2}^2 = \mu^2 + \frac{\mu^4}{4\E^2}. 
\end{equation}

If $\E \ne 1/2$, then there is no way to make $\left |g(x,y) - {1}/{2}
\right|$ constant while keeping $\delta g$ pointwise small. Instead,
we take $\delta g(x,y)$ to be large in a small region. 
We compute 
\begin{equation}
g(x,y) = \begin{cases} \E - \mu \sqrt{2} \sqrt{\frac{\phantom{1}c\phantom{1}}{1-c}} + \nu \frac{1}{1-c}
& x,y > c \cr
\E + \frac{\mu}{\sqrt{2}} \left ( \sqrt{\frac{1-c}{c}} - \sqrt{\frac{\phantom{1}c\phantom{1}}{1-c}}
\right ) - \nu  & x < c < y \hbox{ or }y<c<x \cr 
\E + \mu \sqrt{2} \sqrt{\frac{1-c}{c}} + \nu \frac{1-c}{c} & x,y<c
\end{cases}
\end{equation}

By taking $c = {\mu^2}/{2 (2\E-1)^2} + O(\mu^4)$, we can get the
values of $g$ in the rectangles $x<c<y$ (or $y<c<x$) and the large
square $x,y>c$ to sum to 1.  This makes $I_0[g(x,y)]$ constant except
on a small square of area $c^2 = O(\mu^4)$. Since $g(x,y)$ is equal to
$\E - {\mu^2}/{(1-2\E)} + O(\mu^4)$ when $x, y >c$, we have
\begin{eqnarray}
I(g) &=& I_0\left (\E - \frac{\mu^2}{1-2\E} \right)  + O(\mu^4) 
\cr & =&  I_0(\E) - I_0'(\E) \frac{\mu^2}{1-2\E} + O(\mu^4) \cr 
& = & I_0(\E) - \frac{2 I_0'(\E) \delta \T}{3\E(1-2\E)} + O(\delta t^2)
\end{eqnarray} 
Thus 
\begin{equation}
\frac{\partial s(\E,\T)}{\partial \T} = \frac{\ln \left ( \frac{\E}{1-\E}\right )}{3\E(1-2\E)} + O(\delta \T)
\end{equation}
just above the ER curve.

This ansatz gives an accurate description of our minimizing graphons
near the ER curve, and these extend continuously all the way to the
upper boundary $\T = \E^{3/2}$. Although the behavior in a first-order
neighborhood of the ER curve changes discontinuously when $\E$ passes
through 1/2, the behavior a finite distance above the ER curve appears
to be analytic in $\E$ and $\T$, so that the entire region between the
ER curve and the upper boundary corresponds to a single phase. As can
be seen from these formulas, or from the numerics of
Section~\ref{SEC:numer}, phase I has the following qualitative
features:

\begin{enumerate}
\begin{item}The minimizing graphon is bi\partite\!\!. The smallest value of
$g(x,y)$ is found on one of the two square regions $x,y<c$ or $x,y>c$, 
and the largest value is found on the other, with the value on the 
rectangles $x<c<y$ and $y<c<x$ being intermediate.
\end{item}
\begin{item}If $\E \ne 1/2$, then as $\T$ approaches $\E^3$ from above, 
$c$ goes to 0. The value of
$g(x,y)$ is then approximately $\E$ on the large square $x,y>c$ and $1-\E$ on
the rectangles. If $1/3 < \E < 2/3$, then the value of $g(x,y)$ 
on the small square 
$x,y < c$ is approximately $2-3\E$; the value is close to 1 if $\E < 1/3$
and close to zero if $\E>2/3$. 
\end{item}
\begin{item}Conversely, as $\T$ moves away from $\E^3$, $c$ grows quickly, 
to first order in $\delta \T$. When $\E<1/2$, the small-$c$ region is quite
narrow, as one sees in the second row of Fig.~\ref{FIG:Graphons}, and is easy to miss with numerical 
experiments. Although rapid, the transition from small-$c$ to large-$c$ 
appears to be continuous and smooth, with no sign of any phase transitions.
\end{item}
\begin{item}
As $\T$ approaches $\E^{3/2}$ (from below, of course), the value of $g(x,y)$
approaches 1 on a square of size $\sqrt{\E}$, approaches zero relatively 
slowly on the rectangles, and approaches zero more quickly on the other square. 
In a graph represented by such a graphon, there
is a cluster representing a fraction $\sqrt{\E}$ of the vertices. The probability
of two vertices within the cluster being connected is close to 1, the 
probability of a vertex inside the cluster being connected to a vertex outside
the cluster is small, and the probability of two vertices outside the 
cluster being connected is very small. 
\end{item}
\end{enumerate}
 
\subsection{Phases II and III: $\T<\E^3$}
\label{SUBSEC:Phases}

\begin{figure}[ht]
\begin{picture}(400,150) 
\setlength{\unitlength}{2pt}
\put(10,10){\line(0,3){60}}
\put(10,10){\line(3,0){60}}
\put(10,70){\line(3,0){60}}
\put(70,10){\line(0,3){60}}
\put(10,35){\line(3,0){60}}
\put(35,10){\line(0,3){60}}
\put(5,10){$0$}
\put(10,5){$0$}
\put(5,33){$c$}
\put(33,5){$c$}
\put(68,5){$1$}
\put(5,67){$1$}
\put(12,20){$\E + \nu\frac{1-c}{c}$}
\put(15,50){$\E-\nu$}
\put(45,20){$\E-\nu$}
\put(40,50){$\E + \nu \frac{c}{1-c}$}

\put(90,10){\line(0,3){60}}
\put(90,10){\line(3,0){60}}
\put(90,70){\line(3,0){60}}
\put(150,10){\line(0,3){60}}
\put(90,40){\line(3,0){60}}
\put(120,10){\line(0,3){60}}
\put(85,10){$0$}
\put(90,5){$0$}
\put(84,38){$\frac12$}
\put(119,3){$\frac12$}
\put(148,5){$1$}
\put(85,67){$1$}
\put(91,22){$\E \!-\!  (\E^3\!-\! \T)^{\frac{1}{3}}$}
\put(121,52){$\E \!-\!  (\E^3\!-\! \T)^{\frac{1}{3}}$}
\put(91,52){$\E \!+\!  (\E^3\!-\! \T)^{\frac{1}{3}}$}
\put(121,22){$\E \!+\!  (\E^3\!-\! \T)^{\frac{1}{3}}$}

\put(170,10){\line(0,3){60}}
\put(170,10){\line(3,0){60}}
\put(170,70){\line(3,0){60}}
\put(230,10){\line(0,3){60}}
\put(170,35){\line(3,0){60}}
\put(195,10){\line(0,3){60}}
\put(165,10){$0$}
\put(170,5){$0$}
\put(165,33){$c$}
\put(193,5){$c$}
\put(228,5){$1$}
\put(165,67){$1$}
\put(180,20){$0$}
\put(180,50){$1$}
\put(210,20){$1$}
\put(210,50){$\alpha$}

\end{picture}
\caption{The graphons for the expressions in ~\eqref{EQ:Graphon1} (left), ~\eqref{EQ:Graphon2} (middle) and ~\eqref{EQ:Graphon3} (right) respectively.}
\label{FIG:Graphon Anal}
\end{figure}
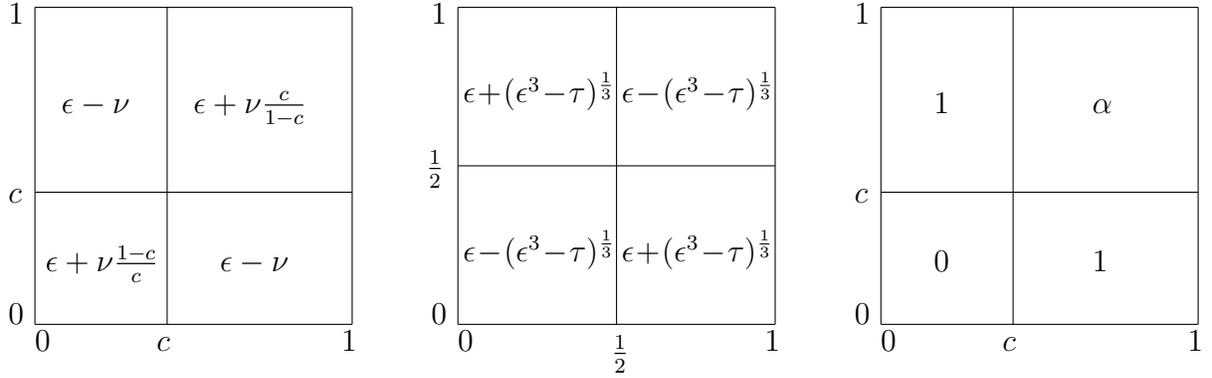

When $\T < \E^3$, we want to minimize the $\mu^2$ term in $\delta \T$ and
make the $\nu^3$ term as negative as possible. This is done by simply taking
$\mu=0$ and $\nu = (\delta \T)^{1/3}$.
This gives a graphon of the form

\begin{equation}\label{EQ:Graphon1}
 g(x,y) = \begin{cases} \E + \nu \frac{1-c}{c} & x,y < c \cr
\E -\nu & x<c<y \hbox{ or } y<c<x \cr 
\E + \nu \frac{\phantom{1}c\phantom{1}}{1-c} & x,y>c\end{cases}
\end{equation}
which is displayed in the left of Fig.~\ref{FIG:Graphon Anal}.

Note that $\nu$ is negative, so the value of the graphon is less than
$\E$ on the squares and greater than $\E$ on the rectangles. The two
vertex clusters
of fraction $c$ and $1-c$ prefer to connect with each other than with 
themselves. 

When $c \ne 1/2$, this gives a graphon that takes on 3 different values.
When $c=1/2$, the graphon takes on only two values, but these values
do not sum to 1. Only when $\E=c=1/2$ can we simultaneously solve both
optimization problems. 

To understand which problem ``wins'' the competition, we look at the 
situation when
$\E<1/2$ and $\T$ is slightly less than $\E^3$. In this region, the 
entropy cost for having $\mu \ne 0$ is lower order in $\E^3-\T$
than the entropy cost for violating the second problem. This shows that
the optimal bipodal graphon must have $c$ close to $1/2$, validating
the use of perturbation theory around $c=1/2$. 

In subsection \ref{transition}, we compute the
second variation of the rate function with respect to $c$ at $c=1/2$,
and see that it is positive for all values of $\T<\E^3$ when $\E<1/2$, 
and for some values of $\T$ when $1/2 < \E < 0.629497839$. In this
regime, the best we can do is to pick $c=1/2$
and $\mu$ exactly zero. Our graphon then takes the form
\begin{equation}\label{EQ:Graphon2}
g(x,y) = \begin{cases} \E - (\E^3-\T)^{1/3} & x,y < 1/2 \hbox{ or } x,y > 1/2 \cr
\E + (\E^3-\T)^{1/3} & x< \frac{1}{2} < y \hbox{ or } y < \frac{1}{2} < x
\end{cases}
\end{equation}
These symmetric bi\partite graphons, displayed in the middle of
Fig.~\ref{FIG:Graphon Anal}, are the optimizers for phase II.


The symmetric bipodal phase II has a natural boundary.  
If $\E>1/2$, then the smallest
possible value of $\T$ is when $\nu = \E-1$, in which case $\T = 2\E^3
-3\E^2+3\E-1$. It is possible to get lower values of $\T$ with $c \ne
1/2$, and still lower values if the graphon is not
bi\partite\!\!. Thus there must be a phase transition between phase II
and the asymmetric bi\partite phase III.

Phase III has its own natural boundary. Among bi\partite graphons,
the minimum triangle density 
for a given edge density (with $\E>1/2$) is given by a graphon of the
form 
\begin{equation}\label{EQ:Graphon3}
 g(x,y) = \begin{cases} 0 & x,y<c \cr 
1 & x<c<y \hbox{ or } y<x<c \cr 
\alpha & x,y > c 
\end{cases}
\end{equation}
displayed in the right of Fig.~\ref{FIG:Graphon Anal}, where $c,\alpha \in [0,1]$ are free parameters. 


The edge and triangle densities are then
\begin{equation} 
\E = 2c(1-c) + (1-c)^2 \alpha; \qquad \T = (1-c)^3 \alpha^3 + 3c(1-c)^2\alpha.
\end{equation}
Minimizing $\T$ for fixed $\E$ yields a sixth order polynomial equation in $c$,
which we can solve numerically. 

The main qualitative features of phases II and III are
\begin{enumerate}
\begin{item}The minimizing graphon consists of two squares, of side
$c$ and $1-c$, and two $c \times (1-c)$ rectangular
regions. The value of the graphon is largest in the rectangular regions and
smaller in the squares. In phase II we have $c=1/2$, and the value of the
graphon is the same in both squares. In phase III we have $c < 1/2$, and
the value of the graphon is smallest in the $c \times c$ square.
\end{item}
\begin{item}
All points with $\E<1/2$ and $\T<\E^3$ lie in phase II. As $\T \to 0$, the graphon
describes graphs that are close to being bipartite\!\!, with two clusters of equal
size such that edges {\em within} a cluster appear with probability 
close to zero,
and edges {\em between} clusters appear with probability close to $2\E$. 
\end{item}
\begin{item}
When $\E>1/2$, phase II does {\bf not} extend to its natural boundary. 
There comes a point
where one can reduce the rate function by making $c$ less than 1/2. 
For instance, when $\E=0.6$ one can construct symmetric bi\partite graphons
with any value of $\T$ between $2\E^3 -3\E^2+3\E-1=0.152$ and $\E^3=0.216$. 
However, phase II only extends between $\T= 0.152704753$ and $\T=0.20651775$,
with the remaining intervals being phase III.
The actual boundary of phase II is computed in the next subsection. 
\end{item}
\begin{item}
Phase III, by constrast, appears to extend all the way to its natural
boundary, shown in Fig.~\ref{FIG:Phase-Boundary}, to the current numerical resolution in the $\T$ variable. The union of phases I, II, and III appears to be all values of $(\E,\T)$ that can be 
achieved with a bi\partite graphon. 
\end{item}
\begin{item}When we cross the natural boundary of phase III, the larger of 
the two clusters breaks into two equal-sized pieces, with slightly different
probabilities of having edges within and between those sub-clusters. 
The resulting tri\partite graphons have $g(x,y)$ strictly between 0 and 1,
and (for $\E<2/3$) extend continuously down to the $k=3$ scallop. 
\end{item}
\begin{item}Throughout phases I, II, and III, there appears to be 
a unique graphon that minimizes the rate function. We conjecture that this 
is true for all $(\E,\T)$, and that the minimizing graphon for each $(\E,\T)$ 
is $k$-\partite with $k$ as small as possible. 
\end{item}
\end{enumerate}

\subsection{The boundary between phase II and 
phase III}\label{transition}

In both phase II and phase III, we can express our graphons in the
form (\ref{ansatz1}). For each value of $c$ we can vary $\mu$ and $\nu$
to minimize the rate function, while preserving the constraint on $\T$. 
Let $\mu(c)$ and $\nu(c)$ be the values of $\mu$ and $\nu$ that achieve
this minimum for fixed $c$. Note that changing $c$ to $1-c$ and changing
$\mu$ to $-\mu$ results in the same graphon, up to reparametrization of 
$[0,1]$. By this symmetry, $\mu(c)$ must be an odd function of 
$c - {1}/{2}$ while $\nu(c)$ must be an even function. 

Since $\T = \E^3 + \nu^3 + 3(\E+\nu)\mu^2/2$ is independent of $c$, we 
obtain constraints on derivatives of $\mu(c)$ and $\nu(c)$ by computing
\begin{eqnarray}
0=\T'(c) & = & 3 \nu^2 \nu' + \frac{3\mu^2\nu'}{2}+ 3(\E+\nu)\mu\mu' \cr
0 = \T''(c) & = & 6 \nu (\nu')^2 + 3 \nu^2 \nu'' + \frac{3\mu^2 \nu''}{2}
+ 6 \mu \mu'\nu' + 3(\E+\nu)(\mu')^2 + 3(\E+\nu) \mu \mu'',
\end{eqnarray}
where $'$ denotes $d/dc$.  Since $\mu(1/2)=0$, the first equation implies
that $\nu'(1/2)=0$, while the second implies that 
\begin{equation}\label{nuequation}
\nu'' = - \frac{\E+\nu}{\nu^2}(\mu')^2
\end{equation}
at $c=1/2$. 

We next compute the rate function. We expand the formula (\ref{ansatz1}) as
\begin{equation} g(x,y) = \begin{cases}
\E + \nu \frac{1-c}{c} + \mu \sqrt{2} \sqrt{\frac{1-c}{c}} & x,y < c \cr 
\E - \nu + \frac{\mu}{\sqrt{2}} \left ( \sqrt{\frac{1-c}{c}} 
- \sqrt{\frac{\phantom{1}c\phantom{1}}{1-c}} \right ) & x < c < y \hbox{ or } 
y<c<x \cr 
\E + \nu \frac{c}{1-c} - \mu \sqrt{2} \sqrt{\frac{\phantom{1}c\phantom{1}}{1-c}} & x,y < c. 
\end{cases}
\end{equation}
This makes the rate function 
\begin{eqnarray}
I(g) &=& c^2I_0\left (\E+\nu\frac{1-c}{c} + \sqrt{2}\mu \sqrt{\frac{1-c}{c}}
\right ) \cr 
&&+ 2c(1-c)I_0\left (\E - \nu + \frac{\mu}{\sqrt{2}}
\left ( \sqrt{\frac{1-c}{c}} 
- \sqrt{\frac{\phantom{1}c\phantom{1}}{1-c}} \right ) \right )   \cr 
&& + (1-c)^2 I_0\left (\E + \nu \frac{c}{1-c} 
- \mu \sqrt{2} \sqrt{\frac{c}{1-c}}\right) 
\end{eqnarray}
Taking derivatives at $c=1/2$ and plugging in (\ref{nuequation}) yields 
\begin{equation}
\left . \frac{d^2 I}{dc^2}\right |_{c=1/2} = A (\mu')^2 + B \mu' + C,
\end{equation}
where
\begin{eqnarray}
A & = & I_0''(\E+\nu) + \frac{(\E+\nu)[I_0'(\E-\nu)-I_0'(\E+\nu)]}{2\nu^2} \cr
B & = & 2\sqrt{2} \left [I_0'(\E+\nu)-I_0'(\E-\nu)-2\nu I_0''(\E+\nu)\right] \cr
C & = & 4[I_0(\E+\nu) -2\nu I_0'(\E+\nu) + 2\nu^2 I_0''(\E+\nu) - I_0(\E-\nu)],
\end{eqnarray}
and $\nu = -(\E^3-\T)^{1/3}$.  The actual value of $\mu'$ is the one that
minimizes this quadratic expression, namely $\mu'=-B/(2A)$. At this value 
of $\mu'$, $d^2I/dc^2$ equals ${(4AC-B^2)}/{(4A)}$. 
As long as the discriminant
$B^2-4AC$ is negative, $d^2I/dc^2>0$ and the phase II graphon is stable 
against changes in $c$. When $B^2-4AC$ goes positive, the phase II graphon
becomes unstable and the minimizing graphon has $c$ different from $1/2$. 

Since the function $I_0$ is transcendental, it is presumably impossible 
to express the solution to $B^2-4AC=0$ in closed form. However, it is easy
to compute this boundary numerically, as shown in Fig.~\ref{FIG:Phase-Boundary}.

\section{Exponential Random Graph Models}
\label{SEC:ERGM}

There is a large, rapidly growing literature on the analysis of large
graphs, much of it devoted to the analysis of some {\it specific}
large graph of particular interest; see \cite{N} and the references
therein. There is also a significant
fraction of the literature which deals with the asymptotics of graphs
$g$ as the vertex number diverges (see \cite{Lov}), and in these works a popular tool is
`exponential random graph models' (`ERGMS') in which a few densities
$t_j(g),\ j=1,2,\ldots$, such as edge and triangle densities, are
selected, conjugate parameters $\beta_j$ are introduced, and the
family of relative probability densities:
\begin{equation}
\rho_{\beta_1,\beta_2,\ldots}=\exp[{n^2(\beta_1 t_1(g)+\beta_2 t_2(g)+\ldots)}]
\end{equation}
is implemented on the space $G(n)$ of all simple graphs $g$ on $n$
vertices. Densities are assumed normalized to have value 1 in the
complete graph, so the analysis is relevant to `dense graphs', which
have order $n^2$ edges.

$\rho_{\beta_1,\beta_2,\ldots}$ has the obvious form of a grandcanonical 
ensemble, with free energy density
\begin{equation}
\psi_{\beta_1,\beta_2,\ldots}={1\over n^2} \ln\Big[\sum_{g\in G(n)} \exp[{n^2(\beta_1 t_1(g)+\beta_2 t_2(g)+\ldots)}]\Big]
\end{equation}
and suggests the introduction of other ensembles. This paper is
the third in a sequence in which we have studied dense networks in 
the microcanonical
ensemble, in this paper specializing to edge and triangle densities. 

One goal in the study of the asymptotics of graphs is {\it emergent
  phenomena}. The prototypical emergent phenomena are the
thermodynamic phases (and phase transitions) which `emerge' as volume
diverges in the statistical mechanics of particles interacting through
short range forces. The results motivated by these -- the mathematics of phase
transitions in large dense graphs \cite{CD,RY,AR,LZ,RS1,Y,RS2,YRF}
 -- are a significant
instance of emergence within mathematics. We conclude with some
observations on the relationship between the phase transitions in dense
networks as they appear in different ensembles.

In the statistical mechanics of particles with short range forces
there is a well-developed theory of the equivalence of ensembles, and
in particular any phase transition seen in one ensemble automatically
appears in any other ensemble because corresponding free energies are
Legendre transforms of one another and transitions are represented by
singularies in free energies \cite{Ru}. For particles with long range forces
(mean-field models) the equivalence of ensembles can break down
\cite{TET}. For dense networks this is known to be the case \cite{RS1}. So the
precise relationship between the phase decomposition in the various
ensembles of dense networks is an open problem of importance. For instance
in ERGMS one of the key results is a transition within a single phase,
structurally similar to a gas/liquid transition,
whose optimal graphons are all Erd\"os-R\'enyi \cite{CD,RY}. This phenomenon does
not seem to have an obvious image in the microcanonical ensemble of
this paper, but this should be explored.

Another important open problem is the phase structure associated with the
remaining scallops.

\section*{Acknowledgment}
The authors gratefully acknowledge useful discussions with Professor Rick Kenyon (Brown University).
The computational codes involved in this research were developed and
debugged on the computational cluster of the Mathematics Department of
UT Austin. The main computational results were obtained on the
computational facilities in the Texas Super Computing Center
(TACC). We gratefully acknowledge these computational supports.
This work was partially supported by NSF grants DMS-1208941,
DMS-1321018 and DMS-1101326.


\end{document}